\documentclass{article}
\usepackage{cite}
\usepackage{setspace}
\usepackage[colorinlistoftodos,prependcaption]{todonotes}

\usepackage{amssymb}
\usepackage{amsfonts}
\usepackage{amsmath}
\usepackage{amscd}
\usepackage{accents}

\usepackage{bm,bbm}
\usepackage{comment}
\usepackage{rotating}
\usepackage{color}

\usepackage{graphicx}

\usepackage{multirow}
\usepackage{multicol}
\usepackage{numprint}
\usepackage{overpic}

\usepackage{multicol}
\usepackage{float}
\usepackage{mdframed}
\usepackage{outlines}

\newtheorem{theorem}{Theorem}[section]
\newtheorem{lemma}[theorem]{Lemma}

\newtheorem{remark}[theorem]{Remark}

\newtheorem{problem}[theorem]{Problem}

\newcommand{\EOD}{\end{document}}
\newcommand{\TERM}[1]{\mathbf{(#1)}}

\definecolor{MyDarkGreen}{rgb}{0,0.45,0}

\newcommand{\REAL}  {\mathbbm{R}}

\newcommand{\av}{\bm{a}}
\newcommand{\bv}{\bm{b}}

\newcommand{\nv}{\bm{n}}

\newcommand{\qv}{\bm{q}}

\newcommand{\uv}{\bm{u}}
\newcommand{\vv}{\bm{v}}
\newcommand{\wv}{\bm{w}}
\newcommand{\xv}{\bm{x}}

\newcommand{\Bv}{\bm{B}}

\newcommand{\Ev}{\bm{E}}
\newcommand{\Fv}{\bm{F}}

\newcommand{\Jv}{\bm{J}}

\newcommand{\as}{a}
\newcommand{\bs}{b}
\newcommand{\cs}{c}

\newcommand{\gs}{g}

\newcommand{\ns}{n}

\newcommand{\ps}{p}
\newcommand{\qs}{q}

\renewcommand{\ss}{s}
\newcommand{\ts}{t}
\newcommand{\us}{u}
\newcommand{\vs}{v}
\newcommand{\ws}{w}
\newcommand{\xs}{x}
\newcommand{\ys}{y}

\newcommand{\Bs}{B}
\newcommand{\Cs}{C}

\newcommand{\Es}{E}

\newcommand{\Ns}{N}

\newcommand{\Ts}{T}

\newcommand{\Eszr}{\Es_0}
\newcommand{\Bvzr} {\Bv^0}
\newcommand{\Bvhzr}{\Bv^0_{\hh}}

\newcommand{\matA}{\mathbbm{A}}

\newcommand{\calE}{\mathcal{E}}
\newcommand{\calF}{\mathcal{F}}

\newcommand{\calH}{\mathcal{H}}
\newcommand{\calI}{\mathcal{I}}

\newcommand{\calM}{\mathcal{M}}

\newcommand{\calP}{\mathcal{P}}

\newcommand{\calS}{\mathcal{S}}

\newcommand{\calV}{\mathcal{V}}
\newcommand{\calX}{\mathcal{X}}

\newcommand{\bzeta}  {\bm\zeta}
\newcommand{\etav}   {\bm\eta}

\newcommand{\KER}{\textrm{ker}}

\newcommand{\restrict}[2]{{#1}_{|#2}}

\newcommand{\ABS}    [1]{|#1|}
\newcommand{\vABS}   [1]{\left|#1\right|}
\newcommand{\abs}    [1]{\big|#1\big|}
\newcommand{\Abs}    [1]{\Big|#1\Big|}

\newcommand{\scal}   [2]{\big (#1,#2\big)}
\newcommand{\Scal}   [2]{\Big (#1,#2\Big)}

\newcommand{\scalP}  [2]{\big (#1,#2\big)_{\P}}

\newcommand{\SCALVh}  [2]{(#1,#2)_{\Vh}}

\newcommand{\scalVh}  [2]{\big (#1,#2\big)_{\Vh}}

\newcommand{\SCALEh}  [2]{(#1,#2)_{\Eh}}

\newcommand{\scalEh}  [2]{\big (#1,#2\big)_{\Eh}}
\newcommand{\ScalEh}  [2]{\Big (#1,#2\Big)_{\Eh}}
\newcommand{\scallEh} [2]{\bigg(#1,#2\bigg)_{\Eh}}
\newcommand{\ScallEh} [2]{\Bigg(#1,#2\Bigg)_{\Eh}}

\newcommand{\scalVhP} [2]{\big (#1,#2\big)_{\Vh(\P)}}

\newcommand{\scalEhP}  [2]{\big (#1,#2\big)_{\Eh(\P)}}

\newcommand{\nlen}{\hspace{-0.2mm}}

\newcommand{\SNORM}  [2]{\vert#1\vert_{#2}}

\newcommand{\NORM}   [2]{\|#1\|_{#2}}

\newcommand{\norm}   [2]{\big\|#1\big\|_{#2}}

\newcommand{\TNORM}  [2]{\vert\nlen\vert\nlen\vert#1\vert\nlen\vert\nlen\vert_{{}_{#2}}}

\newcommand{\dV}{dV}
\newcommand{\dS}{dS}

\newcommand{\Dt}{\Delta\ts}

\def\trait #1 #2 #3 {\vrule width #1pt height #2pt depth #3pt}
\def\fin{\hfill
        \trait .3 5 0
        \trait 5 .3 0
        \kern-5pt
        \trait 5 5 -4.7
        \trait 0.3 5 0
\medskip}
\newcommand{\ENDPROOF}{\fin}
\newcommand{\BEGINPROOF}{\emph{Proof}.~~}

\usepackage{scalerel,stackengine}
\stackMath
\newcommand\reallywidehat[1]{%
\savestack{\tmpbox}{\stretchto{%
  \scaleto{%
    \scalerel*[\widthof{\ensuremath{#1}}]{\kern-.6pt\bigwedge\kern-.6pt}%
    {\rule[-\textheight/2]{1ex}{\textheight}}
  }{\textheight}%
}{0.5ex}}%
\stackon[1pt]{#1}{\tmpbox}%
}

\newcommand{\DIV} {\textrm{div}\,}
\newcommand{\ROT} {\textrm{rot}\,}

\newcommand{\ROTv}{\textbf{rot}\,}

\newcommand{\pdt}[1]{\frac{\partial#1}{\partial\ts}}

\newcommand{\HONE}  {H^1}
\newcommand{\HONEzr}{H^1_0}

\newcommand{\HTWO}  {H^2}

\newcommand{\LTWO}  {L^2}

\newcommand{\LINF}  {L^{\infty}}

\newcommand{\PS}[1] {\mathbbm{P}_{#1}}

\newcommand{\RT} [1]{\textrm{RT}_{#1}}

\newcommand{\LS}[1] {L^{#1}}
\newcommand{\HS}[1] {H^{#1}}

\newcommand{\CS}[1] {C^{#1}}
\newcommand{\VS}[1] {V^{#1}}
\newcommand{\Vspace}{\VS{}}

\newcommand{\HDIV}   [1]{H  (\textrm{div}; #1)}
\newcommand{\HROT}   [1]{H  (\textrm{rot}; #1)}
\newcommand{\HROTv}  [1]{H  (\textbf{rot}; #1)}

\newcommand{\HROTvzr}[1]{H_0(\textbf{rot}; #1)}

\newcommand{\PinP}{\Pi^{\nabla}_{\P}}

\renewcommand{\P} {\textsf{P}}            
\newcommand  {\E} {\textsf{e}}            
\newcommand  {\V} {\textsf{v}}            
\newcommand  {\T} {\textsf{T}}            

\newcommand{\Vp}  {\V_{1}}
\newcommand{\Vpp} {\V_{2}}

\newcommand{\hh}{h}
\newcommand{\Th}{\Omega_{\hh}}

\newcommand{\xvV}{\xv_{\V}}        

\newcommand{\DIM} {d}              

\newcommand{\hP}{\hh_{\P}}

\newcommand{\hE}{\hh_{\E}}

\newcommand{\mP}{\ABS{\P}}

\newcommand{\mE}{\ABS{\E}}

\newcommand{\Vset}{\mathcal{V}}    

\newcommand{\NMB}{N}
\newcommand{\NV}{\NMB^{\Vset}}   

\newcommand{\NPV}{\NMB^{\Vset}_{\P}}      

\newcommand{\norE} {\mathbf{n}_{\E}}

\newcommand{\tngE} {\mathbf{t}_{\E}}

\newcommand{\qsh}{\qs_{\hh}}

\newcommand{\vsh}{\vs_{\hh}}

\newcommand{\wsh}{\ws_{\hh}}

\newcommand{\Esh}{\Es_{\hh}}

\newcommand{\hEs}{\widehat{\Es}}
\newcommand{\hEsh}{\widehat{\Es}_{\hh}}

\newcommand{\vvh}{\vv_{\hh}}

\newcommand{\wvh}{\wv_{\hh}}

\newcommand{\qvh}{\qv_{\hh}}

\newcommand{\Bvh}{\Bv_{\hh}}

\newcommand{\ash}{\as_{\hh}}

\newcommand{\ashzr}{\as_{\hh,0}}

\newcommand{\bil} [2]{\big<#1,#2\big>}

\newcommand{\Vh}{\calV_{\hh}}
\newcommand{\Eh}{\calE_{\hh}}
\newcommand{\Ph}{\calP_{\hh}}

\newcommand{\Vhzr}{\calV_{\hh,0}}

\newcommand{\IVh}{\calI^{\Vh}}
\newcommand{\IEh}{\calI^{\Eh}}
\newcommand{\IPh}{\calI^{\Ph}}
\newcommand{\IVhP}{\calI^{\Vh}_{\P}}
\newcommand{\IEhP}{\calI^{\Eh}_{\P}}
\newcommand{\IPhP}{\calI^{\Ph}_{\P}}

\newcommand{\PiVhPLS}{\Pi^{LS}_{\P}}
\newcommand{\PiVhPpw}{\Pi^{pw}_{\P}}    
\newcommand{\PiVhP}{\Pi^{\Vh}_{\P}}
\newcommand{\PiVh} {\Pi^{\Vh}}

\newcommand{\PiEhP}{\Pi^{\Eh}_{\P}}

\newcommand{\PiEhRT}{\Pi^{RT}}
\newcommand{\PiEhPRT}{\Pi^{RT}_{\P}}

\newcommand{\Xh}{\calX_{\hh}}

\newcommand{\CurlNorm}[1]{\TNORM{#1}{\Dt,\ROTv}}
\newcommand{\DivNorm} [1]{\TNORM{#1}{\Dt,\DIV}}

\newcommand{\CurlsigmaNorm}[1]{\NORM{#1}{\HS{}_{\sigma}(\ROTv;\Omega)}}

\newcommand{\vshP}{\vs_{\P}}

\newcommand{\HALF}{\frac12}
\newcommand{\Cst}{\widetilde{\Cs}}
\newcommand{\Csb}{\widehat{\Cs}}

\usepackage{authblk}
\usepackage{geometry}
\title {The virtual element method for resistive magnetohydrodynamics.}
\author{S. Naranjo Alvarez$^{a}$, V.A. Bokil$^{b}$, V. Gyrya$^{c}$, G. Manzini$^{d}$}
\date{{\small $^a$Department of Mathematics, Oregon State University, Corvallis, OR, 97331,\\ 
\emph{ e-mail: naranjos@math.oregonstate.edu}}\\
{\small $^b$Department of Mathematics, Oregon State University, Corvallis, OR, 97331,\\
\emph{ e-mail: bokilv@math.oregonstate.edu}}\\
{\small $^c$Group T-5, Theoretical Division, Los Alamos National Laboratory, Los Alamos, 87545 NM, USA,\\\emph{e-mail: vitaliy\_gyrya@lanl.gov}}\\ 
{\small $^d$Group T-5, Theoretical Division, Los Alamos National Laboratory, Los Alamos, 87545 NM, USA,\\\emph{e-mail: gmanzini@lanl.gov}}}
\begin{document}
\maketitle
  

  \begin{abstract}
    We present a virtual element method (VEM) for the numerical approximation of 
    the electromagnetics subsystem of the resistive magnetohydrodynamics (MHD) model in two spatial dimensions.
    The major advantages of the virtual element method include great flexibility of polygonal meshes 
    and automatic divergence-free constraint on the magnetic flux field.
    In this work, 
    we rigorously prove the well-posedness of the method
    and the solenoidal nature of the discrete magnetic flux field.
    We also derive stability energy estimates.
    The design of the method includes three choices for the construction of the nodal mass matrix and criteria to more alternatives.
    We present a set of numerical experiments that independently validate theoretical results.
    The numerical experiments include the convergence rate study,
    energy estimates 
    and 
    verification of the divergence-free condition on the magnetic flux field.
    All these numerical experiments have been performed on triangular, perturbed quadrilateral and Voronoi meshes.
    Finally, We demonstrate the development of the VEM method on a numerical model for the Hartmann flow.
  \end{abstract}

  \begin{keywords}
    Maxwell equations, 
    resistive MHD,
    virtual element method,
    polytopal mesh,
    energy stability analysis,
    Hartmann flow.
  \end{keywords}

\section{Introduction}
Interest in the behavior of plasmas has skyrocketed in the modern age
with applications ranging from fusion-based nuclear power to low power
thrusters for contemporary spacecraft. 
Since the late 1930s, efforts have been devoted to the development of
models for plasmas and discretizations that are faithful to the
physics and dynamics.
An approach that has proven successful and has become standard is to
consider plasmas as magnetized fluids, an area called
magnetohydrodynamics (MHD).  
Therefore, the description of these plasmas follow from a blending
together of electromagnetic theory and fluid flow. 
The precise details of how these two theories can be coupled can be
found in \cite{davidson2002introduction,schindler2006physics,moreau2013magnetohydrodynamics}.
Research in MHD is driven by applications that are important to several communities including, astrophysicists that study accretion discs and the dynamics that
govern evolution of stars; planetary scientists that are interested in
the generation of magnetic fields at the core of planets;
plasma physicists whose interest lies in the confinement of plasmas by
means of external magnetic fields and engineers who have found that
with external magnetic fields they can control the motion of liquid
metals leading to a revolution in metallurgical techniques in
industry.

The development of numerical methods for MHD is an active area of
research, being developed over the last few decades.
In \cite{hu2017stable,hiptmair2018fully}, two finite element methods
are presented that use different techniques in order to preserve the
divergence condition on the magnetic field.
In \cite{hu2017stable}, the condition is attained automatically,
similar to how it is done in this article, whereas in
\cite{hiptmair2018fully} the scheme includes the magnetic vector
potential under the temporal gauge, and the magnetic field is obtained
as its curl.
In \cite{torrilhon2003non}, the convergence of finite volume methods
for MHD is studied and in
\cite{powell1999solution,crockett2005unsplit} the classic upwind and
Godunov methods are adapted to ideal MHD.
In \cite{corti2012stable}, the author presents a finite difference
method based on summation by parts (SBP) to
mimic the integration by parts formula in the discrete setting, in
order to preserve important energy conservation properties and attain
an approximate-divergence free scheme.
Finally, in \cite{liu2001energy} the authors develop a MAC scheme for
the fluid flow sub-system of the incompressible MHD equations,
coupling it to the Yee-scheme for the electromagnetic sub-system.

Although models in MHD come about from a coupling between the
equations that govern the fluid flow and Maxwell's equations for
electromagnetism, in this article we will focus on modeling the
evolution of the electric and magnetic fields in a plasma for a
prescribed fluid flow.
Thus, we focus on the Maxwell subsystem of MHD, which combines
Faraday's law, Ampere's law, Ohm's law and Gauss's law for the
electromagnetic fields under a prescribed fluid flow.

The main aim of this article is to present a novel numerical
discretization of Maxwell's equations for resistive MHD in a two
dimensional setting using the virtual element method (VEM).
The VEM was originally proposed
in~\cite{BeiraodaVeiga-Brezzi-Cangiani-Manzini-Marini-Russo:2013} as a
variational reformulation of the \emph{nodal} high-order mimetic
finite difference (MFD) method~\cite{%
  BeiraodaVeiga-Manzini-Putti:2015,%
  BeiraodaVeiga-Lipnikov-Manzini:2011,%
  Brezzi-Buffa-Lipnikov:2009}, for the numerical treatment of elliptic
problems on unstructured polygonal and polyhedral meshes.
The word ``mimetic'' reflects the nature of the method, which mimics
the duality and self-adjointness of differential operators as well as
identities of vector and tensor calculus.
Due to such feature, mimetic methods are often dubbed as
\emph{compatible methods} or \emph{compatible discretizations}.
In particular, satisfying Gauss's law on the divergence of the
magnetic field in the discrete setting requires careful discretization
of the Maxwell curl equations, i.e., Faraday's law and Maxwell-Amp\`{e}re
law.
This fact is in contrast to the continuous setting in which the
divergence-free nature of the magnetic field is a direct consequence of the Maxwell
curl equations when the initial conditions properly satisfy the
Gauss's law.

The violation of Gauss's law is a serious
source of error in the numerical discretization of Maxwell's equations,
causing the appearance of fictitious forces or magnetic monopoles,
which are non-physical, thus rendering the numerical
simulations unfaithful to the real physics.
Over time, mimetic methods were extended from the Support Operator
Method
(SOM)~\cite{Shashkov:1995,Shashkov-Steinberg:1995,Hyman-Shashkov:1997}),
which works on regular tensor grids, to the MFD method, which works
on fairly general polygonal and polyhedral meshes.
The MFD method is, in practice, a family of schemes depending on a set
of parameters.
These parameters can be optimized to satisfy additional properties
such as maximum principles and low dispersion errors.
This process goes by the name of mimetic adaptation or M-adaptation
and it is outlined in \cite{Gyrya-Lipnikov-Manzini:2016}.
Previous work in M-adaptation shows that the process can be
implemented for problems in wave propagation, see
\cite{bokil2015dispersion}, and in the study of cold plasmas, as shown
in \cite{bokil2016dispersion}.
Readers interested in historical perspective on the 50-year long development of mimetic and compatible methods are referred to
\cite{Lipnikov-Manzini-Shashkov:2014}.
Development of mimetic and compatible methods are referred to a recent review in~\cite{Lipnikov-Manzini-Shashkov:2014}.

The VEM can also be interpreted as a generalization of the FEM to
general polygonal and polyhedral meshes that inherits the great
flexibility of the MFD method with respect to the admissible meshes
used in the numerical formulation.
The major difference when compared to a regular FEM is that in the VEM
the shape functions are defined in an implicit manner and never
explicitly constructed.
The name ``virtual element'' stems from the fact that such shape
functions and the finite element space generated by their linear
combinations are, in this sense, ``virtual''.

The VEM was
originally proposed for solving diffusion problems
in~\cite{BeiraodaVeiga-Brezzi-Cangiani-Manzini-Marini-Russo:2013} as a
conforming FEM, and later extended to the nonconforming formulation
in~\cite{AyusodeDios-Lipnikov-Manzini:2016} and the mixed BDM-like and
RT-like formulations in~\cite{Brezzi-Falk-Marini:2014}
and~\cite{BeiraodaVeiga-Brezzi-Marini-Russo:2016c}, respectively.
Generalizations to
convection-reaction-diffusion problems with variable coefficients
can be found 
in~\cite{Ahmad-Alsaedi-Brezzi-Marini-Russo:2013,BeiraodaVeiga-Brezzi-Marini-Russo:2016b,Cangiani-Manzini-Sutton:2017,Berrone-Borio-Manzini:2018:CMAME:journal}.
In a series of
papers~\cite{BeiraodaVeiga-Brezzi-Marini-Russo:2016a,BeiraodaVeiga-Brezzi-Dassi-Marini-Russo:2017,BeiraodaVeiga-Brezzi-Dassi-Marini-Russo:2018a,BeiraodaVeiga-Brezzi-Dassi-Marini-Russo:2018b},
$H(div)$- and $H(curl)$-conforming virtual element spaces on general
polygonal and polyhedral elements have been proposed to
generalize the well known Raviart-Thomas and N\'{e}d\'{e}lec finite
elements to unstructured polytopal meshes.
These methods, combined with the serendipity strategy that reduces the
total number of degrees of freedom,
see~\cite{BeiraodaVeiga-Brezzi-Marini-Russo:2016d,BeiraodaVeiga-Brezzi-Marini-Russo:2017},
have successfully been applied to the numerical resolution of the
magnetostatic Kikuchi's model.
In these papers, exact virtual de
Rham sequences with commuting-diagram interpolation operators are
built and the solenoidal nature of the discrete magnetic flux field is ensured. Finally, VEMs have also been designed for hyperbolic problems (see ~\cite{vacca2017virtual,adak2020virtual}).

In our work, we utilize the low-order spaces proposed
in~\cite{BeiraodaVeiga-Brezzi-Marini-Russo:2016a}, which makes it
possible to obtain the combined approximation of the $H^1$-conforming space
(0-forms) by a nodal-type virtual element space, the
$H(curl)$-conforming space (1-forms) by an edge-type virtual element
space, and the $H(div)$-conforming space (2-forms)
by discontinuous piecewise constant polynomials.

To derive our virtual element approximation, we first reformulate the
MHD equations in a variational framework, and then, approximate all
$\LTWO$-type integrals by using suitably defined inner products for
nodal-, edge- and cell-type virtual functions.
The standard way to build such inner products is through the
orthogonal projection of the virtual element functions onto the
subspace of linear polynomials.
However, as was already noted
in~\cite{BeiraodaVeiga-Brezzi-Marini-Russo:2016a}, the nodal virtual
element space that we consider in this work does not provide enough
information to construct such projections.
Our approach in this paper is to substitute the orthogonal projector
with the elliptic
projector in~\cite{Ahmad-Alsaedi-Brezzi-Marini-Russo:2013}, since in the
low-order case we can always consider these two projection operators as
equal by redefining the virtual element space appropriately.
This strategy is usually referred to as the ``enhanced VEM'' by the
VEM developers and practitioners.

A major issue occurs here, because changing the definition of the nodal virtual
element space requires also changing the definition of the edge and
cell virtual element spaces in order to maintain the exact de Rham commuting diagrams.
This issue has led to the different virtual element space formulations that were used in the magnetostatics application mentioned above.
Instead, in this work we prefer to adopt a different approach, which
consists in designing a special reconstruction operator that is
computable from the degrees of freedom and 
stable and bounded as discussed in the following sections.
Applying the reconstruction operator makes it possible to
recover an approximation of the nodal virtual element functions 
inside each mesh element and then we integrate directly these reconstructed functions.
The choice of the elemental reconstruction operator is not unique. 
In this work, we considered three different options: the elliptic projection; a Least-Squares interpolation of the nodal values; and the piecewise linear Galerkin interpolation on a triangular sub partition of each element. 
Our numerical experiments show that these three options are all quite effective and the resulting scheme's implementations have comparable accuracy. 

\medskip
This article is structured as follows.
The rest of this section includes a brief overview of notation and
some basic mathematical definitions relevant to the rest of the paper.
In Section~\ref{sec:Mathematical:Formulation}, we present the set of
governing equations to be discretized in the continuous setting and
introduce the semi-discrete and fully discrete variational
formulations in the virtual element framework.
Next, in Section~\ref{sec:virtual:element:method}, we define the
virtual element spaces and detail the construction of the inner
products that are used for the numerical approximation of the MHD
model equations.
We also discuss the exactness and commutativity properties of the De-Rham complex and prove that the divergence free condition of the
numerical approximation of the magnetic flux field is preserved over
time.
In Section~\ref{sec:wellposedness:VEM}, we prove that the fully
discrete variational formulation is well posed.
In Section~\ref{sec:stability:energy:estimate}, we derive stability
energy estimates for the continuous and fully discrete models.
In Section~\ref{Sec:NumericalExperiments}, we present the results of a
series of numerical experiments that provide evidence regarding the
convergence rate of the numerical method.
Plots demonstrating that the method preserves the divergence free
condition of the magnetic flux field are available as well as a numerical
study of the energy estimates that are derived theoretically in
Section~\ref{sec:stability:energy:estimate}.
We conclude this section by presenting a simulation of the solution to
the Hartmann Flow problem.
Then, finally we summarize our findings in
section~\ref{sec:conclusions}.
\subsection{Notation, functional spaces and technicalities}
We use the standard definition and notation of Sobolev spaces, norms
and seminorms, cf.~\cite{Adams-Fournier:2003}.
Let $k$ be a non-negative integer.
Consider an open bounded connected subset $\omega$ of $\REAL^{2}$ with
polygonal boundary $\partial\omega$.
Subset $\omega$ can be the whole computational domain $\Omega$, or one
of the polygonal cells $\P$ of the mesh partitioning, $\Th$, covering
$\Omega$.

The Sobolev space $\HS{k}(\omega)$ consists of all square integrable
functions with all square integrable weak derivatives up to order $k$
that are defined on $\omega$.
As usual, if $k=0$, we prefer the notation $\LTWO(\omega)$.
Norm and seminorm in $\HS{k}(\omega)$ are denoted by
$\NORM{\cdot}{k,\omega}$ and $\SNORM{\cdot}{k,\omega}$, respectively.
We denote the inner product in $\LTWO(\omega)$ by
$(\cdot,\cdot)_{\omega}$, but we omit the subscript when $\omega$ is
the whole computational domain $\Omega$.
We denote the norm of an operator $\Pi$, which is a norm in the dual
space, by the general notation $\NORM{\Pi}{}$, regardless of the
spaces where range and image of $\Pi$ are defined.

On $\omega$, we consider the functional spaces:
\begin{subequations}
  \label{eq:spaces}
  \begin{align}
    \LTWO (\omega)  & := \left\{\vs:\omega\to\REAL:\int_{\omega}\ABS{\vs}^2\dV<\infty\right\},                \\[0.5em]
    \HROTv{\omega}  & := \left\{\vs\in\LTWO(\omega):\ROTv\vs\in\left(\LTWO(\omega)\right)^2\right\},     \\[0.5em]
    \HDIV{\omega}   & := \left\{\wv\in\left(\LTWO(\omega)\right)^2:\DIV\wv\in\LTWO(\omega)\right\},      \\[0.5em]
    \LINF(\omega)   & := \{\wv:\omega\to\REAL:\exists\,\Cs>0\,;\ABS{\wv}<\Cs\;\mbox{almost everywhere}\}, 
  \end{align}
\end{subequations}
where
$\ROTv\vs=(\partial\vs\slash{\partial\ys},-\partial\vs\slash{\partial\xs})^T$,
and
$\ROT\wv=(\partial\ws_x\slash{\partial\ys}-\partial\ws_y\slash{\partial\xs})$
for the vector field $\wv=(\ws_x,\ws_y)^T$.
If $\omega=\Omega$ denotes the computational domain, we consider the
functional spaces:
\begin{subequations}
  \label{eq:spaces:global}
  \begin{align}
    \Vspace          & := \left\{\wv\in\HDIV{\Omega}\,:\,\wv\in\left(\LS{2+s}(\Omega)\right)^2, \ \textrm{for~some~real~}\,s>0\right\},\\[0.5em]
    \HROTvzr{\Omega} & := \left\{\vs\in\HROTv{\Omega}\,:\,\vs=0\mbox{~on~}\partial\Omega\right\}.
  \end{align}
\end{subequations}
Space $\Vspace$ is slightly more regular than $\HDIV{\Omega}$ to
ensure that the trace of the normal component
$\restrict{\vvh\cdot\nv}{\E}$ on each mesh edge $\E$ exists and is
continuous across all the internal
edges~\cite{Boffi-Brezzi-Fortin:2013}.

For an open bounded connected subset $\omega\subset\REAL^d$ with $d=1$
or $2$, we denote the linear space of polynomials of degree up to
$\ell$ defined on $\omega$ by $\PS{\ell}(\omega)$, with the useful
conventional notation that $\PS{-1}(\omega)=\{0\}$.
We denote the space of two-dimensional vector polynomials of degree up
to $\ell$ on $\omega$ by $\big[\PS{\ell}(\omega)\big]^2$.
Space $\PS{\ell}(\omega)$ is the span of the finite set of
\emph{scaled monomials of degree up to $\ell$}, that are given by
\begin{align*}
  \calM_{\ell}(\omega) =
  \bigg\{\,
    \left( \frac{\xv-\xv_{\omega}}{\hh_{\omega}} \right)^{\alpha}
    \textrm{~with~}\ABS{\alpha}\leq\ell
    \,\bigg\},
\end{align*}
where 
\begin{itemize}
\item $\xv_{\omega}$ denotes the center of gravity of $\omega$ and
  $\hh_{\omega}$ its characteristic length, as, for instance, the edge
  length or the cell diameter for $\DIM=1,2$;
\item $\alpha=(\alpha_1,\alpha_2)$ is the
  two-dimensional multi-index of nonnegative integers $\alpha_i$
  with degree $\ABS{\alpha}=\alpha_1+\alpha_{2}\leq\ell$ and
  such that
  $\xv^{\alpha}=\xs_1^{\alpha_1}\xs_{2}^{\alpha_{2}}$ for
  any $\xv\in\REAL^{2}$.
\end{itemize}
We will also use the set of \emph{scaled monomials of degree exactly
  equal to $\ell$}, denoted by $\calM_{\ell}^{*}(\omega)$ and obtained
by setting $\ABS{\alpha}=\ell$ in the definition above.

Finally, we use the letter $\Cs$ in many inequalities to denote a
strictly positive constant whose value can change at any instance.
Constant $\Cs$ may depend on the constants of the model equations or
the variational problem, like the coercivity and continuity constants,
or constants that are uniformly defined for the family of meshes of
the approximation while $\hh\to0$, such as the mesh regularity
constant, the stability constants of the discrete bilinear forms, etc.
However, constant $\Cs$ will never depend on the discretization
parameters such as the mesh size $\hh$ and the timestep $\Dt$.

\setcounter{section}{1}
\section{The mathematical formulation}
\label{sec:Mathematical:Formulation} 
Let $\Omega$ be an open, bounded, and polygonal subset of $\REAL^2$
with boundary $\Gamma=\partial\Omega$ and $\Ts$ a positive real
number.
For a given fluid flow described by the velocity vector field
$\uv=(\us_x,\us_y)^T\in\big[\LINF(\Omega)\big]^2$, we consider the
Maxwell problem for the electric and magnetic fields,
respectively denoted by $\Es$ and $\Bv=(\Bs_x,\Bs_y)^T$, that reads as:
\begin{subequations} 
  \label{eq:IntroSystem}
  \begin{align}
    \pdt{\Bv}         &= -\ROTv\Es\phantom{\nu\;\ROT\Bv\Bvzr\Es_0(\cdot,t)\quad\textrm{with}\,\,\,\DIV\Bvzr=0}\hspace{-2cm}\mbox{in~}\Omega\times (0,\Ts],       \label{eq:IntroSystemFaraday}     \\[0.5em]
    \Es+\uv\times\Bv  &= \nu\;\ROT\Bv\phantom{-\ROTv\Es\Bvzr\Es_0(\cdot,t)\quad\textrm{with}\,\,\,\DIV\Bvzr=0}\hspace{-2cm}\mbox{in~}\Omega\times (0,\Ts],       \label{eq:IntroSystemAmpereOhm}   \\[0.5em]
    \Bv(\cdot,0)      &= \Bvzr\quad\textrm{with}\,\,\,\DIV\Bvzr=0\phantom{-\ROTv\Es\nu\;\ROT\Bv\Es_0(\cdot,t)}\hspace{-2cm}\mbox{in~}\Omega,                     \label{eq:IntroSystemBoundaryCond}\\[0.5em]
    \Es(\cdot,t)      &= \Es_0(\cdot,t)\phantom{-\ROTv\Es\nu\;\ROT\Bv\Bvzr\quad\textrm{with}\,\,\,\DIV\Bvzr=0}\hspace{-2cm}\mbox{on~}\partial\Omega\times(0,\Ts],
  \end{align}
\end{subequations}
where $\nu$ is the resistivity of the medium
and $\uv\times\Bv=\us_{x}\Bs_{y}-\us_y\Bs_x$.
We introduce $\sigma=\nu^{-1}$ and assume that it is bounded by two
positive constants $\sigma_*\leq\sigma(\xv)\leq\sigma^*$ for almost
every $\xv\in\Omega$.
The system of partial differential equations
(PDEs)~\eqref{eq:IntroSystem} couples Faraday, Ampere and Ohm
laws.
As discussed in the introduction, an important property of the MHD
system \eqref{eq:IntroSystem}, which we will address in the virtual
element discretization, is the solenoidal nature of the magnetic flux
field $\Bv$.
By taking the divergence of \eqref{eq:IntroSystemFaraday} we find that
the divergence of $\Bv$ does not change in time, so $\Bv$ is
divergence free if the initial field $\Bvzr$
in~\eqref{eq:IntroSystemBoundaryCond} is divergence free.

\medskip
The variational formulation of problem~\eqref{eq:IntroSystem} reads
as:

\medskip
\emph{Find
  $(\Bv,\hEs)\in\CS{1}\left([0,\Ts],\HDIV{\Omega}\right)\times\CS{}\left([0,\Ts],\HROTvzr{\Omega}\right)$,
  such that:}
\begin{subequations}\label{eq:VarForEq2}
  \vspace{-0.2em}
  \begin{align}
    &\Scal{\pdt{\Bv}}{\wv} + \Scal{\ROTv\Es}{\wv} = 0
    \qquad\forall\wv\in\HDIV{\Omega},
    \label{eq:VarForEq2Faraday}
    \\[0.75em]
    &\Scal{\sigma\hEs}{\vs}+\Scal{\sigma\uv\times\Bv}{\vs}-\Scal{\Bv}{\ROTv\vs} = -\Scal{\sigma\Es_0}{\vs}
    \qquad\forall\vs\in\HROTvzr{\Omega},
    \label{eq:VarForEq2AmpereOhm}
    \\[0.50em]
    &\Es = \hEs + \Es_0,
    \label{eq:electric:field:update}
    \\[0.75em]
    &\Bv(\cdot,0) = \Bvzr\;\;\mbox{with}\;\;\DIV\Bvzr = 0.
    \label{eq:magnetic:field:update}
  \end{align}
\end{subequations}
The boundary conditions on $\Es$ are set through the known function
$\Es_0$, so we seek for the solution $\hEs$ with zero trace on
$\Gamma$.
As is the case of any conforming Galerkin method, 
we first select subspaces of $\HDIV{\Omega}$ and $\HROTvzr{\Omega}$ defined 
on the mesh partition $\Th$ of the computational domain $\Omega$. 
The requirements on the mesh partition $\Th$ will be specified in
Section~\ref{sec:virtual:element:method}.
We respectively denote them by $\Eh$ and $\Vh$, and assume that they
are equipped by the inner products $\SCALEh{\cdot}{\cdot}$ and
$\SCALVh{\cdot}{\cdot}$ and suitable interpolation operators, e.g.,
$\IVh$ and $\IEh$, or projection operator, e.g., $\PiEhRT$.
The coefficient $\sigma$ is incorporated in the definition of $\Vh$.
We also use the space $\Vhzr$, the subspace of the functions in $\Vh$
vanishing at the boundary of $\Th$.
The definition and construction of all these mathematical objects are
left for the next section.
The semi-discrete virtual element discretization of
Problem~\eqref{eq:VarForEq2} reads as:

\medskip
\emph{Find
  $(\Bvh,\hEsh)\in\CS{1}\left([0,\Ts],\Eh\right)\times\CS{}\left([0,\Ts],\Vhzr\right)$
  such that for all $(\wvh,\vsh)\in\Eh\times\Vhzr$ it holds:}
\begin{subequations}
  \label{eq:VarForEq3}
  \vspace{-0.2em}
  \begin{align}
    &\ScallEh{\pdt{\Bvh}}{\wvh} + \scalEh{\ROTv\Esh}{\wvh} = 0,
    \label{eq:VarForEq3Faraday}
    \\[0.75em]
    &\scalVh{\hEsh}{\vsh}+\scalVh{\IVh\big(\uv\times\PiEhRT\Bvh\big)}{\vsh}-\scalEh{\Bvh}{\ROTv\vsh} = -\scalVh{\IVh(\Es_0)}{\vsh},
    \label{eq:VarForEq3AmpereOhm}
    \\[0.75em]
    &\Esh = \hEsh + \IVh(\Es_0),
    \label{eq:electric:field:update:VEM}
    \\[0.75em]
    &\Bvh(\cdot,0) = \Bvhzr = \IEh(\Bvzr)\;\;\mbox{with}\;\;\DIV\Bvhzr = 0.
    \label{eq:magnetic:field:update:VEM}
  \end{align}
\end{subequations}

\medskip
Let $\Dt$ denote the timestep that splits the time interval $[0,T]$
into $N=T\slash{\Dt}$ subintervals.
The virtual element solution pair
$\big(\Bvh(\cdot,t^n),\hEsh(\cdot,t^n+\theta\Dt)\big)$, with $t^n = n\Delta t$,  is approximated by
the pair $(\Bvh^n,\hEsh^{n+\theta})$, which is the solution of the
discrete time-dependent problem parameterized by the scalar factor
$\theta\in [0,1]$:

\medskip
\emph{ Find $\big\{\Bvh^n\big\}_{n=0}^N\subset\Eh$ and
  $\big\{\hEsh^{n+\theta}\big\}_{n=0}^{N-1}\subset\Vhzr$ such that
  for all $(\wvh,\vsh)\in\Eh\times\Vhzr$ it holds:}
\begin{subequations} 
  \label{eq:VarForEq4}
  \vspace{-0.2em}
  \begin{align}
    &\ScalEh{\frac{\Bvh^{n+1}-\Bvh^n}{\Dt}}{\wvh}+
    \ScalEh{\ROTv\Esh^{n+\theta}}{\wvh}=0
    \label{eq:VarForEq4Faraday}
    \\
    &
    \scalVh{\hEsh^{n+\theta}}{\vsh}
    +\scalVh{\IVh\big(\uv\times\PiEhRT\Bvh^{n+\theta}\big)}{\vsh}-
    \scalEh{\Bvh^{n+\theta}}{\ROTv\vsh}
    =
    -\scalVh{\IVh\big(\Es_0^{n+\theta}\big)}{\vsh}
    \label{eq:VarForEq4AmpereOhm}
    \\[0.75em]
    &\Esh^{n+\theta} = \hEsh^{n+\theta} + \IVh\big(\Eszr^{n+\theta}\big),\\[0.75em]
    &\Bvh^{n+\theta}=\theta\Bvh^{n+1}+(1-\theta)\Bvh^n,                \\[0.75em]
    &\Bvh(\cdot,0) = \IEh\big(\Bvzr\big).
  \end{align}
\end{subequations}
It is worth noting that $\Bvh^n$ is defined at the time instants
$t^n=n\Dt$, $n=0,\ldots,N$, while $\hEsh^{n+\theta}$ is defined on the
``staggered'' grid at time instants $t^{n+\theta}=(n+\theta)\Dt$.
According to the $\theta$ parameterization, for $\theta=0$ we recover
the explicit or forward Euler scheme, for $\theta=1$ the implicit or
backward Euler scheme and for $\theta=1/2$ the (semi) implicit 
leap-frog scheme.
\section{The virtual element method}
\label{sec:virtual:element:method}





\subsection{Assumptions on mesh regularity }
\label{subsec:mesh-regularity}
Let $\Th=\{\P\}$ be a mesh decomposition of $\Omega$ into polygonal
element (or cell) $\P$ with boundary $\partial\P$, area $\mP$ and diameter
$\hP$.
As usual, $\hh=\max_{\P\in\hh}\hP$ is the \textit{mesh size} parameter.
We denote the edges of $\partial\P$ by $\E$ and its length by
$\mE=\hE$.

We assume that $\hh$ belongs to $\calH\subset(0,+\infty)$, which is a
countable set of mesh sizes having $0$ as its unique accumulation
point.
A family of meshes $\{\Th\}_{\hh}$ is said to be \emph{regular} if
there exists a non-negative real number $\rho$ independent of $\hh$
(and, hence, of $\Th$), such that

\medskip
\begin{description}
\item\textbf{(M1)} (\emph{star-shapedness}): every polygonal cell $\P$
  of every mesh $\Th$ is star-shaped with respect to every point of a
  disk of radius $\rho\hP$;

  \medskip
\item\textbf{(M2)} (\emph{uniform scaling}): every edge
  $\E\in\partial\P$ of cell $\P\in\Th$ satisfies $\hE\geq\rho\hP$.
\end{description}

\medskip
The regularity assumptions \textbf{(M1)}-\textbf{(M2)} allow us to
use meshes with cells having quite general geometric shapes.
For example, nonconvex cells or cells with hanging nodes on their
edges are admissible.
Nonetheless, these assumptions have some important implications such
as: $(\mathrm{i})$ every polygonal element is \textit{simply connected}; $(\mathrm{ii})$
the number of edges of each polygonal cell in the mesh family
$\{\Th\}_{\hh}$ is uniformly bounded; $(\mathrm{iii})$ a polygonal element
cannot have \emph{arbitrarily small} edges with respect to its
diameter $\hP\leq\hh$ for $\hh\to0$ and inequality
$\hP^2\leq\Cs(\rho)\mP\hP^2$ holds, with the obvious dependence of
constant $\Cs(\rho)$ on the mesh regularity factor $\rho$.
It is worth mentioning that virtual element methods on polygonal or
polyhedral meshes possibly containing ``small edges'' in 2D or ``small
faces'' in 3D have been considered in~\cite{Brenner-Sung:2018} for the
numerical approximation of the Poisson problem.
The work in~\cite{Brenner-Sung:2018} extends the results
in~\cite{BeiraodaVeiga-Lovadina-Russo:2017}
for the original two-dimensional virtual element method to the version
of the virtual element method in
\cite{Ahmad-Alsaedi-Brezzi-Marini-Russo:2013}
that can also be applied to problems in three dimensions.

Finally, we note that assumptions $(\mathrm{i})$-$(\mathrm{iii})$ above also imply that
the classical polynomial approximation theory in Sobolev spaces
holds~\cite{Brenner-Scott:2008}.

\subsection{Nodal virtual element space}

On every element $\P\in\Th$, we consider the local virtual element
space:
\begin{align}
  \Vh(\P) := \Big\{\,\vsh\in\HONE(\P)\,:\,
  \restrict{\vsh}{\partial\P}\in\CS{0}(\partial\P),\,
  \restrict{\vsh}{\E}\in\PS{1}(\E),\,\forall\E\in\partial\P,\,
  \Delta\vsh=0~\textrm{in~}\P
  \,\Big\}.
  \label{eq:VhP:def}
\end{align}
Then, we define the global virtual element space:
\begin{align}
  \Vh := \Big\{\,\vsh\in\HONE(\Omega)\,:\,
  \restrict{\vsh}{\P}\in\Vh(\P),\,\,
  \forall\P\in\Th
  \,\Big\};
  \label{eq:Vh:def}
\end{align}
The local and global spaces $\Vh(\P)$ and $\Vh$ were first proposed in
~\cite{BeiraodaVeiga-Brezzi-Cangiani-Manzini-Marini-Russo:2013}.
Space $\Vh$ is a subspace of $\HONE(\Omega)$, so every virtual element
function $\vsh\in\Vh(\P)$ is continuous over the computational domain
$\Omega$.
Every function $\vsh\in\Vh(\P)$ is uniquely determined by its values at
the vertices of $\P$, i.e., by the set
$\{\vsh(\xvV)\}_{\V\in\partial\P}$.
Similarly, a virtual element function in the global space $\Vh$ is
defined by its values at all the mesh vertices.
The unisolvence of such degrees of freedom is proved
in~\cite{BeiraodaVeiga-Brezzi-Cangiani-Manzini-Marini-Russo:2013}.

The virtual element schemes~\eqref{eq:VarForEq3}
and~\eqref{eq:VarForEq4} require an approximation of the $\LTWO$-inner
product in $\Vh$.
The usual approach to build such an approximation would be through the
local orthogonal projection onto the space of linear polynomials,
which is a subspace of $\Vh(\P)$, and by adding a suitable
stabilization term.
However, the orthogonal projection is not computable from the degrees
of freedom of the virtual element functions, namely, the vertex
values, unless we change the definition of the elemental space
according to the construction proposed
in~\cite{Ahmad-Alsaedi-Brezzi-Marini-Russo:2013}.
Here, we prefer not to modify the definition of space $\Vh(\P)$ since
otherwise we would lose the property that $\Vh$ is in a de Rham
complex with space $\Eh$ (which will be defined in the next
subsection).
This topic will be discussed in Section~\ref{subsec:VirtualdeRham}.

Therefore, for the construction of the approximate $\LTWO$-inner
product in $\Vh(\P)$, we proceed in two steps.
First, for each function $\vsh\in\Vh(\P)$, we introduce the linear
polynomial approximation $\PiVhP\vsh$, where operator
$\PiVhP:\HONE(\P)\to\PS{1}(\P)$ has these properties:
\begin{itemize}
\item[] \textbf{(V1)} the linear polynomial $\PiVhP\vsh$ is
  \emph{computable from the degrees of freedom} of $\vsh$;
\item[] \textbf{(V2)} operator $\PiVhP$ is \emph{invariant on linear
    polynomials}, i.e., $\PiVhP\qs=\qs$ whenever $\qs$ belongs to
  $\PS{1}(\P)$;
\item[] \textbf{(V3)} operator $\PiVhP$ is \emph{ uniformly bounded}
  independently of the characteristics of the polygonal element $P$, i.e. 
  there exists a real constant $C>0$ independent of the number of nodes, edges or diameter
  of $P$ such that for every $\vsh\in\Vh(\P)$ one has
  $\NORM{\PiVhP\vsh}{0,\P}\leq C\,\NORM{\vsh}{0,\P}$.
\end{itemize}

\medskip
\begin{remark}\label{remark:codomain}
  In Section~\ref{sec:subsub:construction:oblique:projection}, the
  third operator $\PiVhP\vsh$ is defined as the Galerkin piecewise
  linear interpolant of $\vsh$ on a triangle subpartition of $\P$.
  Such subpartition, which we denote by $\P_{\hh}$, is built by
  connecting the barycenter of $\P$ with its vertices.
  Therefore, conditions \textbf{(V1)}-\textbf{(V3)} above are set for
  $\PiVhP:\HONE(\P)\to\PS{1}(\P_{\hh})$, where $\PS{1}(\P_{\hh})$ is
  the space of continuous piecewise linear polynomials defined on
  $\P_{\hh}$.
\end{remark}

\medskip
In view of mesh regularity assumptions \textbf{(M1)}-\textbf{(M2)} and
according to a Bramble-Hilbert argument~\cite{Bramble-Hilbert:1970,Dupont-Scott:1980} and property \textbf{(V2)},
the \emph{approximation error} satisfies the upper bound estimate
\begin{align}
  \NORM{(1-\PiVhP)\vsh}{0,\P}\leq C\hP\SNORM{\vsh}{1,\P},
  \label{eq:upper-bound:hP}
\end{align}
for every function $\vsh\in\Vh(\P)\subset\HONE(\P)$ and
\begin{align}
  \NORM{(1-\PiVhP)\vsh}{0,\P}+\hP\SNORM{(1-\PiVhP)\vsh}{1,\P} \leq
  C\hP^2\SNORM{\vsh}{2,\P}
  \label{eq:upper-bound:hP2}
\end{align}
whenever $\vsh\in\Vh(\P)\cap\HTWO(\P)$. 



Second, we consider the bilinear form on $\Vh\times\Vh$ given by the
formula
\begin{align}
  \scalVh{\vsh}{\wsh} := \sum_{\P\in\Th}\scalVhP{\vsh}{\wsh}
  \qquad\forall\vsh,\wsh\in\Vh,
  \label{eq:Vh:inner:product:def}
\end{align}
where each local term $\scalVhP{\vsh}{\wsh}$ is computed by using the
elementwise approximations of $\vsh$ and $\wsh$ on $\P$ according to
\begin{align}
  \scalVhP{\vsh}{\wsh} 
  := \scalP{\sigma\PiVhP\vsh}{\PiVhP\wsh}
  +  \calS^{\Vh}_{\P}\big( (1-\PiVhP)\vsh, (1-\PiVhP)\wsh \big).
  \label{eq:Vh:local:inner:product:def}
\end{align}
Here, $\calS^{\Vh}_{\P}(\cdot,\cdot)$ is a symmetric and nonnegative
bilinear form for which there exist two positive constant $\ss_*$ and
$\ss^*$ such that
\begin{align}
  \ss_*\NORM{\vsh}{0,\P}^2\leq\calS^{\Vh}_{\P}(\vsh,\vsh)\leq\ss^*\NORM{\vsh}{0,\P}^2
  \qquad\forall\vsh\in\Vh(\P)\cap\KER\big(\PiVhP{}\big).
  \label{eq:SP:stability}
\end{align}
Constants $\ss_*$ and $\ss^*$ are independent of $\hh$, but may depend
on the regularity parameter $\rho$ and the bounds on $\sigma$, namely,
the two constant factors $\sigma_*$ and $\sigma^*$.
Effective choices for $\calS^{\Vh}_{\P}(\cdot,\cdot)$ are available from the
virtual element literature~\cite{Mascotto:2018,Dassi-Mascotto:2018}

 \subsubsection{Properties of the inner product~\eqref{eq:Vh:local:inner:product:def}}
\medskip
In the rest of this section, we investigate the properties of the
inner product defined in~\eqref{eq:Vh:local:inner:product:def}.
First, we note that the local bilinear form
$(\,\cdot,\cdot\,)_{\Vh(\P)}$ satisfies the  \textbf{consistency} condition with respect to
the linear polynomials in the sense that
$(\,\qs,\ps\,)_{\Vh(\P)}=(\qs,\ps)_{\LTWO(\P)}$ for every pair of
linear polynomials $\qs,\ps$.
This property is more stringent that the usual consistency of the
typical virtual element constructions, where the exactness property
meaning consistency is true if \emph{at least} one of the entries is a
linear polynomials but not necessarily both simultaneously.

The property that is characterized in the next lemma is the
\textbf{stability} of $(\,\cdot,\cdot\,)_{\Vh(\P)}$ with respect to
the $\LTWO$ inner product.

\begin{lemma}
  \label{lemma:Vh:inner:product:stability}
  There exist two positive constants $\alpha_*$ and $\alpha^*$, which
  are independent of $\hh$ (and $\Dt$), but may depend on the mesh
  regularity parameter $\rho$ and the bounds on $\sigma$, such that
  \begin{align}
    \alpha_*\NORM{\vsh}{0,\P}^2\leq\scalVhP{\vsh}{\vsh}\leq\alpha^*\NORM{\vsh}{0,\P}^2
    \label{eq:Vh:inner:product:stability}
  \end{align}
  for every mesh element $\P$.
\end{lemma}
\BEGINPROOF
Stability is strictly interconnected with the fact that
$(\,\cdot,\cdot\,)_{\Vh(\P)}$ is an inner product in $\Vh(\P)$.
First, we note that $\calS^{\Vh}_{\P}(\cdot,\cdot)$ is a symmetric bilinear
form; hence, the bilinear form $(\,\cdot,\cdot\,)_{\Vh(\P)}$
in~\eqref{eq:Vh:local:inner:product:def} is also symmetric.
The lower bound in~\eqref{eq:SP:stability} implies that
$(\,\cdot,\cdot\,)_{\Vh(\P)}$ is bounded from below by the
$\LTWO(\P)$-norm.
Indeed, note that
\begin{align*}
  \norm{\vsh}{0,\P}^2
  \leq\left(\norm{\PiVhP\vsh}{0,\P}+\norm{(1-\PiVhP)\vsh}{0,\P}\right)^2 
  \leq
  2\Big(\norm{\PiVhP\vsh}{0,\P}^2+\norm{(1-\PiVhP)\vsh}{0,\P}^2\Big).
\end{align*}
Then, a straightforward calculation yields the chain of inequalities:
\begin{align*}
  \scalVhP{\vsh}{\vsh}
  &\geq \sigma_*\norm{\PiVhP\vsh}{0,\P}^2 + \ss_*\norm{(1-\PiVhP)\vsh}{0,\P}^2               \\
  &\geq \min(\sigma_*,\ss_*)\Big(\norm{\PiVhP\vsh}{0,\P}^2+\norm{(1-\PiVhP)\vsh}{0,\P}^2\Big) \nonumber
   \geq \alpha_*\norm{\vsh}{0,\P}^2,
\end{align*}
where we set $\alpha_*=\min(\sigma_*,\ss_*)\slash{2}$.

\medskip
The inequality from above is proved in a similar way:
\begin{align}
  \scalVhP{\vsh}{\vsh}
  &= \scal{\sigma\PiVhP\vsh}{\PiVhP\vsh} + \calS^{\Vh}_{\P}( (1-\PiVhP)\vsh, (1-\PiVhP)\vsh )
  \nonumber\\[0.5em]
  &\leq (\sigma^*+\ss^*)\Big(
  \NORM{\PiVhP\vsh}{0,\P}^2 + \NORM{(1-\PiVhP)\vsh}{0,\P}^2 
  \Big)
  \nonumber
  \leq \alpha^*\NORM{\vsh}{0,\P}^2,
  \label{eq:Vh:local:induced:norm}
\end{align}
where we set $\alpha^*=(\sigma^*+\ss^*)\big(1+\NORM{\PiVhP}{}\big)^2$
in the final step.
\ENDPROOF

\begin{remark}
  A suitable choice of $\calS^{\Vh}_{\P}$ and its scaling factor may
  allow us to have $\ss^*=\sigma^*$. Also, we can define $\PiVhP$ so
  that $\NORM{\PiVhP}{}\leq 1$. This implies that
  $\alpha^*\leq2\sigma^*$ and we can use this bound in the
  inequalities of the next sections to have an explicit dependence on
  $\sigma^*$.
\end{remark}

\medskip
The two properties of symmetry and non-negativity imply that
$(\,\cdot,\cdot\,)_{\Vh(\P)}$ is an inner product in $\Vh(\P)$ for any
element $\P\in\Th$, so that the quantity
\begin{align*}
  \TNORM{\vsh}{\Vh(\P)}^2 := \scalVhP{\vsh}{\vsh}
\end{align*}
is the induced local norm and the Cauchy-Schwarz inequality must hold
\begin{align}
  \scalVhP{\vsh}{\wsh}\leq\TNORM{\vsh}{\Vh(\P)}\TNORM{\wsh}{\Vh(\P)}
  \qquad\forall\vsh,\wsh\in\Vh(\P).
\end{align}
By summing over all the mesh elements, we find that the symmetric
bilinear form defined in~\eqref{eq:Vh:inner:product:def} is bounded
from below by the $\LTWO(\Omega)$-norm.
Therefore, equation~\eqref{eq:Vh:inner:product:def} defines an inner
product on the global virtual element space $\Vh$, with induced
norm given by
\begin{align*}
  \TNORM{\vsh}{\Vh}^2 := \sum_{\P\in\Th}\scalVhP{\vsh}{\vsh}.
\end{align*}
We readily see that such inner product is continuous with respect to its 
induced norm
\begin{align}
  \scalVh{\vsh}{\wsh}
  \leq\TNORM{\vsh}{\Vh}\TNORM{\wsh}{\Vh}
  \qquad\forall\vsh,\wsh\in\Vh,
  \label{eq:Vh:continuity:global}
\end{align}
and such norm is bounded from below by the $\LTWO$ norm
\begin{align}
  \TNORM{\vsh}{\Vh}^2 = \scalVh{\vsh}{\vsh} \geq
  \alpha_*\norm{\vsh}{0,\Omega}^2.
  \label{eq:Vh:nonnegativity:global}
\end{align}
Likewise, in view of Lemma~\ref{lemma:Vh:inner:product:stability}, the
global inner product is also continuous with respect to the
$\LTWO(\P)$-inner product.
In fact, on starting from~\eqref{eq:Vh:continuity:global} and using
the upper bound in~\eqref{eq:Vh:inner:product:stability}, we find that
\begin{align}
  \scalVhP{\vsh}{\wsh}
  \leq\TNORM{\vsh}{\Vh(\P)}\TNORM{\wsh}{\Vh(\P)}
  \leq \alpha^*\NORM{\vsh}{0,\P}\,\NORM{\wsh}{0,\P},
\end{align}
where we recall that $\alpha^*=(\sigma^*+\ss^*)\big(1+\NORM{\PiVh}{}\big)^2$.
By summing all the local terms and noting that $\hh\geq\hP$ for every
$\P$ yields:
\begin{align}
  \scalVh{\vsh}{\wsh}\leq \,\alpha^*
  \norm{\vsh}{0,\Omega}\,\NORM{\wsh}{0,\Omega}.
\end{align}
Therefore, the local inner product is continuous with respect to the
$\LTWO(\P)$-norm for every $\P\in\Th$ and the global inner product in
$\Vh$ is continuous with respect to the $\LTWO(\Omega)$-norm.

\subsubsection{Construction of operator $\PiVhP$}
\label{sec:subsub:construction:oblique:projection}
We discuss three different choices for the approximation operator
$\PiVhP$.

\medskip\noindent
(\textrm{I}).\textbf{Elliptic projection operator}(E). 
The most obvious example of such a computable
approximation operator is the elliptic projection of a virtual element
function $\vsh\in\Vh(\P)$, which is the linear polynomial $\PinP\vsh$
solving the variational problem:
\begin{align}
  \int_{\P}\nabla\PinP\vsh\cdot\nabla\qsh\dV 
  &= \int_{\P}\nabla\vsh\cdot\nabla\qsh\dV
  \qquad\forall\qsh\in\PS{1}(\P),\\[0.5em]
  \frac{1}{N_\V}\sum_{v\in P}\PinP\vsh&=\frac{1}{N_\V} \sum_{v\in P}\vsh.
\end{align}
The elliptic projection $\PinP\vsh$ clearly provides a linear
polynomial approximation of $\vsh$, which is computable from the
degrees of freedom, \textbf{(V1)}, and invariant on linear
polynomials, \textbf{(V2)}, cf.
Ref.~\cite{Ahmad-Alsaedi-Brezzi-Marini-Russo:2013}.
Property \textbf{(V3)} is proved in the appendix,
see Section~\ref{app:sec:proof:V3:PinP}.

\medskip\noindent
(\textrm{II}). \textbf{Least Squares reconstruction operator}(LS).
An alternative to the elliptic projection operator is
provided by the linear interpolant
\begin{align}
  \PiVhPLS\vsh(\xs,\ys) = \as+\bs\frac{\xs-\xs_{\P}}{\hP} +
  \cs\frac{\ys-\ys_{\P}}{\hP},
  \label{eq:least-squares-polynomials}
\end{align}
where the three real coefficients $\as,\bs,\cs$ are determined by
imposing that
\begin{align}
  \PiVhPLS\vsh(\xs_{\V},\ys_{\V}) = \as+\bs\frac{\xs_{\V}-\xs_{\P}}{\hP}
  + \cs\frac{\ys_{\V}-\ys_{\P}}{\hP}
  = \vsh(\xs_{\V},\ys_{\V})
  \qquad\forall\V\in\partial\P,
  \label{eq:least-squares-conditions}
\end{align}
where $\xv_{\V}=(\xs_{\V},\ys_{\V})^T$ is the coordinate position
vector of vertex $\V$.
We solve the resulting system using the Least Squares method.
Indeed, this system has $\NPV$ equations where $\NPV$ is the number of
vertices of the polygonal element and only three unknowns, and is
overdetermined unless $\P$ is a triangular cell.
The linear polynomial $\PiVhPLS\vsh$ only depends on the vertex values
of $\vsh$ and is clearly computable (property~\textbf{(V1)}) and is
invariant on the linear polynomials (property~\textbf{(V2)}).
Property \textbf{(V3)} is proved in the appendix, 
see Section~\ref{app:sec:proof:V3:PiVhPLS}.

\medskip\noindent
(\textrm{III}). \textbf{Galerkin interpolation operator}(GI).
The third alternative that we consider in this paper is
given by a finite element-like piecewise linear interpolant on the
polygonal cell $\P$.
Assumptions~\textbf{(M1)}-\textbf{(M2)} imply the existence of an
internal point $\V^*$ with respect to which $\P$ must be star-shaped
(e.g., the center of the disk in \textbf{(M1)}).
We assume that this point is described by the coordinate vector
\begin{align*}
  \xv^*_{\P} = \sum_{\V\in\partial\P}\omega_{\P,\V}\xvV,
  \qquad\textrm{with}\quad0<\omega_{\P,\V}<1
  \quad\textrm{and}\quad\sum_{\V\in\partial\P}\omega_{\P,\V}=1,
\end{align*}
where the weights $\omega_{\P,\V}$ are known.
For example, if $\P$ is convex, we can choose the arithmetic average
of the vertex positions, so $\omega_{\P,\V}=1\slash{\NPV}$, or the
baricenter of $\P$.
Then, we approximate $\vsh(\xv^*_{\P})$ by the average of the vertex
values using the same weights $\omega_{\P,\V}$:
\begin{align}
  \vsh(\xv^*_{\P})\approx\vshP^* = \sum_{\V\in\partial\P}\omega_{\P\V}\vsh(\xvV).
\end{align}
We note that $\vsh(\xv^*_{\P})=\vshP^*$ if $\vsh$ is a linear
polynomial, which is crucial to ensure that property \textbf{(V2)} is
satisfied.
We connect the internal point $\V^*$ to all the vertices $\V\in\Th$,
thus splitting $\P$ in $\NPV$ subtriangles $\T$ that form a patch
around $\V^*$.
The patch nodes are the vertices of the polygonal boundary of $\P$ and
vertex $\V^*$.
Let $\phi_{\V}(\xvV)$ be the continuous piecewise linear function
defined on the patch that is one at a given patch node (including
vertex $\V^*$) and zero at the other nodes.
Finally, we define the operator $\PiVhP:\HONE(\P)\to\PS{1}(\P_{\hh})$
(see Remark~\ref{remark:codomain}) by
\begin{align}
  \PiVhPpw\vsh(\xv) = \sum_{\V\in\partial\P}\vsh(\xv_{\V})\phi_{\V}(\xv) + \vshP^*\phi_{\V^*}(\xv),
\end{align}
which is the continuous piecewise linear interpolant of $\vsh$ on the
set of values
$\{(\xv_{\V},\vsh(\xv_{\V}))\}_{\V\in\Th}\cup\{(\xv_{\V^*},\vshP^*)\}$. 
From this construction it is obvious that $\PiVhPpw\vsh$ is computable
from the vertex values of $\vsh$ (property~\textbf{(V1)});
$\PiVhPpw\qs=\qs$ if $\qs$ is a linear polynomial
(property~\textbf{(V2)}); $\PiVhPpw\vsh$ is bounded
(property~\textbf{(V3)}) since $0\leq\phi_{\V}(\xv)\leq1$ for all the
patch functions $\phi_{\V}$ at every patch node $\xv_{\V}$ and
$\phi_{\V^*}$ at $\xv^*$.\\
\subsection{Edge virtual element space}
On every element $\P\in\Th$, we consider the following
finite-dimensional space:
\begin{align}
  \Eh(\P) := \big\{
  \vvh\in\HDIV{\P}\cap\HROT{\P}:
  &\restrict{\vvh\cdot\nv}{\E}\in\PS{0}(\E)\,\forall\E\in\partial\P,\,
  \DIV\vvh\in\PS{0}(\P)~\textrm{and}~
  \ROT\vvh=0~\textrm{in}~\P
  \big\}.
  \label{eq:EhP:def}
\end{align}
The local virtual element space $\Eh(\P)$ was introduced in the VEM
literature in Ref.~\cite{BeiraodaVeiga-Brezzi-Marini-Russo:2016c}.
It is worth noting that on a triangular cell, space $\Eh(\P)$
coincides with the space of vector-valued polynomials
$\RT{0}(\P)=(\PS{0}(\P))^2+\PS{0}(\P)\xv$, i.e., those vector-valued
fields that are of the form $\wv(\xv)=\av+\bs\xv$ for some vector and
scalar coefficients $\av\in\REAL^2$ and $\bs\in\REAL$, respectively;
see~\cite{Boffi-Brezzi-Fortin:2013}.
In the case of a general polygonal cell, $\big(\PS{0}(\P)\big)^2$ and
$\RT{0}(\P)$ are clearly subspaces of $\Eh(\P)$.
In view of this elemental definition, we have the corresponding
global virtual element space:
\begin{align}
  \Eh := \Big\{\,\vvh\in\Vspace\,:\,
  \restrict{\vvh}{\P}\in\Eh(\P),\,
  \forall\P\in\Th
  \,\Big\}.
  \label{eq:Eh:def}
\end{align}
By definition, space $\Eh$ is a subspace of $\Vspace$.
Each virtual element function $\vvh\in\Eh(\P)$ is uniquely defined by
the values of its normal components at the edges of $\P$,
$\{\restrict{\vvh\cdot\nv}{\E}\}_{\E\in\partial\P}$.
Similarly, a virtual function in the global space $\Eh$ is defined by
the values of its normal components at the mesh edges.
The unisolvence of this set of degrees of freedom for $\Eh$ is proved
in~\cite{BeiraodaVeiga-Brezzi-Marini-Russo:2016c}.

In $\Eh$ we can compute the two different orthogonal projection
operators denoted by $\PiEhP$ and $\PiEhPRT$, which respectively
project from $\HDIV{\P}$ onto $\big[\PS{0}(\P)\big]^2$ and
$\RT{0}(\P)$.
The orthogonal projection $\PiEhP\vv$ is the constant vector field
solving the variational problem
\begin{align*}
  \int_{\P}\PiEhP\vvh\cdot\qvh\dV = \int_{\P}\vvh\cdot\qvh\dV,
  \qquad\forall\qvh\in\big[\PS{0}(\P)\big]^2.
\end{align*}
This operator is computable from the degrees of
freedom, cf.~\cite{BeiraodaVeiga-Brezzi-Marini-Russo:2016c}.

Also, $\PiEhPRT\vvh$ is the (unique) solution of the following
variational problem:
\begin{align*}
  \int_{\P}\PiEhPRT\vvh\cdot\wvh\dV = \int_{\P}\vvh\cdot\wvh\dV,
  \qquad\forall\wvh\in\RT{0}(\P).
\end{align*}
We show here that $\PiEhPRT\vvh$ is computable from the degrees of
freedom  $\vvh\in\Eh$.
Since $\wvh(\xv)=\av+\bs\xv$, we write it as the gradient of
a second-degree polynomial, i.e., $\wvh=\nabla\qs$ where
$\qs(\xv)=\av\cdot\xv+(\bs/2)\xv^T\xv$.
Then, we substitute this expression for $\wvh$ in the right-hand side,
integrate by parts and obtain:
\begin{align*}
  \int_{\P}\vvh\cdot\wvh\dV
  = \int_{\P}\vvh\cdot\nabla\qs \ \dV
  = -\int_{\P}(\DIV\vvh)\qs \ \dV
  + \sum_{\E\in\partial\P}\int_{\E}\nv\cdot\vvh\,\qs \ \dS.
\end{align*}
All the integrals on the right-hand side are computable.
In fact, the values $\restrict{\nv\cdot\vvh}{\E}$ for all edges
$\E\in\partial\P$ are known as they are the degrees of freedom of
$\vvh$.
Moreover, the divergence of $\vvh$ is also known as it is constant
over $\P$ and a straightforward application of the Gauss Divergence
theorem yields:
\begin{align*}
  \DIV\vvh = \frac{1}{\mP}\sum_{\E\in\partial\P}\mE\restrict{\nv\cdot\vvh}{\E}.
\end{align*}
A similar argument can be used to prove that $\PiEhP\vvh$ is computable
from the degrees of freedom of $\vvh$ (take $\qsh=\av\cdot\xv$), see
Ref.~\cite{BeiraodaVeiga-Brezzi-Marini-Russo:2016c}.

We use the orthogonal projector onto the constant vector fields to
define the inner product in $\Eh$.
As usual in the VEM, we split it into the sum of local contributions:
\begin{align}
  \scalEh{\vvh}{\wvh} = \sum_{\P\in\Th}\scalEhP{\vvh}{\wvh},
\end{align}
where each local term is the inner product in $\Eh(\P)$ and takes the
form
\begin{align}
  \scalEhP{\vvh}{\wvh} = 
  (\PiEhP\vvh,\PiEhP\wvh)_{\P} + 
  \calS^{\Eh}_{\P}\big( (1-\PiEhP)\vvh, (1-\PiEhP)\wvh \big)
\end{align}
and again we assume that $\calS^{\Eh}_{\P}(\cdot,\cdot)$ is a
symmetric and nonnegative bilinear form for which there exist two
positive constant $\bar{\ss}_*$ and $\bar{\ss}^*$ such that
\begin{align*}
  \bar{\ss}_*\NORM{\vvh}{0,\P}^2\leq\calS^{\Eh}_{\P}(\vsh,\vsh)\leq\bar{\ss}^*\NORM{\vvh}{0,\P}^2
  \qquad\forall\vsh\in\Eh(\P)\cap\KER\big(\PiEhP{}\big).
\end{align*}
Since $\PiEhP$ is the orthogonal projection onto the constant
vector-valued fields defined on $\P$, it is now easy to prove that
this inner product is consistent and stable in the usual VEM sense; namely,
\begin{itemize}
\item \textbf{consistency}:
  \begin{align}
    \scalEhP{\vvh}{\qvh} = \int_{\P}\vvh\cdot\qvh\dV
    \qquad \qvh\in\big(\PS{0}(\P)\big)^2;
    \label{eq:Eh:inner-product:consistency}
  \end{align}
\item \textbf{stability}: there exist two positive constants,
  $\bar{\alpha}_*$ and $\bar{\alpha}^*$, such that
  \begin{align}
    \bar{\alpha}_*\NORM{\vvh}{0,\P}^2\leq\scalEhP{\vvh}{\vvh}\leq\bar{\alpha}^*\NORM{\vvh}{0,\P}^2
    \qquad\forall\vvh\in\Eh(\P).
    \label{eq:Eh:inner-product:stability}
  \end{align}
\end{itemize}

\subsection{Cell space}

On every element $\P\in\Th$, we consider the finite-dimensional space
$\Ph(\P):=\PS{0}(\P)$, which is the space of constant functions
defined on $\P$.
The corresponding global space is
\begin{align}
  \Ph := \Big\{\,\qsh\in\LTWO(\Omega)\,:\,
  \restrict{\qsh}{\P}\in\Ph(\P),\,
  \forall\P\in\Th
  \,\Big\};
  \label{eq:Ph:def}
\end{align}
Space $\Ph$ is the space of piecewise constant functions 
$\qsh\in\LTWO(\P)$
defined on mesh $\Th$.
So, the degrees of freedom of $\qsh$ are the values that $\qsh$
takes in each mesh cell, namely, $\restrict{\qsh}{\P}$.

\subsection{Interpolation operators and approximation of
  $\scal{\sigma\uv\times\Bv}{\vs}$}
\label{subsec:VirtualdeRham}
We define the local interpolation operators
\begin{equation}
  \IVhP:\HS{1}(\P)\to\Vh(\P),\;\;
  \IEhP:\Vspace\to\Eh(P)\;\;\mbox{and}\;\;
  \IPhP:\LTWO (\P)\to\Ph(\P),
\end{equation}
by requiring that 
\begin{itemize}
\item for any scalar function
  $\vs\in\HONE(\P)\cap\CS{0}(\overline{\P})$, it holds
  $\IVhP\vs(\V)=\vs(\V)$, for every vertex $\V\in\partial\P$;
\item for any vector-valued function $\wv\in\HDIV{\P}\cap\HROTv{\P}$,
  it holds
  \begin{align*}
    \norE\cdot\IEhP(\wv) 
    = \frac{1}{\mE}\int_{\E}\norE\cdot\IEhP(\wv)\ \dS 
    = \frac{1}{\mE}\int_{\E}\norE\cdot\wv\ \dS,
  \end{align*}
  for every edge $\E\in\partial\P$;
\item for any scalar function $\qs\in\LTWO(\P)$,
  \begin{align*}
    \int_{\P}\IPhP\qs\dV=\int_{\P}\qs\dV.
  \end{align*}
\end{itemize}
Correspondingly, we define the global interpolation operators by
pasting together the elementwise operators
\begin{equation}
  \restrict{(\IVh\vs)}{\P}=\IVhP\big(\restrict{\vs}{\P}\big),\qquad
  \restrict{(\IEh\vv)}{\P}=\IEhP\big(\restrict{\vv}{\P}\big),\quad\textrm{and}\quad
  \restrict{(\IPh\qs)}{\P}=\IPhP\big(\restrict{\qs}{\P}\big).
\end{equation}
It is easy to see that these interpolation operators are continuous
\begin{align}
  \TNORM{\IVh\vs}{\Vh}&\!\leq \Cs\NORM{\vs}{0,\Omega} \phantom{\NORM{\wv}{0,\Omega}\NORM{\vs}{0,\Omega}} \hspace{-1.5cm}\forall\vs\in\HONE(\Omega),\label{eq:IVh:continuity}\\[0.5em]
  \TNORM{\IEh\wv}{\Eh}&\!\leq \Cs\NORM{\wv}{0,\Omega} \phantom{\NORM{\vs}{0,\Omega}\NORM{\vs}{0,\Omega}} \hspace{-1.5cm}\forall\wv\in\Vspace,      \label{eq:IEh:continuity}\\[0.5em]
  \TNORM{\IPh\vs}{\Ph}&\!\leq \Cs\NORM{\vs}{0,\Omega} \phantom{\NORM{\vs}{0,\Omega}\NORM{\wv}{0,\Omega}} \hspace{-1.5cm}\forall\vs\in\LTWO(\Omega).\label{eq:IPh:continuity}
\end{align}

Finally, we use the interpolation operator $\IVh$ and the orthogonal
projection operator $\PiEhRT$ to approximate the term involving
$\uv\times\Bv$ as follows:
\begin{align}
  \Scal{\sigma\uv\times\Bv}{\vs}\approx\scalVh{\IVh(\uv\times\PiEhRT\Bvh)}{\vsh},
  \label{eq:Vh:inner-product-RT}
\end{align}
where all the terms on the right have been defined except the
$RT$-orthogonal projection of $\Bvh\in\Eh$, which must be such that
$\restrict{\big(\PiEhRT\Bvh\big)}{\P}=\PiEhPRT\big(\restrict{\Bvh}{\P}\big)$
for every mesh cell $\P\in\Th$.
Note that the coefficient $\sigma$ is incorporated into the definition
of the inner product in accordance with
definition~\eqref{eq:Vh:local:inner:product:def}.
We conclude this section with a technical lemma that provides a useful
estimate for the term in~\eqref{eq:Vh:inner-product-RT}.

\medskip
\begin{lemma}
  \label{lemma:Vh:inner-product:estimate}
  There exists a real positive constant $\Cst$ 
  independent of $\hh$ (and $\Dt$) that may depend on $\alpha^*$ 
  and the continuity constants of $\IVh$ and $\PiEhRT$, such that
  \begin{align}
    \scalVh{\IVh\big(\uv\times\PiEhRT\wvh\big)}{\vsh} \leq
    \Cst\NORM{\uv}{\infty} \NORM{\wvh}{0,\Omega}\,\NORM{\vsh}{0,\Omega}
    \label{eq:Vh:inner-product:estimate}
  \end{align}
  for every $\wvh\in\Eh$, $\vsh\in\Vh$, and any assigned velocity
  $\uv\in\LINF(\Omega)$.
\end{lemma}
\BEGINPROOF
\begin{align*}
  \begin{array}{rll}
    &\hspace{-1cm}
    \scalVh{\IVh\big(\uv\times\PiEhRT\wvh\big)}{\vsh}
    \leq \TNORM{\IVh\big(\uv\times\PiEhRT\wvh\big)}{\Vh}\,\TNORM{\vsh}{\Vh}
    & \qquad\mbox{[use~\eqref{eq:Vh:nonnegativity:global}]}
    \\[0.35em]
    &\leq (\alpha^*)^{\HALF}
    \NORM{\IVh\big(\uv\times\PiEhRT\wvh\big)}{0,\Omega}\,\NORM{\vsh}{0,\Omega}           
    & \qquad\mbox{[use~\eqref{eq:IVh:continuity}]}
    \\[0.35em]
    &\leq (\alpha^*)^{\HALF}\NORM{\IVh}{}
    \NORM{\uv\times\PiEhRT\wvh}{0,\Omega}\,\NORM{\vsh}{0,\Omega}           
    & \qquad\mbox{[note that $\NORM{\uv}{\infty}<\infty$]}
    \\[0.35em]
    &\leq (\alpha^*)^{\HALF}\NORM{\IVh}{}\NORM{\uv}{\infty}
    \NORM{\PiEhRT\wvh}{0,\Omega}\,\NORM{\vsh}{0,\Omega}           
    & \qquad\mbox{[note that $\NORM{\PiEhRT}{}\leq1$]}
    \\[0.35em]
    &\leq (\alpha^*)^{\HALF}\NORM{\IVh}{}\NORM{\uv}{\infty}
    \NORM{\wvh}{0,\Omega}\,\NORM{\vsh}{0,\Omega}    ,       
  \end{array}
\end{align*}
which is the assertion of the lemma after setting
$\Cst=(\alpha^*)^{\HALF}\NORM{\IVh}{}$.
\ENDPROOF

\subsection{Commuting properties and the virtual De-Rham complex}
The elementwise interpolation operators $\IVhP$, $\IEhP$ and $\IPhP$
for every mesh element $\P\in\Th$ commute with the differential
operators $\ROTv$ and $\DIV$.
We state this property in the next lemma.
\begin{lemma}[Commutation properties]
  \label{lemma:commutation:properties}
  \begin{align*}
    \begin{array}{rll}
      (\mathrm{i})  & \ROTv\circ\IVhP &= \IEhP\circ\ROTv \textrm{~in~}\Eh(\P),\qquad\forall\P\in\Th,\\[0.5em]
      (\mathrm{ii}) & \DIV \circ\IEhP &= \IPhP\circ\DIV  \textrm{~in~}\Ph(\P),\qquad\forall\P\in\Th.
    \end{array}
  \end{align*}
\end{lemma}
\BEGINPROOF
In view of the unisolvence of the degrees of freedom in
$\Eh(\P)$~\cite{BeiraodaVeiga-Brezzi-Marini-Russo:2016c}, to prove
$(\mathrm{i})$ we only need to show that the degrees of freedom of
$\ROTv\big(\IVhP\vs\big)$ are equal to the degrees of freedom of
$\IEhP\big(\ROTv\vs\big)$.
Consider $\vs\in\HONE(\P)$ and its interpolant
$\vsh=\IVhP\vs\in\Vh(\P)$, whose degrees of freedom are the vertex
values $\vsh(\V)=\vs(\V)$, $\V\in\partial\P$, and recall, for every
edge $\E\in\partial\P$, that
\begin{align*}
  \norE
  =\left(\begin{array}{c} \ns^{\E}_{x}\\\ns^{\E}_{y}\end{array}\right)
  =\left(\begin{array}{c}\ts^{\E}_{y}\\ -\ts^{\E}_{x}\end{array}\right)
  =\left(\begin{array}{rr}0&\ \ 1\\-1 & 0\end{array}\right)\tngE.
\end{align*}
A straightforward calculation shows that
\begin{equation*}
  \norE\cdot\ROTv\vsh
  = \ns^{\E}_{x}\frac{\partial\vsh}{\partial\ys} - \ns^{\E}_{y}\frac{\partial\vsh}{\partial\xs}
  = \ts^{\E}_{x}\frac{\partial\vsh}{\partial\xs} + \ts^{\E}_{y}\frac{\partial\vsh}{\partial\ys}
  = \tngE\cdot\nabla\vsh,
  \nonumber\\[0.5em]
\end{equation*}
which by the fundamental theorem of line integrals yields that
\begin{equation*}
\frac{1}{\mE}\int_{\E}\norE\cdot\ROTv \vs \ \dS=
\frac{1}{\mE}\int_{\E}\tngE\cdot\nabla\vs \ \dS=
\frac{\vsh(\Vpp)-\vsh(\Vp)}{\mE}. 
\end{equation*}
\medskip
\noindent
Similarly, to prove $\mathrm{ii})$, we only need to show that for any
$\wv\in\HDIV{\P}$, the degrees of freedom of $\DIV(\IEhP\wv)$ in
$\Ph(\P)$ are equal to the degrees of freedom of $\IPhP(\DIV\wv)$.
This fact is evident from the following chain of identities:
\begin{align*}
  \restrict{\DIV(\IEhP\wv)}{\P}
  &=\frac{1}{\mP}\int_{\P}\DIV(\IEhP\wv)\ \dV
  = \frac{1}{\mP}\int_{\partial\P}\nv\cdot\IEhP\wv\ \dS
  = \frac{1}{\mP}\sum_{\E\in\partial\P}\int_{\E}\norE\cdot\IEhP\wv\ \dS
  \nonumber\\[0.5em]
  &= \frac{1}{\mP}\sum_{\E\in\partial\P}\int_{\E}\norE\cdot\wv\ \dS
  = \frac{1}{\mP}\int_{\partial\P}\nv\cdot\wv\ \dS
  = \frac{1}{\mP}\int_{\P}\DIV\wv\ \dV
  = \restrict{\IPhP(\DIV\wv)}{\P}.
\end{align*}
\ENDPROOF

\begin{theorem}\label{ExactChainTheorem}
  The de Rham diagram
  \begin{align*}
    \begin{CD}
      \HROTv{\Omega} @> \ROTv >> \HDIV{\Omega} @> \DIV >> \LTWO(\Omega)\\
      @VV \IVh V  @VV \IEh V @VV \IPh V \\
      \Vh  @> \ROTv >> \Eh @> \DIV >> \Ph
    \end{CD}
  \end{align*}
  is commutative and the chain
  \begin{align*}
    \begin{CD}
      \Vh @> \ROTv >> \Eh @> \DIV >> \Ph 
    \end{CD}
  \end{align*}
  is short and exact.
\end{theorem}
\BEGINPROOF
Consider a virtual element function $\wvh\in\Eh$ whose restriction to
every element $\P\in\Th$ has zero divergence, i.e.,
$\DIV(\restrict{\wvh}{\P})=0$.
Since from Assumption~\textbf{(M1)}-\textbf{(M2)}, element $\P$ is
simply connected, there exists a function $\vs$ in $\HONE(\P)$ such
that $\wvh=\ROTv\vs$.
Let $\vsh=\IVhP\vs$.
Lemma~\ref{lemma:commutation:properties}-$(\mathrm{i})$, and the fact that
$\restrict{\wvh}{\P}=\IEhP(\wvh)$ and
$\restrict{\vsh}{\P}=\IVhP(\vsh)$, imply that
\begin{align*}
  \restrict{\wvh}{\P} = \IEhP(\wvh) = \IEhP\big(\ROTv\vs\big) =
  \ROTv\big(\IVhP\vs\big) = \ROTv(\vsh),
\end{align*}
for every $\P\in\Th$.
The left-most part of the de Rham complex follows by considering
together all the elemental commuting relations.

Similarly, consider a piecewise constant function $\qsh\in\Ph$, and
let $\wv\in\HDIV{\Omega}$ be the vector-valued field whose divergence
reproduces the elemental values of $\qsh$ when restricted to the mesh
elements, i.e., $\restrict{\qsh}{\P}=\DIV(\restrict{\wv}{\P})$.
Let $\wvh=\IEh(\wv)$.
Lemma~\ref{lemma:commutation:properties}-$(\mathrm{ii})$, and the fact that
$\restrict{\wvh}{\P}=\IEhP(\wvh)$ and
$\restrict{\qsh}{\P}=\IPhP(\qsh)$ imply that
\begin{align*}
  \restrict{\qsh}{\P} = \IPhP(\qsh) = \IPhP\big(\DIV\wv\big) =
  \DIV\big(\IEhP\wv\big) = \DIV(\restrict{\wvh}{\P}),
\end{align*}
for every $\P\in\Th$.
The right-most part of the de Rham complex follows by considering
together all the elemental commuting relations.

\ENDPROOF

\section{Wellposedness of the Virtual Element Method}
\label{sec:wellposedness:VEM}

Inspired by~\cite{hu2017stable}, in this section, we investigate the
wellposedness of the virtual element method that we presented in the
previous section.
The major result of this section is stated by the following theorem.

\medskip
\begin{theorem}\label{theorem:wellposedness:problem:A}
  If $\theta>0$, then solution to Problem~\ref{problem:A}
  exists and is unique.
  Moreover, the map $(\Fv,\gs)\to(\Bvh^{n+1},\hEsh^{n+\theta})$ is
  uniformly continuous independently of $\hh$ and $\Dt$ in the norm
  defined in $\Xh$.
\end{theorem}

\medskip
The definition of the space $\Xh$ and its norm will be presented in the next section whereas the proof of this theorem will be postponed at the end of the
section since it requires some further investigation about the properties of
the VEM.
In particular, we will follow this roadmap.
First, we prove that the approximation of the magnetic flux field is
divergence free provided that such condition is satisfied at the
initial time.
Second, we reformulate the $(n+1)$-step of scheme
\eqref{eq:VarForEq4} in a suitable way, cf. Problem~\ref{problem:A}
below, and introduce two additional problems, namely,
Problem~\ref{problem:B} and Problem~\ref{problem:C}.
Third, we prove that these three problems are equivalent,
cf. Theorem~\ref{theorem:equivalence-problems:ABC}, and, finally, that
Problem~\ref{problem:C} is wellposed as a consequence of
Babuska-Lax-Milgram Theorem~\cite{Babuska:1970},
These facts eventually imply the wellposedness of
Problem~\ref{problem:A}.

\medskip
To prove the equivalence of Problems~\ref{problem:A} \ref{problem:B}
and~\ref{problem:C} we need two additional theorems stating that
$\DIV\Bvh^{n+1}=0$ whenever $\DIV\IEh\Bvzr=0$.
These intermediate results confirm that the virtual element
approximation $\Bvh$ to the magnetic flux field satisfies the
divergence free condition.

\medskip
We start by reformulating the $(n+1)$-th step of scheme
\eqref{eq:VarForEq4} as follows.
\begin{problem}\label{problem:A}
  Suppose that $\Bvh^{n}$ and $\hEsh^{n-1+\theta}$ are known.
  Then, the $(n+1)$-th step of scheme \eqref{eq:VarForEq4} can be
  written as:
  \emph{Find $(\Bvh^{n+1},\hEsh^{n+\theta})\in\Eh\times\Vhzr$ such that
    for all $(\wvh,\vsh)\in\Eh\times\Vhzr$ it holds:}
  \begin{align}
    &\Dt^{-1}\scalEh{\Bvh^{n+1}}{\wvh} + \scalEh{\ROTv\hEsh^{n+\theta}}{\wvh} = \scalEh{\Fv}{\wvh},
    \label{eq:OnstepFaraday}\\[0.25em]
    &\scalVh{\hEsh^{n+\theta}}{\vsh}+\theta\scalVh{\IVh\big(\uv\times\PiEhRT\Bvh^{n+1}\big)}{\vsh}-\theta\scalEh{\Bvh^{n+1}}{\ROTv\vsh} = \bil{\gs}{\vsh},
    \label{eq:OneStepAmpereOhm}
  \end{align}
  where we define
  \begin{align}
    &\Fv = \Dt^{-1}\Bvh^{n}+\ROTv\big(\IVh\Es_0^{n+\theta}\big),\label{eq:F:def}\\
    &\bil{\gs}{\vsh} = (1-\theta)\left(\scalEh{\Bvh^n}{\ROTv\vsh}-\scalVh{\IVh\big(\uv\times\PiEhRT\Bvh^n\big)}{\vsh}\right)-\scalVh{\IVh\Es_0^{n+\theta}}{\vsh}.
    \label{eq:bil:def}
  \end{align}
\end{problem}

Next we show some results regarding the stability of scheme
\eqref{eq:VarForEq4}.

\subsection{Abstract setting and equivalent problems}
To have a setting to analyze Problem~\ref{problem:A}, we introduce the
space $\Xh:=\Esh\times\Vhzr$.
We set $(\Bvh,\Esh)=\xi\in\Xh$ and equip $\Xh$ with the norm
\begin{align}
  &\TNORM{\xi}{\Xh}^2 := \CurlNorm{\Esh}^2 + \DivNorm{\Bvh}^2,\label{eq:norm-Xh:def}\\
  \intertext{where}
  &\CurlNorm{\Esh}^2  := \TNORM{\Esh}{\Vh}^2 + \Dt\TNORM{\ROTv\Esh}{\Eh}^2,\label{eq:CurlNorm:def}\\[0.25em]
  &\DivNorm {\Bvh}^2  : = \Dt^{-1}\TNORM{\Bvh}{\Eh}^2+\NORM{\DIV\Bvh}{0}^2.\label{eq:DivNorm}
\end{align}
The space $\Xh$ is complete in the topology induced by norm
$\TNORM{\,\cdot\,}{\Xh}$.

Next, we introduce two additional variational problems.
To formulate such problems, we define the two bilinear forms
$\ash:\Xh\times\Xh\to\REAL$ and $\ashzr:\Xh\times\Xh\to\REAL$.
Let $\xi=(\Bvh,\Esh)$ and $\eta=(\wvh,\vsh)$.
The first bilinear form is given by
\begin{align}
  \ash(\xi,\eta) 
  = \scalEh{\Dt^{-1}\Bvh+\ROTv\Esh}{\wvh}
  + \scalVh{\Esh+\theta\IVh\big(\uv\times\PiEhRT\Bvh\big)}{\vsh}
  - \theta\scalEh{\Bvh}{\ROTv\vsh}.
  \label{eq:ash:def}
\end{align}
The second bilinear form is given by
\begin{equation}
  \ashzr(\xi,\eta) = \ash(\xi,\eta) + \scal{\DIV\Bvh}{\DIV\wvh}.
  \label{eq:ashzr:def}
\end{equation}

\medskip
\noindent
The first auxiliary variational problem reads as follows.
\begin{problem}\label{problem:B}
  \emph{Find $(\Bvh^{n+1},\hEsh^{n+\theta})=\xi\in\Xh$ such that for any
    $(\wvh,\vsh)=\eta\in\Xh$ it holds:}
  \begin{equation}
    \ash(\xi,\eta) = \scalEh{\Fv}{\wvh} + \bil{\gs}{\vsh},
    \label{eq:Onestepeq2}
  \end{equation}
  where $\Fv$ and $\bil{\gs}{\vsh}$ are given by~\eqref{eq:F:def}
  and~\eqref{eq:bil:def} assuming that $\Bvh^{n}$ (such that
  $\DIV\Bvh^{n}=0$) and $\hEsh^{n-1+\theta}$ are known.
\end{problem}

\medskip
\noindent
The second auxiliary variational problem reads as follows.
\begin{problem}\label{problem:C}
  \emph{ Find $(\Bvh^{n+1},\hEsh^{n+\theta})=\xi\in\Xh$ such that for
    any $(\wvh,\vsh)=\eta\in\Xh$:}
  \begin{equation}
    \ashzr(\xi,\eta) = \scalEh{\Fv}{\wvh} + \bil{g}{\vsh}.
    \label{eq:OneStepeq3}
  \end{equation}
  where $\Fv$ and $\bil{\gs}{\vsh}$ are given by~\eqref{eq:F:def}
  and~\eqref{eq:bil:def} assuming that $\Bvh^{n}$ (such that
  $\DIV\Bvh^{n}=0$) and $\hEsh^{n-1+\theta}$ are known.
\end{problem}

\medskip
\begin{theorem}[Zero-divergence magnetic flux from
  system~\eqref{eq:VarForEq4}]
  \label{theorem:DivFreeScheme}
  Let $\{\Bvh^n\}_{n=0}^{N}\subset\Eh$ and
  $\{\Esh^{n+\theta}\}_{n=0}^{N}\subset\Vhzr$ be the solution of the
  virtual element scheme~\eqref{eq:VarForEq4}, with
  $\Bvh^{0}=\IEh\Bvzr$ and $\DIV\Bvzr=0$.
  Then, $\DIV\Bv^{n}=0$ for every $0\leq\ns\leq\Ns$.
\end{theorem}

\BEGINPROOF
First, Lemma~\ref{lemma:commutation:properties}-$(\mathrm{ii})$ implies that
$\DIV\Bvh^{0}=\DIV\big(\IEh\Bvzr\big)=\IPhP(\DIV\Bvzr)=0$ since we
assume that $\Bvh^{0}=\IEh\Bvzr$ with $\DIV\Bvzr=0$.
Then, we observe that $\ROTv\Esh^{n+\theta}\in\Eh$ for every
$\Esh^{n+\theta}\in\Vh$.
Therefore, equation~\eqref{eq:VarForEq4Faraday} states that
\begin{align}
  \Bvh^{n+1}-\Bvh^n = \Dt \ \ROTv\Esh^{n+\theta} \qquad\textrm{in}\,\,\Eh
  \label{eq:zero-divg:proof:00}
\end{align}
for every $n\geq0$.
Taking the divergence of both sides of \eqref{eq:zero-divg:proof:00},
we find that $\DIV\Bvh^{n+1}=\DIV\Bvh^{n}$.
We apply this relation recursively back to $n=0$ and find that
$\DIV\Bvh^{n}=\ldots=\DIV\Bvh^{0}=0$, which is
the assertion of the theorem.
\ENDPROOF

\medskip
\begin{theorem}[Zero-divergence magnetic flux from
  Problem~\ref{problem:C}]
  \label{theorem:zero-divergence-B:problem:C}.
  If $\xi = (\Bvh^{n+1},\hEsh^{n+\theta})$ solves
  Problem~\ref{problem:C}, then $\DIV\Bvh^{n+1}=0$.
\end{theorem}

\BEGINPROOF
Test~\eqref{eq:OneStepeq3} against $\eta=(\wvh,\vsh)$ with $\vsh=0$,
while leaving $\wvh\in\Eh$ undefined for the moment.
Using definitions~\eqref{eq:ashzr:def}, \eqref{eq:ash:def},
\eqref{eq:F:def}, and~\eqref{eq:bil:def}, and rearranging the terms,
we obtain the identity:
\begin{align}
  \scalEh{\Fv^{n}-\Dt^{-1}\Bvh^{n+1}-\ROTv\Esh^{n+\theta}}{\wvh}
  = \scal{\DIV\Bvh^{n+1}}{\DIV\wvh}.
  \label{eq:prop:00}
\end{align}
Now, we set 
\begin{align*}
  \wvh = \Fv^{n}-\Dt^{-1}\Bvh^{n+1}-\ROTv\,\Esh^{n+\theta}.
\end{align*}
Since $\DIV\Bvh^{n}=0$ by hypothesis and $\DIV\circ\ROTv=0$ we find
that
\begin{align*}
  \DIV\Fv^{n} = \Dt^{-1}\DIV\Bvh^{n}+\DIV\big(\ROTv\hEsh^{n-1+\theta}\big) = 0
  \qquad\textrm{and}\quad
  \DIV\big(\ROTv\Esh^{n+\theta}\big)=0,
\end{align*}
so that
\begin{align*}
  \DIV\wvh = 
  \DIV(\Fv^{n}-\Dt^{-1}\Bvh^{n+1}-\ROTv\hEsh^{n+\theta})
  =-\Dt^{-1}\DIV\Bvh^{n+1}.
\end{align*}
Substituting the expressions of $\wv$ and $\DIV\wv$
in~\eqref{eq:prop:00} yields
\begin{align*}
  0\leq\TNORM{\wvh}{\Eh}^2 = -\Dt^{-1}\NORM{\DIV\Bvh^{n+1}}{0,\Omega}^2,
\end{align*}
which implies that $\NORM{\DIV\Bvh^{n+1}}{0,\Omega}\leq0$, and, thus,
the proposition.
\ENDPROOF

\medskip
\begin{theorem}[Equivalence of
  Problems~\ref{problem:A},~\ref{problem:B}, and~\ref{problem:C}]
  \label{theorem:equivalence-problems:ABC}
  Problems~\ref{problem:A}, \ref{problem:B} and \ref{problem:C} are equivalent.
\end{theorem}

\BEGINPROOF
It is immediate to see that Problem~\ref{problem:A} is equivalent to
Problem~\ref{problem:B}.
In fact, adding \eqref{eq:OnstepFaraday} and
\eqref{eq:OneStepAmpereOhm} yields \eqref{eq:Onestepeq2}, while
testing \eqref{eq:Onestepeq2} against $\eta = (\wvh,0)$ yields
\eqref{eq:OnstepFaraday} and against $\eta = (0,\vsh)$ yields
\eqref{eq:OneStepAmpereOhm}.

\medskip
To prove that Problem~\ref{problem:B} is equivalent to
Problem~\ref{problem:C}, we use the result of
Theorem~\ref{theorem:zero-divergence-B:problem:C}.
In light of this theorem, if $\xi=(\Bvh^{n+1},\hEsh^{n+\theta})$
solves Problem~\ref{problem:C}, then $\DIV\Bvh^{n+1}=0$, and
$\ashzr(\xi,\eta)=\ash(\xi,\eta)$ for every $\eta\in\Xh$, so $\xi$ is
also a solution of Problem~\ref{problem:B}.
Instead, if $\xi=(\Bvh^{n+1},\hEsh^{n+\theta})$ solves
Problem~\ref{problem:B}, then it is also a solution of
Problem~\ref{problem:A}, and $\DIV\Bvh^{n+1}=0$ in view of
Theorem~\ref{theorem:DivFreeScheme}.
Therefore, we can conclude that
$\ashzr(\xi,\eta)=\ash(\xi,\eta)$ for every $\eta\in\Xh$ and $\xi$
must be a solution of Problem~\ref{problem:C}.
\ENDPROOF

To prove that Problem~\ref{problem:C} is well-posed, we prove that the
bilinear form $\ashzr(\cdot,\cdot)$ and the linear functionals $\scalVh{\Fv}{\cdot}$, $\bil{g}{\cdot}$ satisfy the hypothesis of the
Babuska-Lax-Milgram theorem~\cite{atkinson2005theoretical}.
First, we prove that $\ashzr(\cdot,\cdot)$ is continuous
\begin{lemma}\label{ContinuityA0}
  There exists a constant $C>0$, independent of $\hh$ and $\Dt$, such
  that
  \begin{align}
    \forall\xi,\eta\in\Xh:\quad&\ashzr(\xi,\eta)\le C\TNORM{\xi}{\Xh}\,\TNORM{\eta}{\Xh}.
  \end{align}
\end{lemma}

\BEGINPROOF
Let $\xi=(\Bvh,\Esh)$ and $\eta=(\wvh,\vsh)$ be arbitrary elements in
$\Xh$. 
A systematic application of the Cauchy Schwartz inequality yields that
\begin{align*}
  \Dt^{-1}\scalEh{\Bvh}{\wvh}
  &\leq \Dt^{-\HALF}\TNORM{\Bvh}{\Eh}\,\Dt^{-\HALF}\TNORM{\wvh}{\Eh}
  \leq \DivNorm{\Bvh}\,\DivNorm{\wvh},
  \\[0.5em]
  \scalEh{\ROTv\Esh}{\wvh} 
  &\leq \Dt^{\HALF}\TNORM{\ROTv\Esh}{\Eh}\,\Dt^{-\HALF}\TNORM{\wvh}{\Eh}
  \leq \CurlNorm{\Esh}\,\DivNorm{\wvh},
  \\[0.5em]
  \scalVh{\Esh}{\vsh}
  &\leq \TNORM{\Esh}{\Vh}\,\TNORM{\vsh}{\Vh}
  \leq \CurlNorm{\Esh}\,\CurlNorm{\vsh},
  \\[0.5em]
  \scal{\DIV\Bvh}{\DIV\wvh}
  &\leq \NORM{\DIV\Bvh}{0,\Omega}\,\NORM{\DIV\wvh}{0,\Omega}
  \leq \DivNorm{\Bvh}\,\DivNorm{\wvh}.
\end{align*}
We recall that the Friedrichs-Poincar\'e inequality holds so that
$\NORM{\vsh}{0,\Omega}\leq\Cs\NORM{\nabla\vsh}{0,\Omega}$ for every
$\vsh\in\Vhzr\subset\HONEzr(\Omega)$ and note that
$\NORM{\nabla\vsh}{0,\Omega}=\NORM{\ROTv\vsh}{0,\Omega}$.
In view of Lemma~\ref{lemma:Vh:inner-product:estimate}, we find that
\begin{align*}
  \begin{array}{rll}
    &\hspace{-1cm}
    \scalVh{\IVh\big(\uv\times\PiEhRT\Bvh\big)}{\vsh}
    \leq \Cst\NORM{\uv}{\infty}
    \NORM{\Bvh}{0,\Omega}\,\NORM{\vsh}{0,\Omega}           
    & \qquad\mbox{[use Poincar\'e inequality]}
    \\[0.35em]
    &\leq \Cst\NORM{\uv}{\infty}
    \NORM{\Bvh}{0,\Omega}\,\NORM{\ROTv\vsh}{0,\Omega}      
    & \qquad\mbox{[use stability condition~\eqref{eq:Eh:inner-product:stability}]}
    \\[0.35em]
    &\leq \Cst\NORM{\uv}{\infty}
    \TNORM{\Bvh}{\Eh}\,\TNORM{\ROTv\vsh}{\Eh}          
    & \qquad\mbox{[multiply and divide by $\Dt^{\HALF}$]}
    \\[0.35em]
    &\leq \Cst\NORM{\uv}{\infty}
    \Dt^{-\HALF}\TNORM{\Bvh}{\Eh}\,\Dt^{\HALF}\TNORM{\ROTv\vsh}{\Eh}
    & \qquad\mbox{[use definitions~\eqref{eq:CurlNorm:def} and~\eqref{eq:DivNorm}]}
    \\[0.35em]
    &\leq \Cst\NORM{\uv}{\infty}
    \DivNorm{\Bvh}\,\CurlNorm{\vsh}\,
    & \qquad\mbox{[use definition~\eqref{eq:norm-Xh:def} ]}
    \\[0.35em]
    &\leq \Cst\NORM{\uv}{\infty}\TNORM{\xi}{\Xh}\,\NORM{\eta}{\Xh}.
  \end{array}
\end{align*}
\medskip
The assertion of the lemma follows from the definition of the norm in
$\Xh$ and the above estimates.
\ENDPROOF
 
The next lemma will show that $\ashzr(\cdot,\cdot)$ satisfies the
\emph{inf-sup} condition.
\begin{lemma}\label{lemma:infsup}
  Let $\theta>0$. 
  Then, for a sufficiently small $\Dt$, there exists a real positive
  constant $\Csb$, independent of $\hh$ and $\Dt$, such that:
  \begin{equation}
    \inf_{\xi\in\Xh}\sup_{\eta\in\Xh}\frac{\ashzr(\xi,\eta)}{\TNORM{\xi}{\Xh}\TNORM{\eta}{\Xh}}\geq\Csb>0.
  \end{equation}
  The constant $\Csb$ depends on parameter $\theta$ (and the mesh
  regularity parameter $\rho$).
\end{lemma}

\BEGINPROOF
The assertion of the lemma follows from proving that for every
$\xi=(\Bvh,\Esh)\in\Xh$ there exists a $\eta_{\xi}\in\Xh$ such that
$\TNORM{\eta_{\xi}}{\Xh}\leq\Cs\TNORM{\xi}{\Xh}$, and
\begin{align}
  \ashzr(\xi,\eta_{\xi})\geq\Csb\TNORM{\xi}{\Xh}\TNORM{\eta_{\xi}}{\Xh},
\end{align}
where both $\Cs$ and $\Csb$ are real positive constants independent
of $\hh$ and $\Dt$.
To this end, we first split the bilinear form in~\eqref{eq:ashzr:def}
as follows
\begin{align}
  \ashzr(\xi,\eta) = \TERM{T1} + \TERM{T2},
  \label{eq:ash0:split}
\end{align}
where 
\begin{align}
  \TERM{T1} &= \scalEh{\Dt^{-1}\Bvh+\ROTv\Esh}{\wvh} + \scal{\DIV\Bvh}{\DIV\wvh},\\[0.5em]
  \TERM{T2} &= \scalVh{\Esh+\theta\IVh\big(\uv\times\PiEhRT\Bvh\big)}{\vsh} - \theta\scalEh{\Bvh}{\ROTv\vsh}.
\end{align}
Then, for an arbitrary pair $\big(\Bvh,\Esh\big)=\xi\in\Xh$, we
consider the pair $\big(\wvh,\vsh\big)=\eta_{\xi}\in\Xh$ with
$\wvh=(\theta\slash{2})\,\big(\Bvh+\Dt\ROTv\Esh\big)$ and $\vsh=\Esh$.
Note that $\DIV\wvh=(\theta\slash{2})\DIV\Bvh$ because
$\DIV(\ROTv\Esh)=0$.
Substituting $\xi$ and $\eta$ we transform the first term
in~\eqref{eq:ash0:split} as follows:
\begin{align*}
  \TERM{T1} 
  &= \frac{\theta}{2}\Big( 
  \scalEh{\Dt^{-1}\Bvh+\ROTv\Esh}{\Bvh+\Dt\ROTv\Esh} +
  \scal{\DIV\Bvh}{\DIV\Bvh}
  \Big)\\[0.5em]
  &= \frac{\theta}{2}
  \Big( 
  \Dt^{-1}\TNORM{\Bvh}{\Eh}^2 + \Dt\TNORM{\ROTv\Esh}{\Eh}^2
  + 2\scalEh{\Bvh}{\ROTv\Esh} +
  \NORM{\DIV\Bvh}{0,\Omega}^2
  \Big)\\[0.5em]
  &= \frac{\theta}{2}\DivNorm{\Bvh}^2
  +\frac{\theta}{2}\Dt\TNORM{\ROTv\Esh}{\Eh}^2 + \theta\scalEh{\Bvh}{\ROTv\Esh}.
  \intertext{Similarly, we transform the second term in ~\eqref{eq:ash0:split} as follows:}
  \TERM{T2} 
  &= \scalVh{\Esh}{\Esh} + \theta\scalVh{\IVh\big(\uv\times\PiEhRT\Bvh\big)}{\Esh} 
  - \theta\scalEh{\Bvh}{\ROTv\Esh}\\[0.5em]
  &= \TNORM{\Esh}{\Vh}^2 + \theta\scalVh{\IVh\big(\uv\times\PiEhRT\Bvh\big)}{\Esh} 
  - \theta\scalEh{\Bvh}{\ROTv\Esh}.
\end{align*}
Adding $\TERM{T1}$ and $\TERM{T2}$ we find that
\begin{align}
  \ashzr(\xi,\eta) 
  &= 
  \frac{\theta}{2}\DivNorm{\Bvh}^2 + \frac{\theta}{2}\Dt\TNORM{\ROTv\Esh}{\Eh}^2 +
  \TNORM{\Esh}{\Vh}^2 + \theta\scalVh{\IVh\big(\uv\times\PiEhRT\Bvh\big)}{\Esh}
  \nonumber\\[0.5em]
  &\geq \theta
  \left(
    \frac{1}{2}\DivNorm{\Bvh}^2 + \frac{1}{2}\CurlNorm{\Esh}^2 + \scalVh{\IVh\big(\uv\times\PiEhRT\Bvh\big)}{\Esh}
  \right).
  \label{eq:ash0:split:00}
\end{align}
Now, we prove that the right-hand side of~\eqref{eq:ash0:split:00} can
be bounded from below by $\TNORM{\xi}{\Xh}^2$ for a suitable choice of
$\Dt$.
Using the results of the Lemma~\ref{lemma:Vh:inner-product:estimate}
as an upper bound estimate we have
\begin{align*}
  \begin{array}{rll}
    &\hspace{-1cm}
    \scalVh{\IVh\big(\uv\times\PiEhRT\Bvh\big)}{\Esh}
    \geq -\Cst\NORM{\uv}{\infty}
    \NORM{\Bvh}{0,\Omega}\,\NORM{\Esh}{0,\Omega}           
    & \qquad\mbox{[multiply and divide by $\Dt^{\HALF}$]}
    \\[0.35em]
    &\geq -\Cst\NORM{\uv}{\infty}
    \Dt^{\HALF}\Dt^{-\HALF}\,\NORM{\Bvh}{0,\Omega}\,\NORM{\Esh}{0,\Omega}           
    & \qquad\mbox{[use Young's inequality]}
    \\[0.35em]
    &\geq -\Cst\NORM{\uv}{\infty}\Dt^{\HALF}
    \left( \frac{1}{2}\Dt^{-1}\NORM{\Bvh}{\Eh}^2+\frac{1}{2}\NORM{\Esh}{\Vh}^2\right)           
    & \qquad\mbox{[use definitions~\eqref{eq:CurlNorm:def} and~\eqref{eq:DivNorm}]}
    \\[0.35em]
    &\geq -\Cst\NORM{\uv}{\infty}
    \left(\frac{1}{2}\DivNorm{\Bvh}^2+\frac{1}{2}\CurlNorm{\Esh}^2\right),
    &
    \\[0.35em]
  \end{array}
\end{align*}
where we note that $\Cst=(\alpha^*)^{\HALF}\NORM{\IVh}{}$ is the
constant from Lemma~\ref{lemma:Vh:inner-product:estimate}.
We choose $\Dt$ sufficiently small so that
$\Cs=1-\Cst\NORM{\uv}{\infty}\Dt^{\HALF}>0$ and we write
\begin{align}
  \ashzr(\xi,\eta) 
  \geq \frac{\theta}{2}\left(1-\Cst\NORM{\uv}{\infty}\Dt^{\HALF}\right)\left(\DivNorm{\Bvh}^2 + \CurlNorm{\Esh}^2\right)
  \geq \Cs\frac{\theta}{2}\TNORM{\xi}{\Xh}^2.
  \label{eq:ash0:split:10}
\end{align}
Finally, we note that
\begin{align*}
  \TNORM{\eta_{\xi}}{\Xh}^2 
  &= 
  \DivNorm{ (\theta\slash{2})\,\big(\Bvh+\Dt \ \ROTv\Esh\big) }^2 + \CurlNorm{\Esh}^2
  \\[0.5em]
  &=
  \frac{\theta^2}{4}
  \left(
    \Dt^{-1}\TNORM{\Bvh + \Dt\ROTv\Esh}{\Eh}^2 + \NORM{\DIV\Bvh}{0,\Omega}^2
  \right) + \CurlNorm{\Esh}^2
  \\[0.5em]
  &=
  \frac{\theta^2}{4}
  \left(
    \Dt^{-1}\TNORM{\Bvh}{\Eh}^2 +
    \Dt\TNORM{\ROTv\Esh}{\Eh}^2 + 
    2\ScalEh{\Bvh}{\ROTv\Esh}  +
    \NORM{\DIV\Bvh}{0,\Omega}^2
  \right) + \CurlNorm{\Esh}^2
  \\[0.5em]
  &=
  \frac{\theta^2}{4}
  \left(
    \Dt^{-1}\TNORM{\Bvh}{\Eh}^2 \hspace{-1mm} + \hspace{-1mm}
    \NORM{\DIV\Bvh}{0,\Omega}^2 \hspace{-1mm} + \hspace{-1mm}
    2\ScalEh{\Dt^{-1/2}\Bvh}{\Dt^{1/2}\ROTv\Esh} \hspace{-1mm} + \hspace{-1mm}
    \Dt\TNORM{\ROTv\Esh}{\Eh}^2
  \right) + \CurlNorm{\Esh}^2
  \\[0.5em]
  &\leq 
  \frac{\theta^2}{4}
  \left(
    2\Dt^{-1}\TNORM{\Bvh}{\Eh}^2 + \NORM{\DIV\Bvh}{0,\Omega}^2 +
    2\Dt\TNORM{\ROTv\Esh}{\Eh}^2
  \right) +
  \CurlNorm{\Esh}^2
  \\[0.5em]
  &\leq \frac{\theta^2}{2}\DivNorm{\Bvh}^2 + \left(1+\frac{\theta^2}{2}\right)\CurlNorm{\Esh}^2
  \\[0.5em]
  &\leq \left(1+\frac{\theta^2}{2}\right)\TNORM{\xi}{\Xh}^2.
\end{align*}
The last inequality implies that
\begin{equation}
  \forall\xi\in\Xh\;\;\exists\eta\in\Xh:\quad
  \ashzr(\xi,\eta)\geq\Csb\TNORM{\xi}{\Xh}\TNORM{\eta}{\Xh},
  \qquad\Csb=\Cs\frac{\theta}{2}\left(1+\frac{\theta^2}{2}\right)^{-\HALF},
\end{equation}
from which the \emph{inf-sup} condition stated in the lemma follows
immediately.
Note that for $\Dt$ sufficiently small, we have $0<\Cs<1$.
Hence, we can just set $\Cs=1$.
\ENDPROOF

\textbf{ Proof of Theorem~\ref{theorem:wellposedness:problem:A}}.
According to Lemmas \ref{ContinuityA0} and \ref{lemma:infsup}, the
hypothesis of the Babuska-Lax-Milgram theorem are satisfied for
Problem~\ref{problem:C}.
Since problem~\ref{problem:A} and Problem~\ref{problem:C} are equivalent this will also imply the well posedness of Problem~\ref{problem:A}.


%
\section{Stability energy estimates}
\label{sec:stability:energy:estimate}
In this section we show that \eqref{eq:VarForEq4} satisfies an energy
estimates.
We begin by finding such an estimate for the continuous system
\eqref{eq:IntroSystem}.
The techniques used in the proof are, partially, laid out in
\cite{emmrich1999discrete}.
\begin{theorem}
  \label{theorem:ContEnergyEst}
  Let $\Bv$ and $\hEs$ solve \eqref{eq:VarForEq2} then
  \begin{align}
    \label{eq:FirstEnerEst}
    \frac{d}{dt}\NORM{\Bv}{0,\Omega}^2 +
    \frac{1}{2}\NORM{\sigma^{1/2}\hEs}{0,\Omega}^2
    \le 
    \CurlsigmaNorm{\Es_0}^2 +  
    \left(
      2(\sigma^*)^2\NORM{\uv}{\infty}^2
      + 1
    \right)
    \NORM{\Bv}{0,\Omega}^2,
  \end{align}
  where
  $\CurlsigmaNorm{\Es}^2=\NORM{\sigma^{1/2}\Es}{0,\Omega}^2+\NORM{\ROTv\Es}{0,\Omega}^2$. 
  As a consequence there exists a bounded function
  $\beta:[0,T]\to\mathbbm{R}^+$ such that
  \begin{align}
    \label{eq:SecondEnerEst}
    \beta(t)\NORM{\Bv(\cdot,t)}{0,\Omega}^2 +
    \frac{1}{2}\int_0^t\beta(\tau)\NORM{\sigma^{1/2}\hEs(\cdot,\tau)}{0,\Omega}^2\,d\tau
    \leq \int_0^t\beta(\tau)\CurlsigmaNorm{\Es_0(\cdot,\tau)}^2\,d\tau + \NORM{\Bvzr(\cdot,t)}{0,\Omega}^2.
  \end{align}
\end{theorem}

\BEGINPROOF
Testing equation \eqref{eq:VarForEq2Faraday} against $\wv=\Bv$,
equation \eqref{eq:VarForEq2AmpereOhm} against $\vs=\hEs(\cdot,t)$
and adding the resulting expressions we find that
\begin{equation}
  \label{eq:1}
  \frac{1}{2}\NORM{\Bv}{0,\Omega}^2 + \NORM{\sigma^{1/2}\hEs}{0,\Omega}^2
  \leq 
  -(\sigma\uv\times\Bv,\hEs) 
  -(\sigma\Es_0,\hEs)
  -(\ROTv\Es_0,\Bv).
\end{equation}
We proceed to bound the right-hand side of \eqref{eq:1} as follows
\begin{align}
  -(\sigma\hEs,\Es_0)
  &\leq
  \NORM{\sigma^{1/2}\hEs}{0,\Omega}
  \NORM{\sigma^{1/2}\Es_0}{0,\Omega}
  \leq
  \frac{1}{2}\NORM{\sigma^{1/2}\hEs}{0,\Omega}^2 + 
  \frac{1}{2}\NORM{\sigma^{1/2}\Es_0}{0,\Omega}^2,
  \label{eq:2}
  \\[0.5em]
  -(\ROTv\Es_0,\Bv)
  &\leq
  \NORM{\ROTv\Es_0}{0,\Omega}
  \NORM{\Bv}{0,\Omega}
  \leq
  \frac{1}{2}\NORM{\ROTv\Es_0}{0,\Omega}^2+
  \frac{1}{2}\NORM{\Bv}{0,\Omega}^2,
  \label{eq:3}
  \\[0.5em]
  -(\sigma\uv\times\Bv,\hEs)
  &\leq
  \NORM{\sigma^{1/2}\uv\times\Bv}{0,\Omega}
  \NORM{\sigma^{1/2}\hEs}{0,\Omega}
  \leq 
  \sigma^*
  \NORM{\uv}{\infty}^2
  \NORM{\Bv}{0,\Omega}^2 +
  \frac{1}{4}\NORM{\sigma^{1/2}\hEs}{0,\Omega}^2,
  \label{eq:4}
\end{align}
Estimate \eqref{eq:FirstEnerEst} follows from \eqref{eq:1},
\eqref{eq:2}, \eqref{eq:3} and \eqref{eq:4}.
To prove \eqref{eq:SecondEnerEst} we define
\begin{equation}
  \beta(t) = 
  \exp\left( 
    -\int_0^{t}
    \big(
      2\NORM{\uv}{\infty}^2(\sigma^*)^2+1
    \big)\,d\tau
  \right).
\end{equation}
Multiplication by $\beta$ in \eqref{eq:FirstEnerEst} yields
\begin{equation}
  \frac{d}{dt}\left( 
    \beta\NORM{\Bv}{0,\Omega}^2
  \right) +
  \frac{\beta}{2}\NORM{\sigma^{1/2}\hEs}{0,\Omega}^2
  \leq \beta\CurlsigmaNorm{\Es_0}.
\end{equation}
Integration in time gives \eqref{eq:SecondEnerEst}.
\ENDPROOF

Next Theorem mimics the continuous Theorem~\ref{theorem:ContEnergyEst} in the discrete settings.
\begin{theorem}
  \label{theorem:DiscreteEnergyEstimate}~\\
  \begin{description}
  \item[$(\mathrm{i})$] 
    Let $\theta\in\left[ 0,1\right]$.
    The solution of Scheme \eqref{eq:VarForEq4} satisfies
    \begin{align}
      &\left(\theta-\frac{1}{2}\right)
      \frac{\TNORM{\Bvh^{n+1}-\Bvh^{n}}{\Eh}^2}{\Dt} +
      \frac{ \TNORM{\Bvh^{n+1}}{\Eh}^2 - \TNORM{\Bvh^n}{\Eh}^2}{\Dt} +
      \frac{1}{2}\TNORM{\hEsh^{n+\theta}}{\Vh}^2
      \nonumber\\[0.5em]
      &\qquad\qquad\leq
      \TNORM{\IVh\Es_0^{n+\theta}}{\HROTv{\Omega}}^2 +
      \frac{1}{2}\left( 
        1 + 4\Cst
        \NORM{\uv}{\infty}^2
      \right)\,
      \left(
        \theta    \TNORM{\Bvh^{n+1}}{\Eh}^2+
        (1-\theta)\TNORM{\Bvh^{n}  }{\Eh}^2
      \right),
      \label{eq:DiscreteEstimateone}
    \end{align}
    where
    $\TNORM{\IVh\Es_0^{n+\theta}}{\HROTv{\Omega}}^2=\TNORM{\IVh\Es_0^{n+\theta}}{\Vh}^2+\TNORM{\ROTv\IVh\Es_0^{n+\theta}}{\Eh}^2$,
    and we recall that $\Cst$ is the constant introduced in
    Lemma~\ref{lemma:Vh:inner-product:estimate}.

    \medskip
  \item[$(\mathrm{ii})$] If $\theta\in\left[\frac{1}{2},1\right]$, then we can conclude that
      \begin{align}
      \hspace{-2mm}
      (\beta)^{n+1}\TNORM{\Bvh^{n+1}}{\Eh}^2 
      \hspace{-.5mm} + \hspace{-.5mm}
      \frac{\gamma\Dt}{2}\sum_{\ell=0}^{n}\beta^{n+1-\ell}\TNORM{\hEsh^{n-\ell+\theta}}{\Vh}^2 
      \leq
      \TNORM{\Bvh^{0}}{\Eh}^2 
      \hspace{-.5mm} + \hspace{-.5mm}
      \gamma\Dt\sum_{\ell=0}^{n}\beta^{n+1-\ell}
      \TNORM{\IVh\Es_0^{n-\ell+\theta}}{\HROTv{\Omega}}^2,
      \label{eq:DiscEstTwo}
    \end{align}
    where
        \begin{align}
      \beta = \frac
      {\big(1-Q\theta\big)}
      {\big(1+Q(1-\theta)\big)},
      \qquad
      \gamma = \frac{1}{\big(1-Q\theta\big)}
      \qquad\textrm{and}~~
      Q=\Dt\big(1+4\Cst\NORM{\uv}{\infty}^2\big).
      \label{eq:beta-Q:def}
    \end{align}
    The coefficients in \eqref{eq:DiscEstTwo} are guaranteed to be positive when
    \begin{equation}
      \Dt<\frac{1}{\theta\big(1+4\Cst\NORM{\uv}{\infty}^2\big)},
      \label{eq:Assumption}
    \end{equation}
    making \eqref{eq:DiscEstTwo} an energy estimate.
  \end{description}
\end{theorem}

\BEGINPROOF
$(\mathrm{i})$.  Testing equation \eqref{eq:VarForEq4Faraday} against $\wvh
=\Bvh^{n+\theta}=\theta\Bvh^{n+1}+(1-\theta)\Bvh^n$ and equation
\eqref{eq:VarForEq4AmpereOhm} against $\vsh = \hEsh^{n+\theta}$ and
adding them together we arrive at
\begin{align}
  &\scallEh{\frac{\Bvh^{n+1}-\Bvh^n}{\Dt}}{\Bvh^{n+\theta}} +
  \TNORM{\hEsh^{n+\theta}}{\Vh}^2
  \nonumber\\[0.5em] 
  &=
  -\scalEh{\ROTv\IVh\Es_0^{n+\theta}}{\Bvh^{n+\theta}}
  -\scalVh{\IVh\Es_0^{n+\theta}}{\hEsh^{n+\theta}}
  -\scalVh{\IVh\big(\uv\times\PiEhRT\Bvh^{n+\theta}\big)}{\hEsh^{n+\theta}}
  \nonumber\\[0.5em] 
  &= \TERM{T1} + \TERM{T2} + \TERM{T3}.
  \label{eq:deq1}
\end{align}
We transform the first term of the left-hand side of \eqref{eq:deq1}
using the identity
\begin{align}
  \Bvh^{n+\theta} = 
  \Dt\left(\theta-\frac{1}{2}\right)\frac{\Bvh^{n+1}-\Bvh^{n}}{\Dt} +
  \frac{\Bvh^{n+1}+\Bvh^n}{2}.
  \label{eq:Neweq1}
\end{align}
We obtain:
\begin{align}
  \scallEh{\frac{\Bvh^{n+1}-\Bvh^n}{\Dt}}{\Bvh^{n+\theta}} 
  &=
  \Dt\left(\theta-\frac{1}{2}\right)
  \scallEh{ \frac{\Bvh^{n+1}-\Bvh^n}{\Dt} }{ \frac{\Bvh^{n+1}-\Bvh^n}{\Dt} } 
  \nonumber\\[0.5em]
  &\quad+
  \scallEh{ \frac{\Bvh^{n+1}-\Bvh^n}{\Dt} }{ \frac{\Bvh^{n+1}+\Bvh^n}{2}   } 
  \nonumber\\[0.5em]
  &=
  \Dt\left(\theta-\frac{1}{2}\right)
  \frac{ \TNORM{\Bvh^{n+1}-\Bvh^{n}}{\Eh}^2 }{\Dt^2} +
  \frac{ \TNORM{\Bvh^{n+1}}{\Eh}^2-\TNORM{\Bvh^{n}}{\Eh}^2 }{2\Dt}.
  \label{eq:deq2}
\end{align}
Next, we bound the three terms in the right-hand side of
\eqref{eq:deq1} by using the Young inequality with parameters
$\epsilon_1$, $\epsilon_2$, and $\epsilon_1$.
For the first two terms we obtain the estimates:
\begin{align}
  \TERM{T1} 
  &\leq
  \frac{\epsilon_1}{2}\TNORM{\ROTv\IVh\Es_0^{n+\theta}}{\Eh}^2 + \frac{1}{2\epsilon_1}\TNORM{\Bvh^{n+\theta}}{\Eh}^2
  \nonumber\\[0.5em] 
  &\leq 
  \frac{\epsilon_1}{2}
  \TNORM{\ROTv\IVh\Es_0^{n+\theta}}{\Eh}^2 +
  \frac{1}{2\epsilon_1}
  \left(
    \theta^2    \TNORM{\Bvh^{n+1}}{\Eh}^2+
    (1-\theta)^2\TNORM{\Bvh^{n}  }{\Eh}^2
  \right),
  \nonumber\\[0.5em] 
  &\leq 
  \frac{\epsilon_1}{2}
  \TNORM{\ROTv\IVh\Es_0^{n+\theta}}{\Eh}^2 +
  \frac{1}{\epsilon_1}
  \left(
    \theta    \TNORM{\Bvh^{n+1}}{\Eh}^2+
    (1-\theta)\TNORM{\Bvh^{n}  }{\Eh}^2
  \right),
  \label{eq:deq4}
  \\ 
  \TERM{T2}
  &\leq 
  \frac{\epsilon_2}{2}  \TNORM{\IVh\Es_0^{n+\theta}}{\Vh}^2 +
  \frac{1}{2\epsilon_2} \TNORM{\hEsh^{n+\theta}}{\Vh}^2.
  \label{eq:deq3}
\end{align}

The bound for the third term requires a bit more work.
Since $\theta\in[0,1]$, we note that $\theta^2\leq\theta$ and
$(1-\theta)^2\leq1-\theta$. 
Therefore we have an estimate
\begin{align*}
  \TNORM{\IVh(\uv\times\PiEhRT\Bvh^{n+\theta})}{\Vh}^2
  &\leq \Cs\NORM{\uv}{\infty}^2\TNORM{\theta\Bvh^{n+1}+(1-\theta)\Bvh^{n}}{\Eh}^2
  \nonumber\\
  &\leq 2\Cs\NORM{\uv}{\infty}^2
  \left(
    \theta^2    \TNORM{\Bvh^{n+1}}{\Eh}^2 +
    (1-\theta)^2\TNORM{\Bvh^{n}  }{\Eh}^2
  \right)
  \nonumber\\
  &\leq 2\Cs\NORM{\uv}{\infty}^2
  \left(
    \theta    \TNORM{\Bvh^{n+1}}{\Eh}^2 +
    (1-\theta)\TNORM{\Bvh^{n}  }{\Eh}^2
  \right).
\end{align*}
Next we again use the Young's inequality
\begin{align}
  \TERM{T3}
  &\leq
  \frac{\epsilon_3}{2} \TNORM{\IVh(\uv\times\theta\Bvh^{n+\theta})}{\Vh}^2+
  \frac{1}{2\epsilon_3}\TNORM{\hEsh^{n+\theta}}{\Vh}^2 +
  \nonumber\\[0.5em]
  &\leq
  \Cs\epsilon_3\NORM{\uv}{\infty}^2
  \left(
    \theta    \TNORM{\Bvh^{n+1}}{\Eh}^2 +
    (1-\theta)\TNORM{\Bvh^{n}  }{\Eh}^2
  \right) +
  \frac{1}{2\epsilon_3}\TNORM{\hEsh^{n+\theta}}{\Vh}^2.
  \label{eq:deq5}  
\end{align}
Setting $\epsilon_1=\epsilon_2=\epsilon_3=2$, combining
\eqref{eq:deq2} with the estimates of $\TERM{T1}$, $\TERM{T2}$, and
$\TERM{T3}$, and finally noting that
$\TNORM{\IVh\Es_0^{n+\theta}}{\HROTv{\Omega}}^2=\TNORM{\IVh\Es_0^{n+\theta}}{\Vh}^2+\TNORM{\ROTv\IVh\Es_0^{n+\theta}}{\Eh}^2$
yield~\eqref{eq:DiscreteEstimateone}, 
which is the first assertion of the theorem.

\medskip 
\noindent
$(\mathrm{ii})$.  If $\theta\in [1/2,1]$, the coefficient in the first term on
the left hand side of \eqref{eq:DiscreteEstimateone} is positive and
we can write
\begin{align} 
  \TNORM{\Bvh^{n+1}}{\Eh}^2 -
  \TNORM{\Bvh^n}{\Eh}^2
  &\leq
  \Dt\left(
    -\frac{1}{2}\TNORM{\hEsh^{n+\theta}}{\Vh}^2 
    +\TNORM{\IVh\Es_0^{n+\theta}}{\HROTv{\Omega}}^2
  \right)
  \nonumber\\[0.5em] 
  &\qquad
  +\Dt
  \big(1+4\Cst\NORM{\uv}{\infty}^2\big)\,
  \left(
    \theta    \TNORM{\Bvh^{n+1}}{\Eh}^2+
    (1-\theta)\TNORM{\Bvh^{n}  }{\Eh}^2
  \right).
  \label{eq:DiscreteEstimateonestar}
\end{align}
To simplify the notation, let
$Q=\Dt\big(1+4\Cst\NORM{\uv}{\infty}^2\big)$ and
\begin{align*}
  \calF^{n+\theta}(\hEsh,\Es_0) = \Dt\left(
    -\frac{1}{2}\TNORM{\hEsh^{n+\theta}}{\Vh}^2
    +\TNORM{\IVh\Es_0^{n+\theta}}{\HROTv{\Omega}}^2
  \right)
\end{align*}
Rearranging the terms and dividing by $\big(1-Q\theta\big)$ we find:
\begin{align}
  \TNORM{\Bvh^{n+1}}{\Eh}^2 - 
  \frac{\big(1+Q(1-\theta)\big)}{\big(1-Q\theta\big)}\TNORM{\Bvh^{n} }{\Eh}^2 \leq
  \frac{1}{\big(1-Q\theta\big)}\calF(\hEsh,\Eszr)^{n+\theta}.
  \label{eq:energy:estimate:2:10}
\end{align}
Now, we introduce the quantities
\begin{align*}
  \alpha = \frac
  {\big(1+Q(1-\theta)\big)}
  {\big(1-Q\theta\big)},
  \qquad
  \gamma = \frac{1}{\big(1-Q\theta\big)},
\end{align*}
and note that quantity $\alpha$ is well defined and strictly positive
since Assumption~\eqref{eq:Assumption} guarantees that $1-Q\theta>0$,
and $Q>0$ implies $\big(1+Q(1-(1-\theta)\big)\leq1$ for
$\theta\in[0,1]$, so that $\alpha>0$.
We rewrite~\eqref{eq:energy:estimate:2:10} as
\begin{align*}
  \TNORM{\Bvh^{n+1}}{\Eh}^2 - 
  \alpha\TNORM{\Bvh^{n} }{\Eh}^2 \leq
  \gamma\calF^{n+\theta}(\hEsh,\Eszr).
\end{align*}
Such inequality must be true for any index $n\geq0$.
We express this fact by keeping $n$ fixed and introducing the index
$\ell=0,\ldots,n$ such that
\begin{align*}
  \TNORM{\Bvh^{n+1-\ell}}{\Eh}^2 - 
  \alpha\TNORM{\Bvh^{n-\ell} }{\Eh}^2 \leq
  \gamma\calF^{n-\ell+\theta}(\hEsh,\Eszr).
  \label{eq:energy:estimate:2:10}
\end{align*}
Then, we multiply by $\alpha^{\ell}$ and adding all the resulting
inequalities we find a telescopic sum where all intermediate terms 
like $\Bvh^{n-\ell}$ cancel.
We illustrate this fact by writing the first four inequalities for
$\ell=0,\ldots,3$:
\begin{align*}
  \begin{array}{rllll}
    \mbox{for $\ell=0$:} & \qquad\TNORM{\Bvh^{n+1}}{\Eh}^2 & -\alpha\TNORM{\Bvh^{n}   }{\Eh}^2 &\leq \gamma\calF^{n+\theta}(\hEsh,\Eszr)   & \quad\mbox{\big[multiply by $1$      \big]},\\[0.5em]
    \mbox{for $\ell=1$:} & \qquad\TNORM{\Bvh^{n}  }{\Eh}^2 & -\alpha\TNORM{\Bvh^{n-1} }{\Eh}^2 &\leq \gamma\calF^{n-1+\theta}(\hEsh,\Eszr) & \quad\mbox{\big[multiply by $\alpha$ \big]},\\[0.5em]
    \mbox{for $\ell=2$:} & \qquad\TNORM{\Bvh^{n-1}}{\Eh}^2 & -\alpha\TNORM{\Bvh^{n-2} }{\Eh}^2 &\leq \gamma\calF^{n-2+\theta}(\hEsh,\Eszr) & \quad\mbox{\big[multiply by $\alpha^2$\big]},\\[0.5em]
    \mbox{for $\ell=3$:} & \qquad\TNORM{\Bvh^{n-2}}{\Eh}^2 & -\alpha\TNORM{\Bvh^{n-3} }{\Eh}^2 &\leq \gamma\calF^{n-3+\theta}(\hEsh,\Eszr) & \quad\mbox{\big[multiply by $\alpha^3$\big]},\\[0.5em]
    \ldots & \ldots
  \end{array}
\end{align*}
The sum of these expressions (with coefficients indicated on the right) gives:
\begin{align*}
  \TNORM{\Bvh^{n+1}}{\Eh}^2 - \alpha^4\TNORM{\Bvh^{n-3} }{\Eh}^2\leq \gamma\sum_{\ell=0}^{3}\alpha^{\ell}\calF^{n-\ell+\theta}(\hEsh,\Eszr).
\end{align*}
Adding all inequalities for $\ell=0,\ldots,n$ yields
\begin{align*}
  \TNORM{\Bvh^{n+1}}{\Eh}^2 - 
  \alpha^{n+1}\TNORM{\Bvh^{0} }{\Eh}^2 \leq
  \gamma\sum_{\ell=0}^{n}\alpha^{\ell}\calF^{n-\ell+\theta}(\hEsh,\Eszr).
\end{align*}
Finally, we substitute back the expression for $\calF$ and $\gamma$,
multiply both side of~\eqref{eq:beta-Q:def} by
$\beta^{n+1}=\alpha^{-(n+1)}$, rearrange the terms and obtain the
second assertion of the theorem.
\ENDPROOF

\begin{remark}
   Theorem \ref{theorem:DiscreteEnergyEstimate} above gives sufficient conditions for energy stability,
   but condition \eqref{eq:Assumption} is by no means necessary. 
   Numerical experimentation shows that for $\theta\in[1/2,1]$ the
   method is unconditionally stable.
 \end{remark}


\section{Numerical experiments}
\label{Sec:NumericalExperiments}

In this section we will present the results of  a series of numerical experiments that sheds some light on the performance of the VEM developed and analyzed throughout this article. It is divided in three sections, the first on explores the rate of convergence and the divergence preserving nature of the numerical method. The second section studies the energy estimate that was introduced in theorem \ref{theorem:DiscreteEnergyEstimate}. In the final section we introduce the Hartmann problem and use this novel discretization to approximate its solution.


  




\begin{figure}[ht!]
  \begin{center} 
    \includegraphics[width=.30\textwidth]{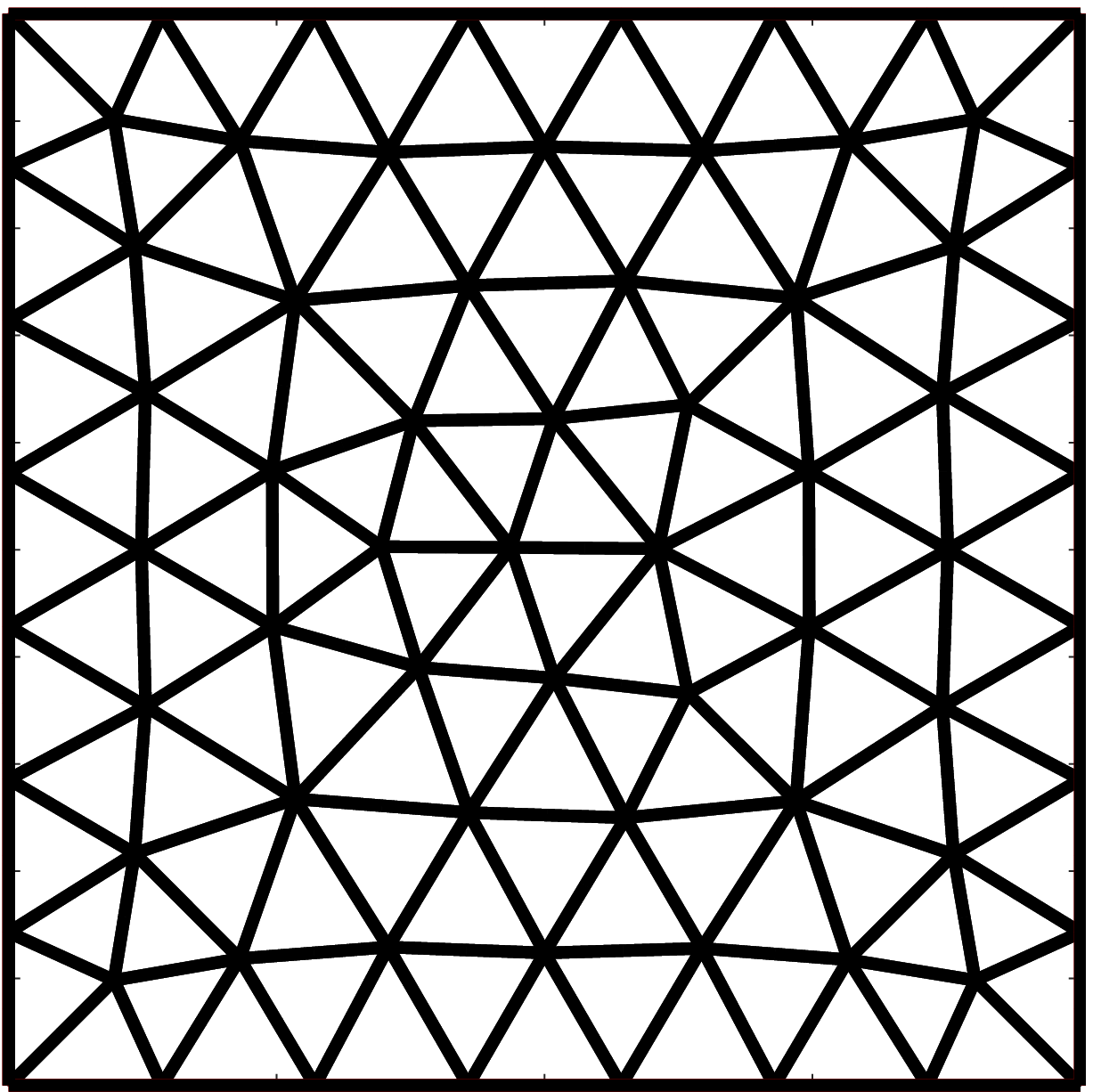} \hspace{2mm}
    \includegraphics[width=.30\textwidth]{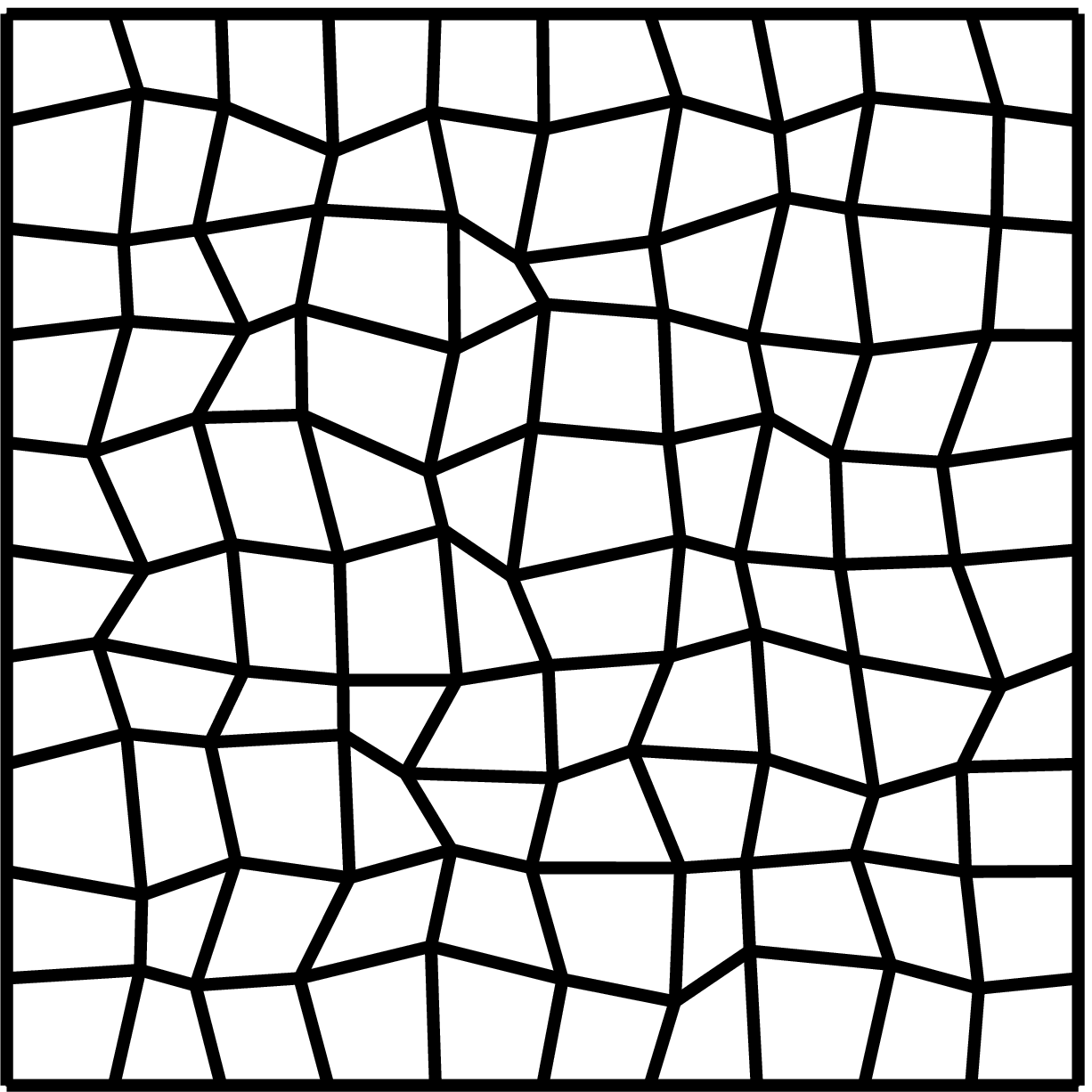}\hspace{2mm}
    \includegraphics[width=.30\textwidth]{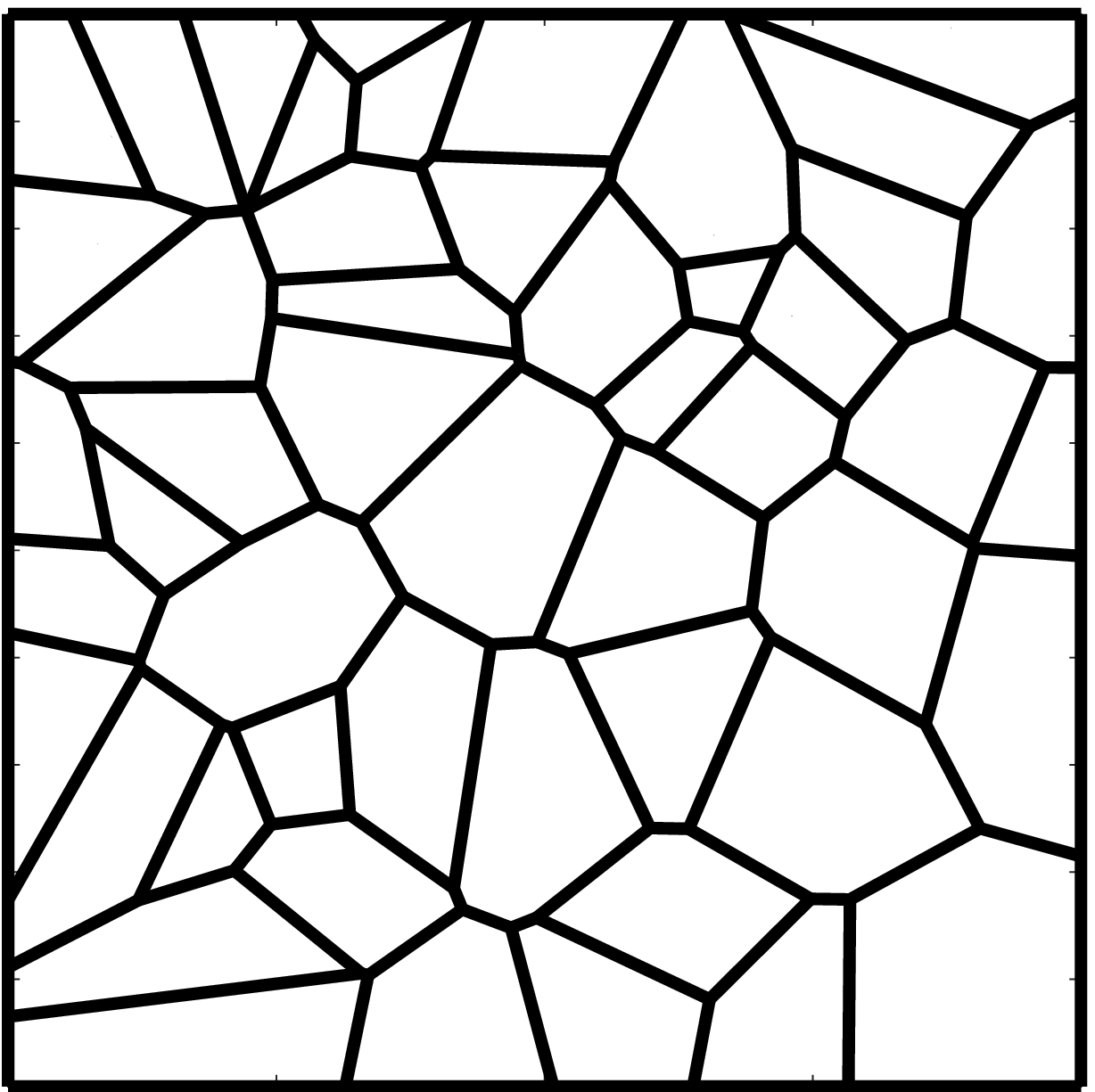}
    \caption{Illustration of the meshes used for testing the rate of
      convergence: triangular mesh (left panel), perturbed square mesh
      (central panel) and Voronoi tesselation (right panel).}
    \label{fig:Meshes}
  \end{center}
\end{figure}

\subsection{Experimental analysis of the rate of convergence and the divergence free condition.}

To assess the performance of the VEM we study the numerical approximations of
Problem~\ref{eq:VarForEq2} on a square domain $\Omega=[-1,1]^2$.
We consider the velocity field $\uv=(\us_x,\us_y)^T$ given by
\begin{align}
  \us_x(x,y)&=-\frac{(x^2+y^2-1)(\sin(xy)+\cos(xy))-100e^{x}+100e^{y}}{2(50e^x-y\sin(xy)+y\cos(xy))},\\[0.25em]
  \us_y(x,y)&=\frac{(x^2+y^2-1)(\sin(xy)+\cos(xy))-100e^{x}+100e^{y}}{2(50e^y+x\sin(xy)-x\cos(xy))}
\end{align}
and the initial and the boundary conditions are set in accordance with the exact solution of the electric and the omagnetic fields:
\begin{align}
  \Bv(x,y,t) &= 
  \begin{pmatrix}
    50e^y+x\sin(xy)-x\cos(xy)\\[0.25em]
    50 e^x-y\sin(xy)+y\cos(xy)
  \end{pmatrix}e^{-t},\\[0.5em]
  \Es(x,y,t) &= -\big( 50(e^x-e^y)+\cos(xy)+\sin(xy) \big) e^{-t}.
\end{align}
To check the robustness of the method we have selected three different
mesh families, including triangular meshes, randomly perturbed square
meshes, and meshes based on Voronoi tessellations.
An example of each mesh family is shown in Figure~\ref{fig:Meshes}.

\begin{figure}[H]
  \begin{center}
    \begin{tabular}{cc}
      \begin{overpic}[width=.475\textwidth]{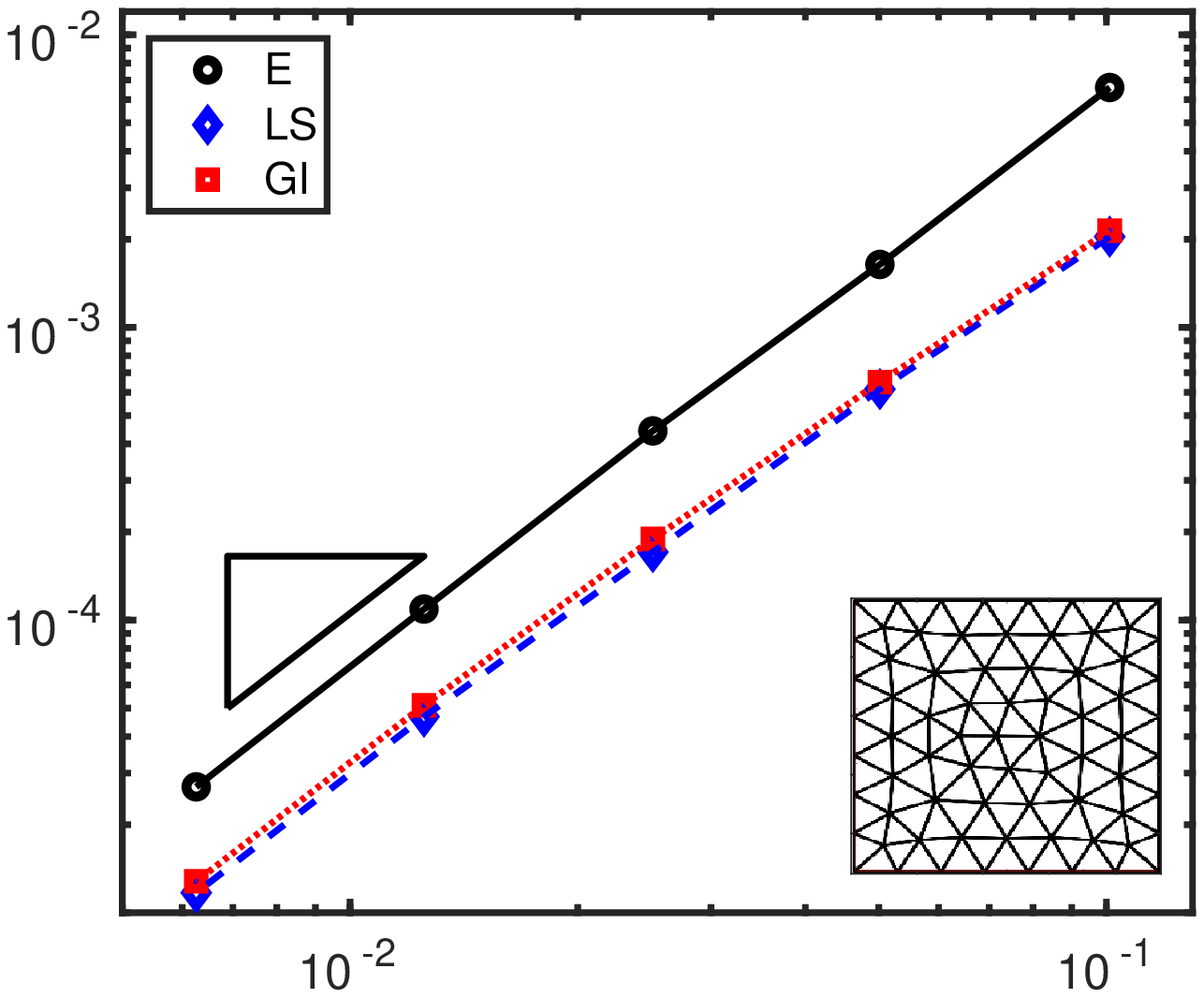}
       \put( 0,8.5){\begin{sideways}\textbf{Electric field relative error}\end{sideways}}
        \put(40,-2){\textbf{Mesh size $\mathbf{h}$}}
        \put(17,29){\textbf{2}}
        \put(26,36){\textbf{1}}
      \end{overpic}
      & \qquad
      \begin{overpic}[width=.475\textwidth]{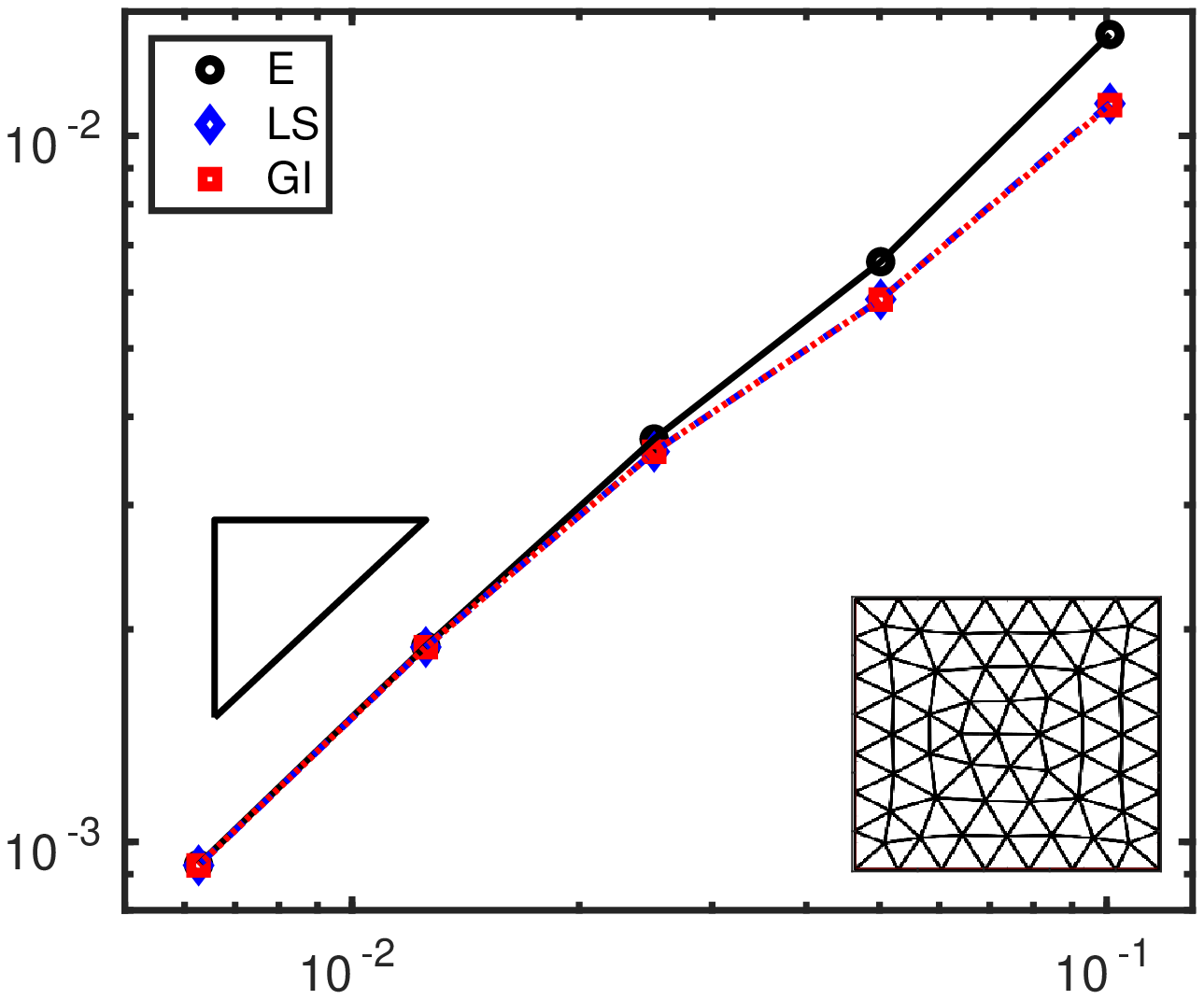}
       \put( 0,7.5){\begin{sideways}\textbf{Magnetic field relative error}\end{sideways}}
        \put(40,-2) {\textbf{Mesh size $\mathbf{h}$}}
        \put(17,29){\textbf{1}}
        \put(25,38){\textbf{1}}
      \end{overpic}
      \vspace{-5mm} \\[1.5em] 
      \begin{overpic}[width=.475\textwidth]{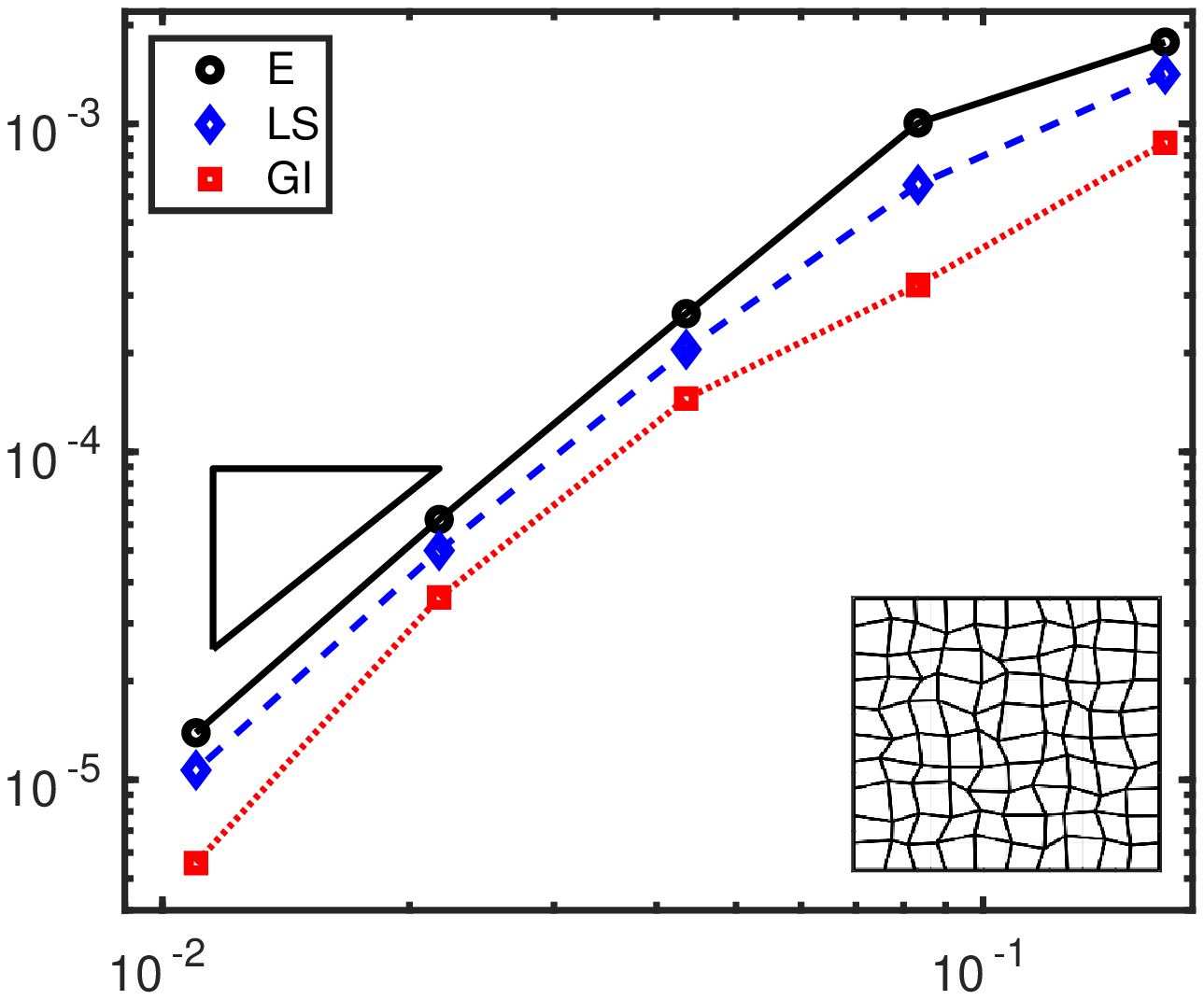}
       \put( 0,8.5){\begin{sideways}\textbf{Electric field relative error}\end{sideways}}
        \put(40,-2) {\textbf{Mesh size $\mathbf{h}$}}
        \put(16,34){\textbf{2}}
        \put(25,42){\textbf{1}}
      \end{overpic}
      & \qquad
      \begin{overpic}[width=.475\textwidth]{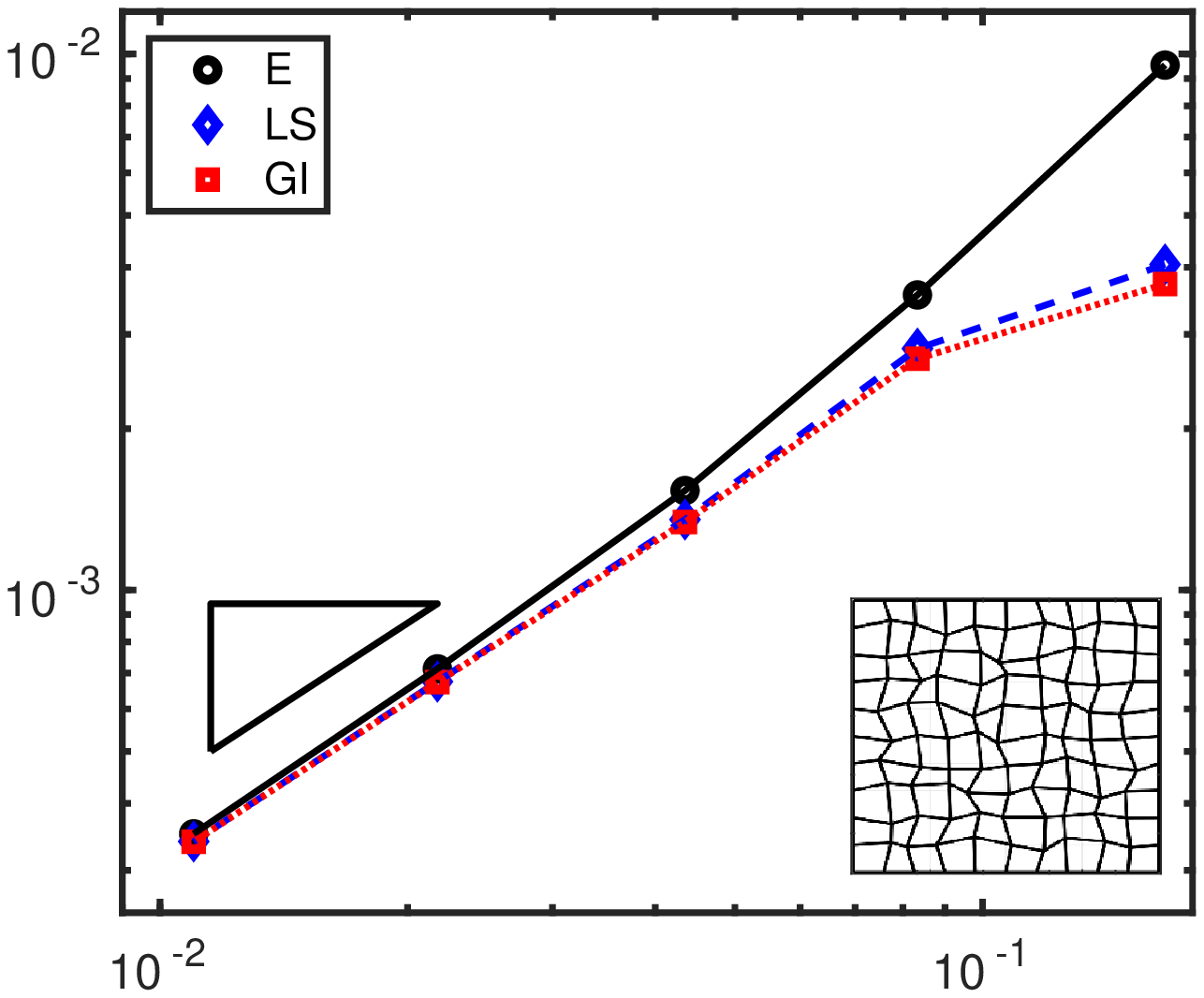}
       \put( 0,7.5){\begin{sideways}\textbf{Magnetic field relative error}\end{sideways}}
        \put(45,-2) {\textbf{Mesh size $\mathbf{h}$}}
        \put(16,25){\textbf{1}}
        \put(25,32){\textbf{1}}
      \end{overpic}
      \vspace{-5mm} \\[1.5em] 
      \begin{overpic}[width=.475\textwidth]{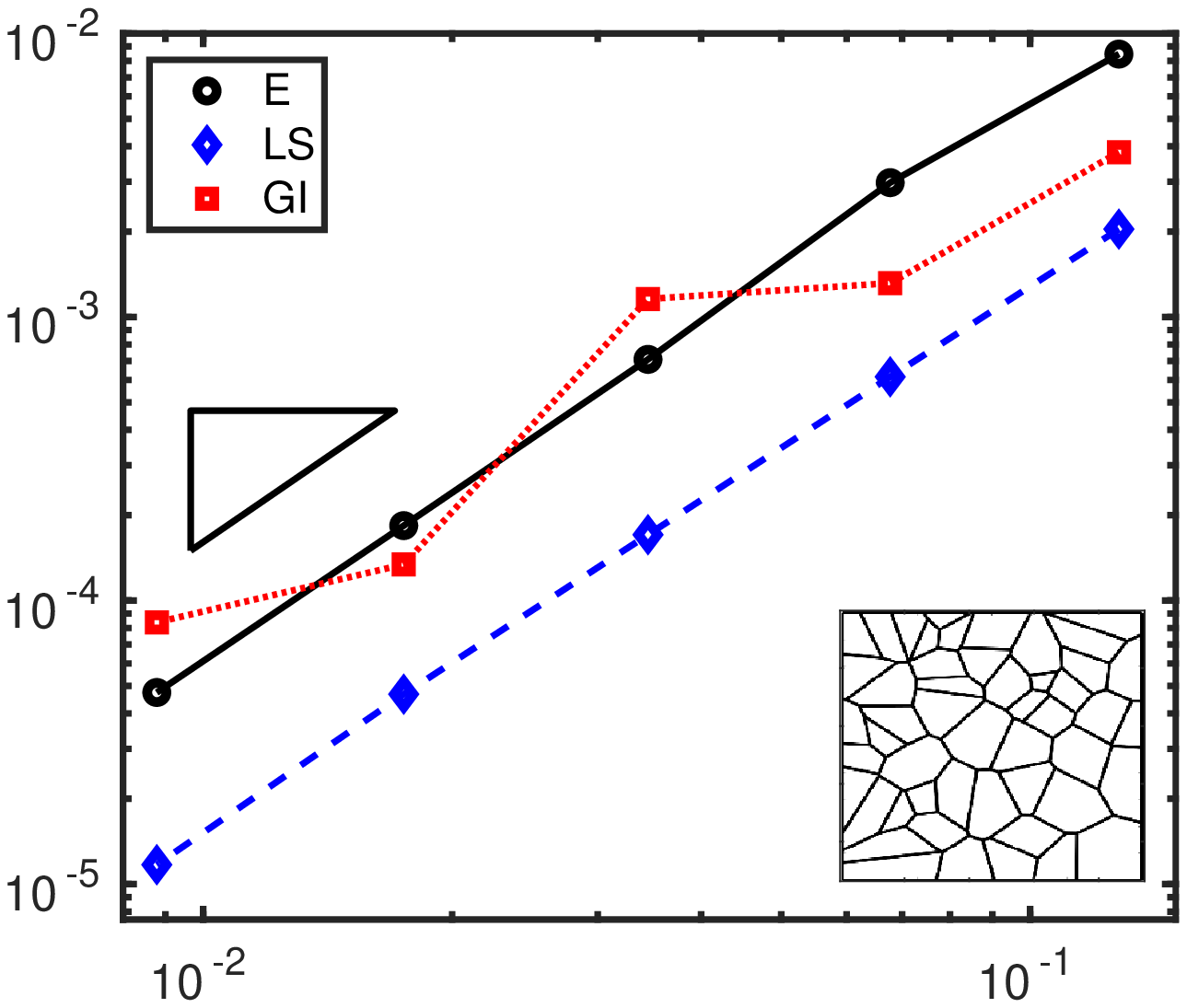}
       \put( 0,8.5){\begin{sideways}\textbf{Electric field relative error}\end{sideways}}
        \put(40,-2) {\textbf{Mesh size $\mathbf{h}$}}
        \put(15,40){\textbf{2}}
        \put(23,47){\textbf{1}}
      \end{overpic}
      & \qquad
      \begin{overpic}[width=.475\textwidth]{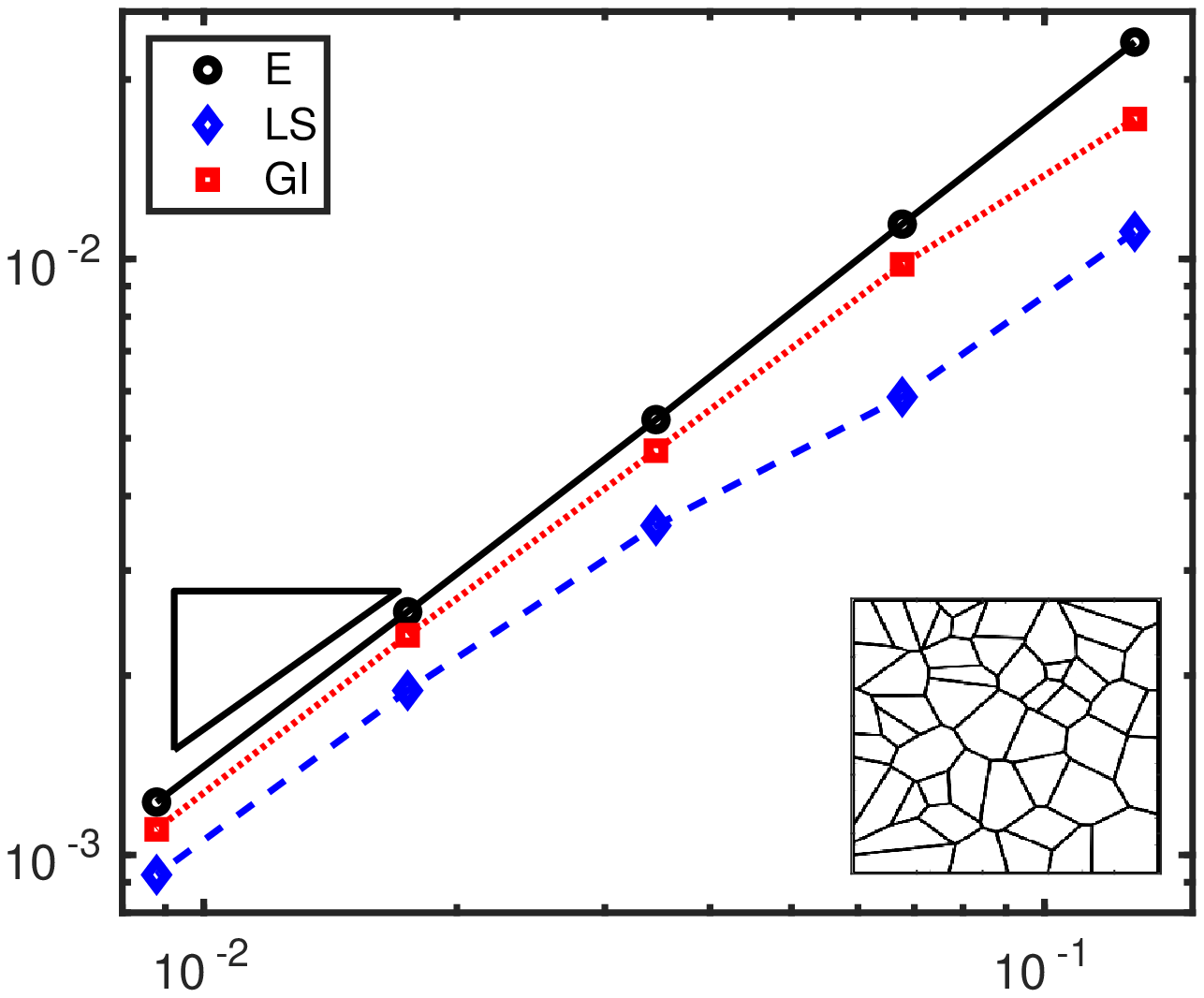}
       \put( 0,7.5){\begin{sideways}\textbf{Magnetic field relative error}\end{sideways}}
        \put(40,-2) {\textbf{Mesh size $\mathbf{h}$}}
        \put(14,25){\textbf{1}}
        \put(23,33){\textbf{1}}
      \end{overpic}
    \end{tabular}
  \end{center}
  
  \medskip
 \caption{ Error curves for the virtual element approximation of the electric and magnetic field (respectively, left and right panels)
    for the three mesh families of Figure~\ref{fig:Meshes}: triangular mesh family (top), quadrilateral mesh family (middle), Voronoi mesh family (bottom).
    The convergence rate is reflected by the slope of the curves in
    the log-log plots; the reference convergence rate is shown by the
    triangle in each plot. The symbols E,LS,GI refers to the three alternatives we have for constructing the nodal mass matrix; the elliptic projector (E), least squares projector (LS) and the Galerkin interpolator (GI), respectively.}
  \label{fig:test1:convergence-curves}
\end{figure}

\begin{figure}[H]
  \begin{center}
    \begin{tabular}{ccc}
      \begin{overpic}[width=.3\textwidth]{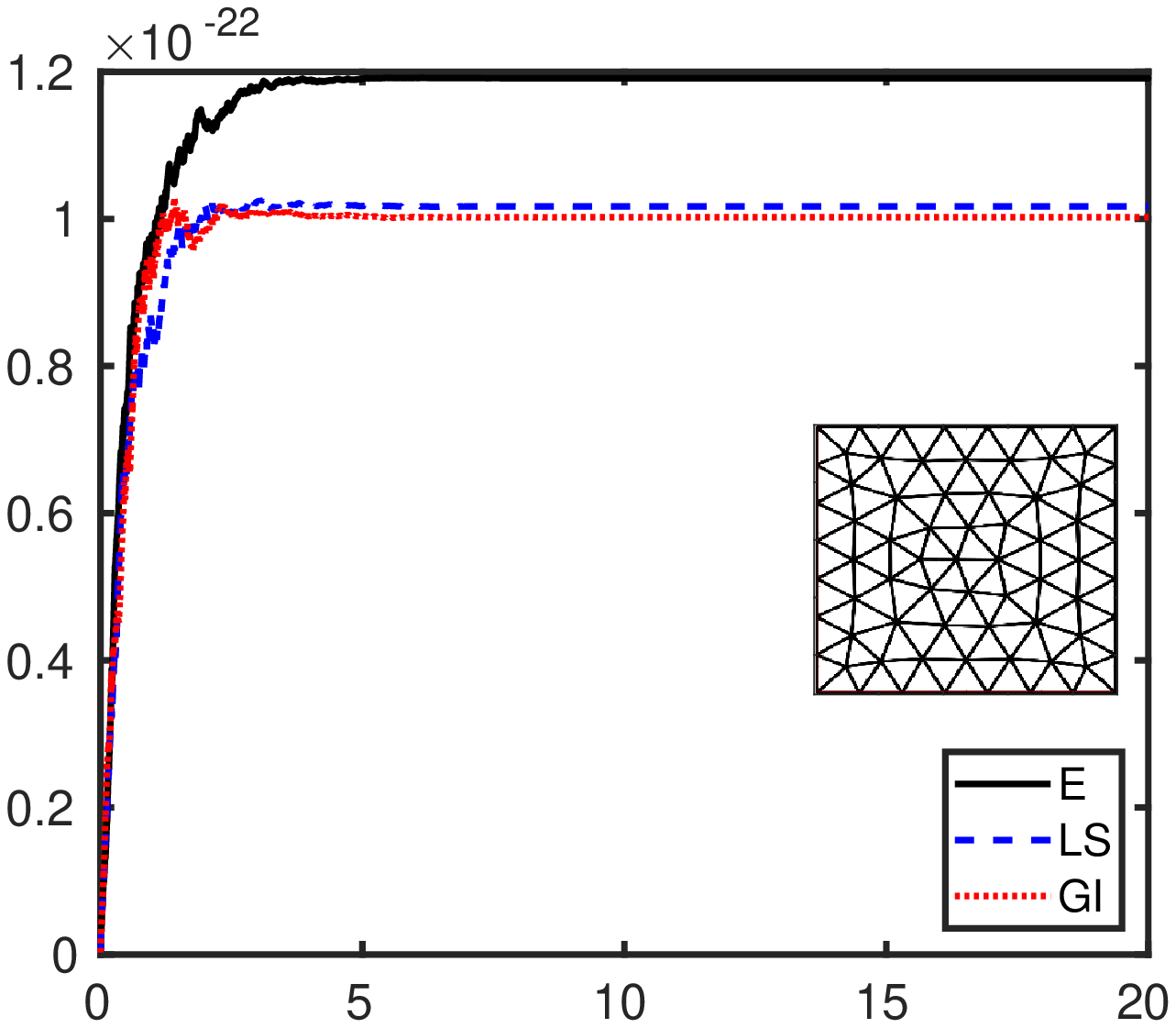}
        \put(45,-2){{\small\textbf{Time}}}
        \put(-6,25){\begin{sideways}{\small $\TNORM{\DIV\Bvh}{\Ph}^2$}\end{sideways}}
      \end{overpic}
      & \qquad
      \begin{overpic}[width=.3\textwidth]{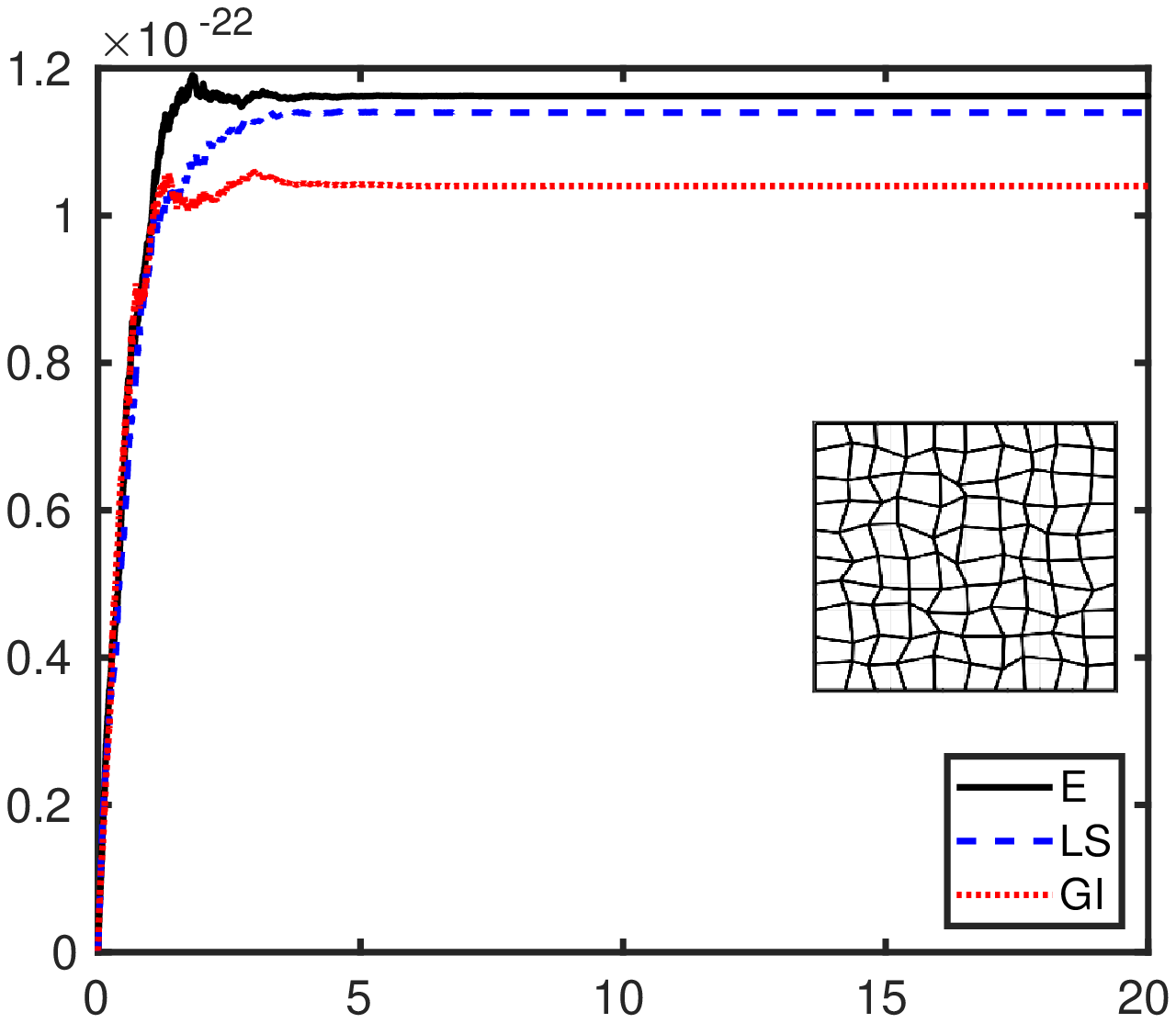}
        \put(45,-2){{\small\textbf{Time}}}
        \put(-6,25){\begin{sideways}{\small $\TNORM{\DIV\Bvh}{\Ph}^2$}\end{sideways}}
      \end{overpic}
      &\qquad
      \begin{overpic}[width=.3\textwidth]{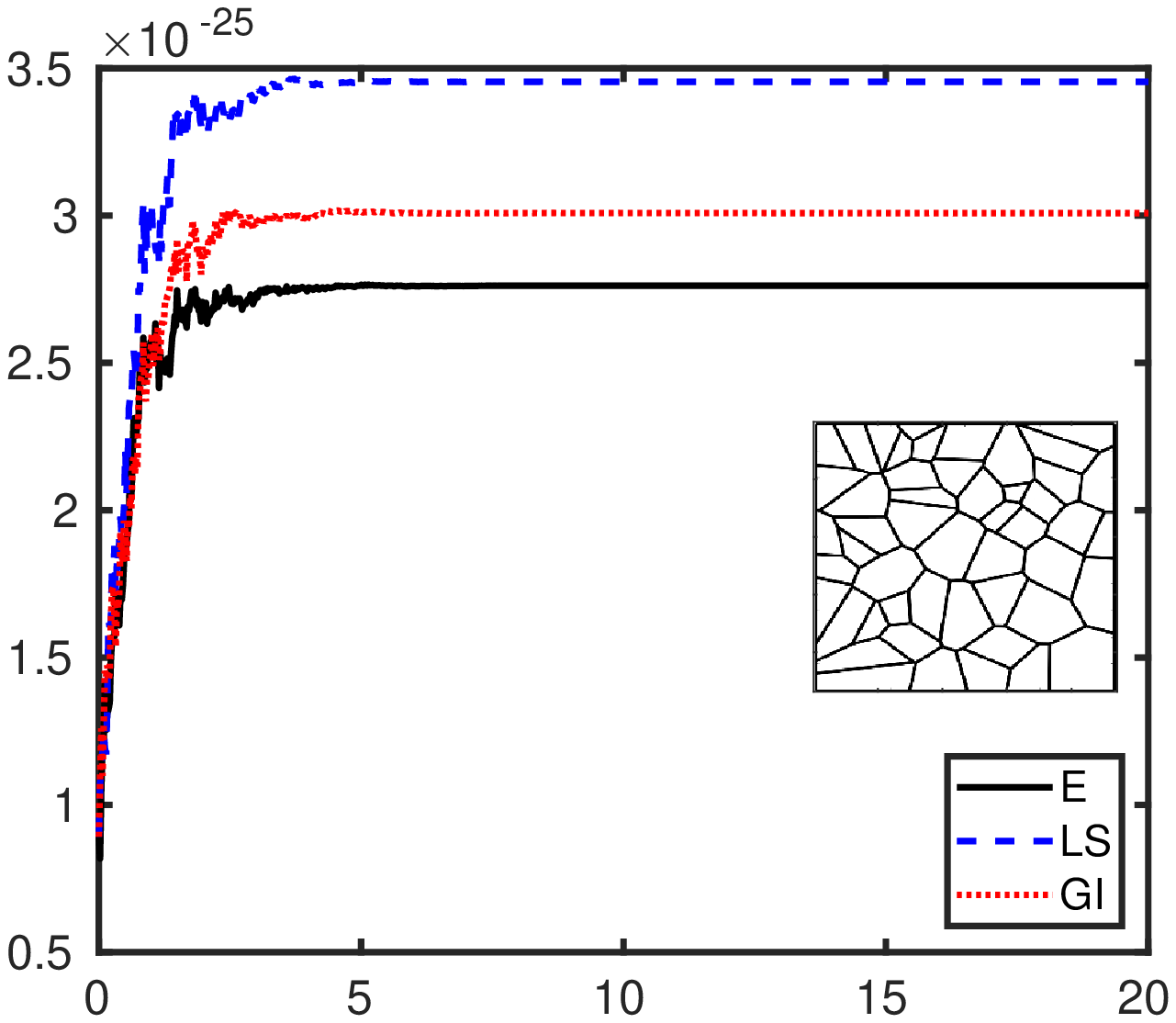}
        \put(45,-2){{\small\textbf{Time}}}
        \put(-6,25){\begin{sideways}{\small $\TNORM{\DIV\Bvh}{\Ph}^2$}\end{sideways}}
      \end{overpic}
     \end{tabular}
  \end{center}
  \caption{Plots of the time evolution of the square of the $L^2$ norm of the divergence of the numerical magnetic field on the three mesh families of Figure~\ref{fig:Meshes}. (Left) triangular mesh, (Middle) quadrilateral cells and (Right) Voronoi tesselation.}
\label{fig:DivPlots}
\end{figure}

The time marching scheme uses $\theta =1/2$. 
Errors with different values of $\theta$ are very similar and we therefore omit them.  
The final time is set at $T=0.25$ and the time step follows the
assignment $\Dt=0.05\hh^2$.
Figure~\ref{fig:test1:convergence-curves} shows the log-log plots of
the error curves for the approximation of the electric and magnetic
fields.
The errors are relative and measured in the $\LTWO$ norms, this is to say they are the $\LTWO$ norm of the difference between numerical and exact solutions divided by the norm of the exact solution.
\newpage
An important feature of the VEM that we have presented is that the magnetic field remains divergence free throughout the simulations. 
Next, we will present the results of numerical experiments aimed at gathering experimental evidence to support our theoretical findings.
In Figure~\ref{fig:DivPlots} we present three simulations, each done in a different type of mesh, the $y-$axis represents the squared $L^2$ norm of the magnetic field.
\subsection{Experimental analysis of the energy estimates}

This section is dedicated to an experimental study of the energy estimate presented in 
Theorem \ref{theorem:DiscreteEnergyEstimate}.
For this purpose we define a normalized version of the right hand side and the left hand side of \eqref{eq:DiscEstTwo} and their difference as:
\begin{align}
    &\calE_R(n) = \frac{\TNORM{\Bvh^{0}}{\Eh}^2 +
      \gamma\Dt\sum_{\ell=0}^{n}\beta^{n+1-\ell}
      \TNORM{\IVh\Es_0^{n-\ell+\theta}}{\HROTv{\Omega}}^2}{\TNORM{\Bvh^{0}}{\Eh}^2}, 
      \\
    &\calE_L(n) = \frac{(\beta)^{n+1}\TNORM{\Bvh^{n+1}}{\Eh}^2 
      +\frac{\gamma\Dt}{2}\sum_{\ell=0}^{n}\beta^{n+1-\ell}\TNORM{\hEsh^{n-\ell+\theta}}{\Vh}^2}{\TNORM{\Bvh^{0}}{\Eh}^2},
      \\
    &\calE(n) = \calE_R(n)-\calE_L(n).
\end{align}
Notice that, by Assumption \eqref{eq:Assumption}, the value of $\beta$ as defined in \eqref{eq:beta-Q:def} is necessarily smaller than $1$, which implies that most of the coefficients in the terms that appear in $\calE$ decay exponentially. Therefore, we can expect that $\calE\to 1$ as $n\to\infty$ unless the growth, in time, of the electric and magnetic fields is fast enough to offset this decay. To illustrate this, we introduce a parameter $C\in\mathbbm{R}$ and the family of solutions
\begin{align}
\Bv^C(x,y,t) &= 
  \begin{pmatrix}
    50e^y-x\sin(xy)+x\cos(xy)\\[0.25em]
    50 e^x+y\sin(xy)+y\cos(xy)
  \end{pmatrix}e^{Ct},\\
  \Es^C(x,y,t) &= C\big( 50(e^x-e^y)-\cos(xy)-\sin(xy) \big) e^{Ct}
\end{align}
and velocity fields $\uv^C=(\us_x^C,\us_y^C)^T$ with
\begin{align}
  \us_x^C(x,y)&=-C\frac{(-x^2-y^2-1)(\sin(xy)+\cos(xy))}{2(50e^x+y\sin(xy)-y\cos(xy))},\\[0.25em]
  \us_y^C(x,y)&=C\frac{(-x^2-y^2-1)(\sin(xy)+\cos(xy))}{2(50e^y-x\sin(xy)+x\cos(xy))}
\end{align}
and define conductivity $\sigma\equiv 1/C$.

\begin{figure}[H]
  \begin{center}
    \begin{tabular}{cc}
      \begin{overpic}[width=.475\textwidth]{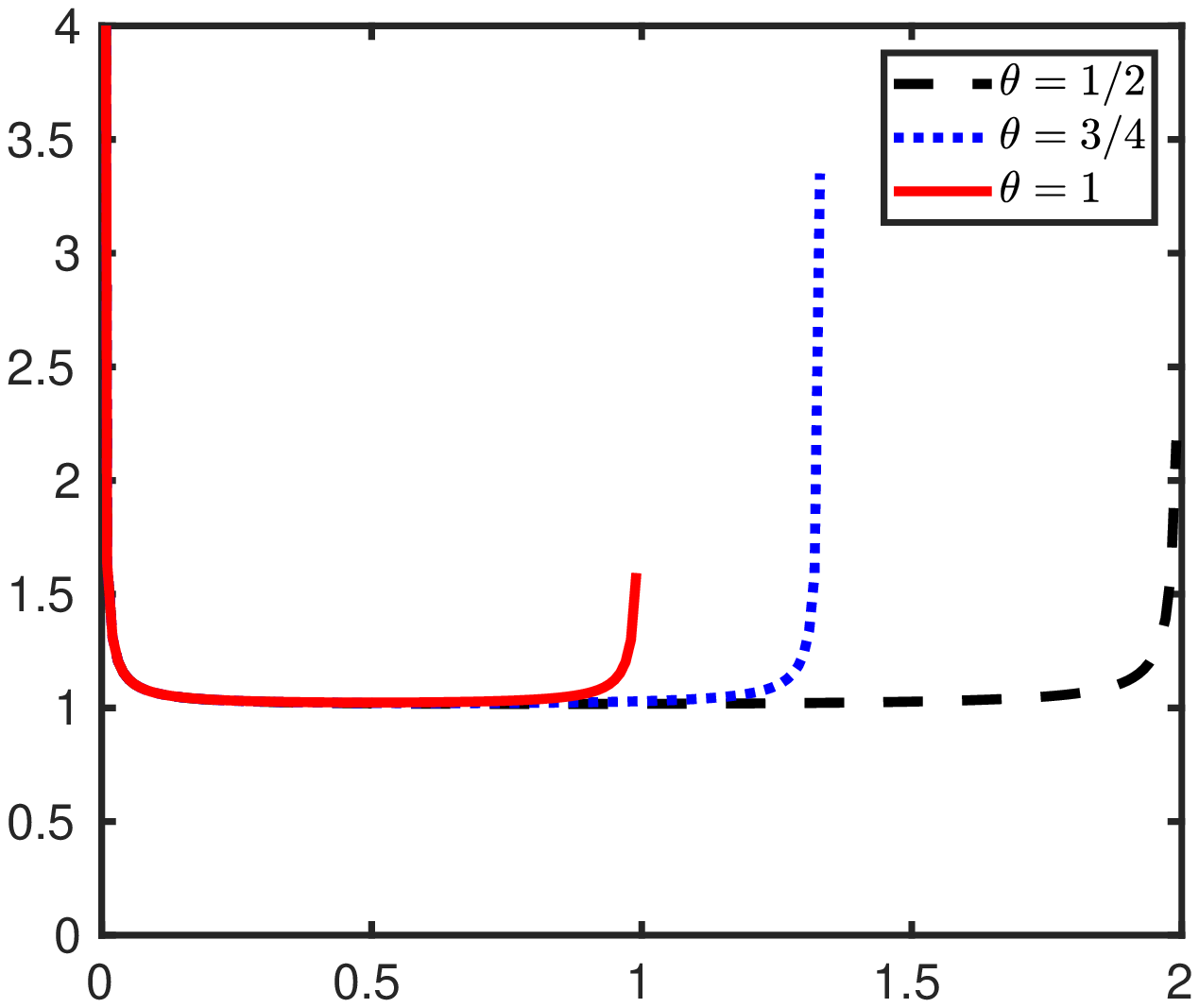}
       \put( 1,35){\begin{sideways}\textbf{$\calE$}\end{sideways}}
        \put(50,0){\textbf{Q}}
      \end{overpic}
      \qquad
    \begin{overpic}[width=.475\textwidth]{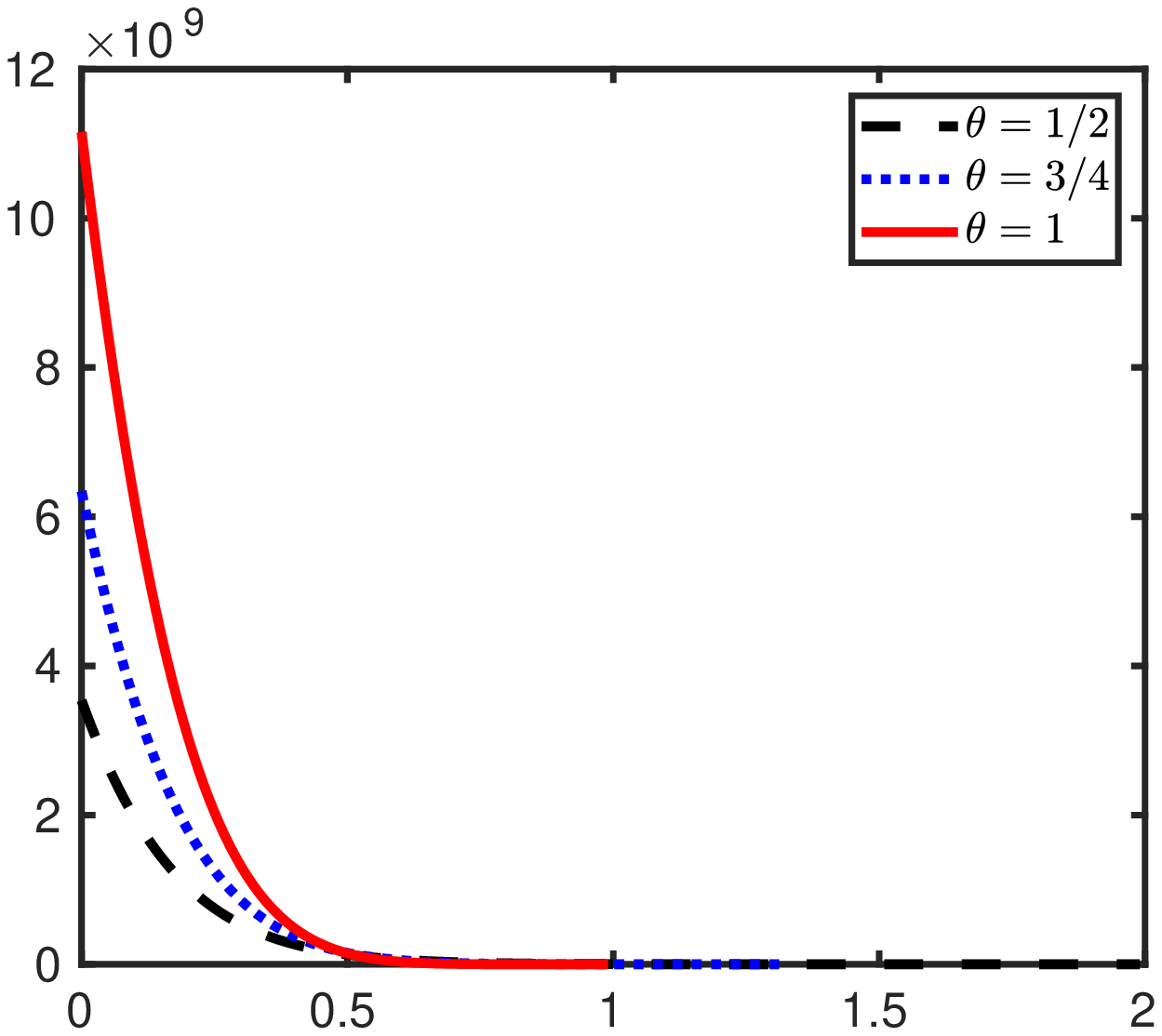}
       \put( 1,35){\begin{sideways}\textbf{$\calE$}\end{sideways}}
        \put(50,0){\textbf{Q}}
      \end{overpic}
     \end{tabular}
  \end{center}
  
  \caption{Plot of Q against the resulting energy estimate at time $T=0.5$. The initial data that yields the plot in the left is that associated with $C = 0.1$ and time step $\Dt=0.001$, whereas the results in the right plot are associated with $C = 5$ and $\Dt=0.21$}
\label{fig:QvsEnergy}
\end{figure}

Note that the Assumption \eqref{eq:Assumption} yields that any choice of $0<Q<\theta^{-1}$, as defined in \eqref{eq:beta-Q:def}, is admissible.
In Figure~\ref{fig:QvsEnergy} we plot the difference between the right and left hand sides of \eqref{eq:DiscEstTwo} normalized by the squared $L^2$-norm of the initial condition on the magnetic field against the value of $Q$ at time $T=0.5$.
The type of mesh or the alternative on the nodal mass matrix do not yield significant difference to the results in this figure.
Thus, we present the results on Voronoi tessalations of the elliptic projector as a representative with mesh size $h=0.0678$.

\begin{figure}[htb]
  \begin{center}
    \begin{tabular}{cc}
      \begin{overpic}[width=.475\textwidth]{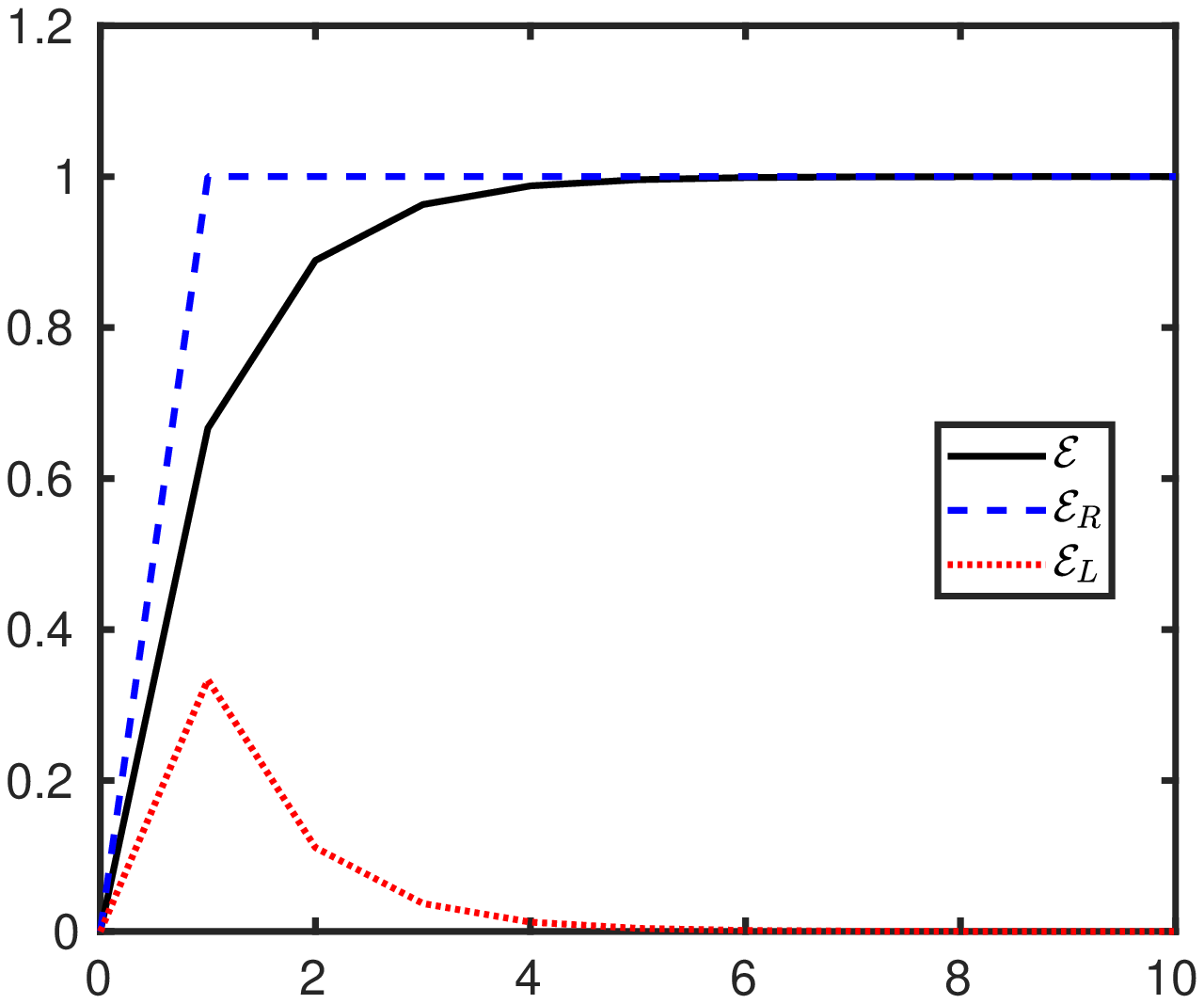}
        \put(20,-2){\textbf{Number of Time Steps}}
      \end{overpic}
      & \qquad
      \begin{overpic}[width=.475\textwidth]{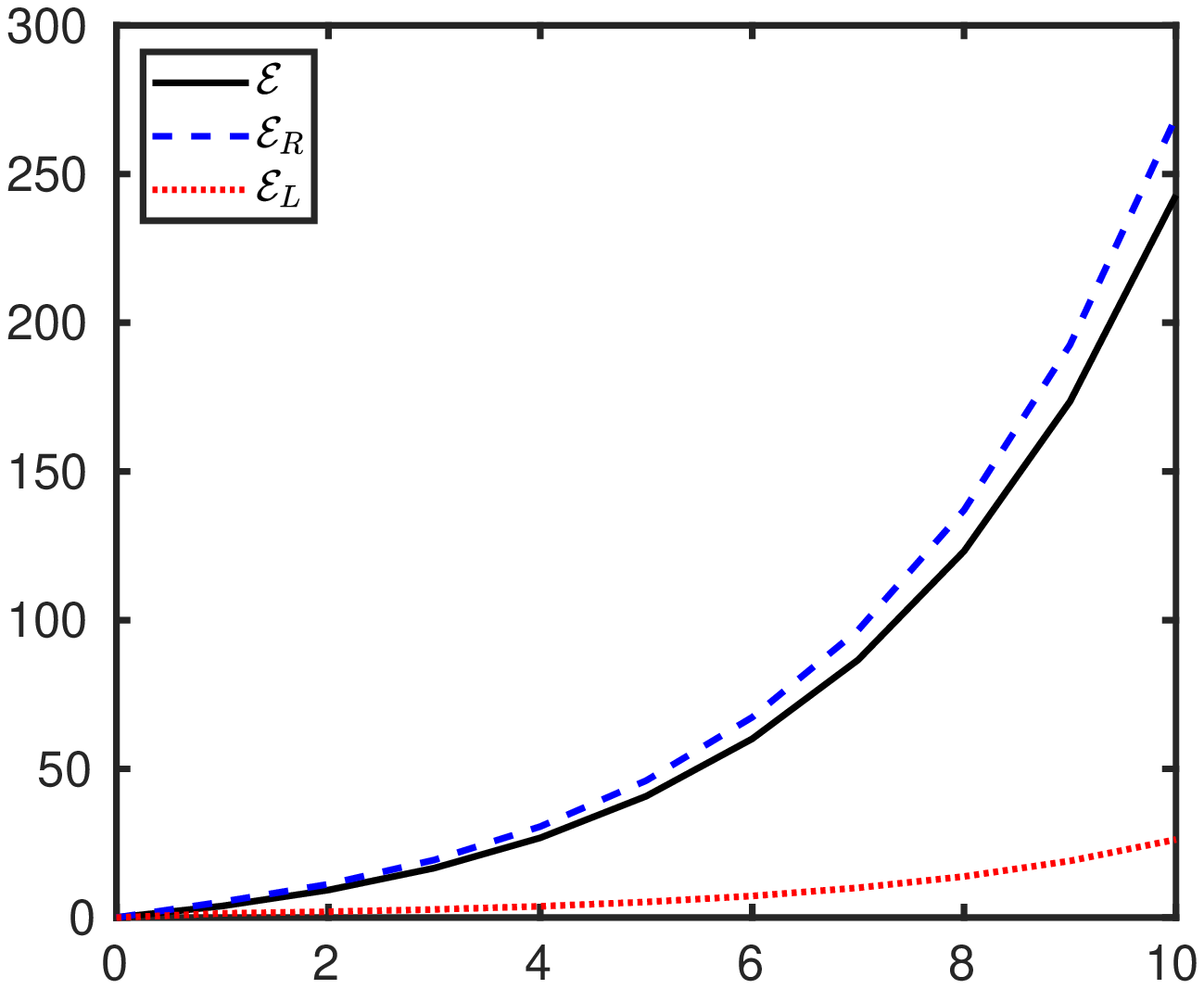}
        \put(20,-2){\textbf{Number Of Time Steps}}
      \end{overpic}
      \\[1.5em]
     \end{tabular}
  \end{center}
  
  \caption{Energy Plots against number of time steps. The initial data that yields the plot in the left is that associated with $C = 0.1$ and time step $\Dt=0.001$, whereas the results on the right plot are associated with $C = 5$ and $\Dt=0.21$. In both cases, $h = 0.0678$.}
\label{fig:EnergyEstimates}\end{figure}

The results of Figure \ref{fig:QvsEnergy} indicate that, in the case that the growth of the solution is relatively small only the values of $Q$ near zero yields $\beta\approx 1$ and the coefficients in $\calE$ will show some exponential growth, if $Q\approx \theta^{-1}$ then the value of $\gamma$ blows up yielding that $\calE$ will be large.
The rest of the values of $Q$ will show convergence towards the norm of the initial conditions on the magnetic field. Since we normalized the error by this value we can expect a flat line of height one.
If, however, the solution grows faster than the decay brought about by the coefficients in $\calE$ then we will see the energy blow up.
Note that the growth in time, at least in our example, of $\calE$ is mainly ruled by terms that look like $\beta^{n} e^{Cn \Dt}$ were $t= n\Dt$, hence a rule of thumb for checking whether the energy will grow or flatten is to check if $\ln{\beta}+C\Dt$ is positive or negative respectively.
This is the reason we picked such a small value for $\Dt$ in the right plot of Figure \ref{fig:QvsEnergy} since large values of $C$ can yield overflow errors.
In Figure \ref{fig:EnergyEstimates} we can clearly see the two different types of behavior that the energy estimates present.
 \subsection{Hartmann Flow}
Consider a square duct of infinite length  containing a conducting fluid.
Assume that this fluid is subjected to a magnetic field that runs along a direction perpendicular to the duct. 
This is the set up for the Hartmann Flow problem which is regarded as a benchmark in MHD.
The behavior of the fluid will depend on the ratio of the Laplace force and the viscous forces, a dimensionless quantity that goes by the name of Hartmann number. 
There is a set of known formulas that describes the solution to this problem, a proof of which can be found in \cite{moreau2013magnetohydrodynamics}.
It is for this reason that researchers use the Hartmann flow problem to test the performance of their simulations, see e.g. \cite{hu2017stable,codina2006stabilized,shadid2016scalable}.

In this section we consider a square computational domain $[-1,1]^2$ as cross section of the aforementioned duct and consider a fluid with conductivity $1$ filling this duct.
The magnetic field is applied in the direction of the $y-$axis.
Consider the case where the viscous forces and Laplace forces are of equal strength, so that the Hartmann number is $1$.
Then, we can expect the fluid to behave in accordance to the solution $\Bv = (B_x,1,0)$, $\uv=(u_x,0,0)$ and $\Ev = (0,0,E_z)$ with 
\begin{equation}
    \begin{aligned}
    &u_x(x,y) = \frac{\cosh{1/2}-\cosh{y}}{2\sinh{1/2}},\\
    &B_x(x,y) = \frac{\sinh{y}-2y\sinh{1/2}}{2\sinh{1/2}},\\
    &E_z(x,y) = \frac{2\sinh{1/2}-\cosh{1/2}}{2\sinh{1/2}}\approx -0.0820.
    \end{aligned}
\end{equation}
Note that the $y-$component of the magnetic field is $1$ by assumption.
Therefore, our main interest in this section is in checking if we can recover approximations to the $x-$component.
To do this we feed the analytical solution for the initial and boundary conditions and evolve the system until $T = 10$ with step size $\Dt =0.005$.
\begin{figure}[H]
  \begin{center}
    \begin{tabular}{cc}
      \begin{overpic}[width=.45\textwidth]{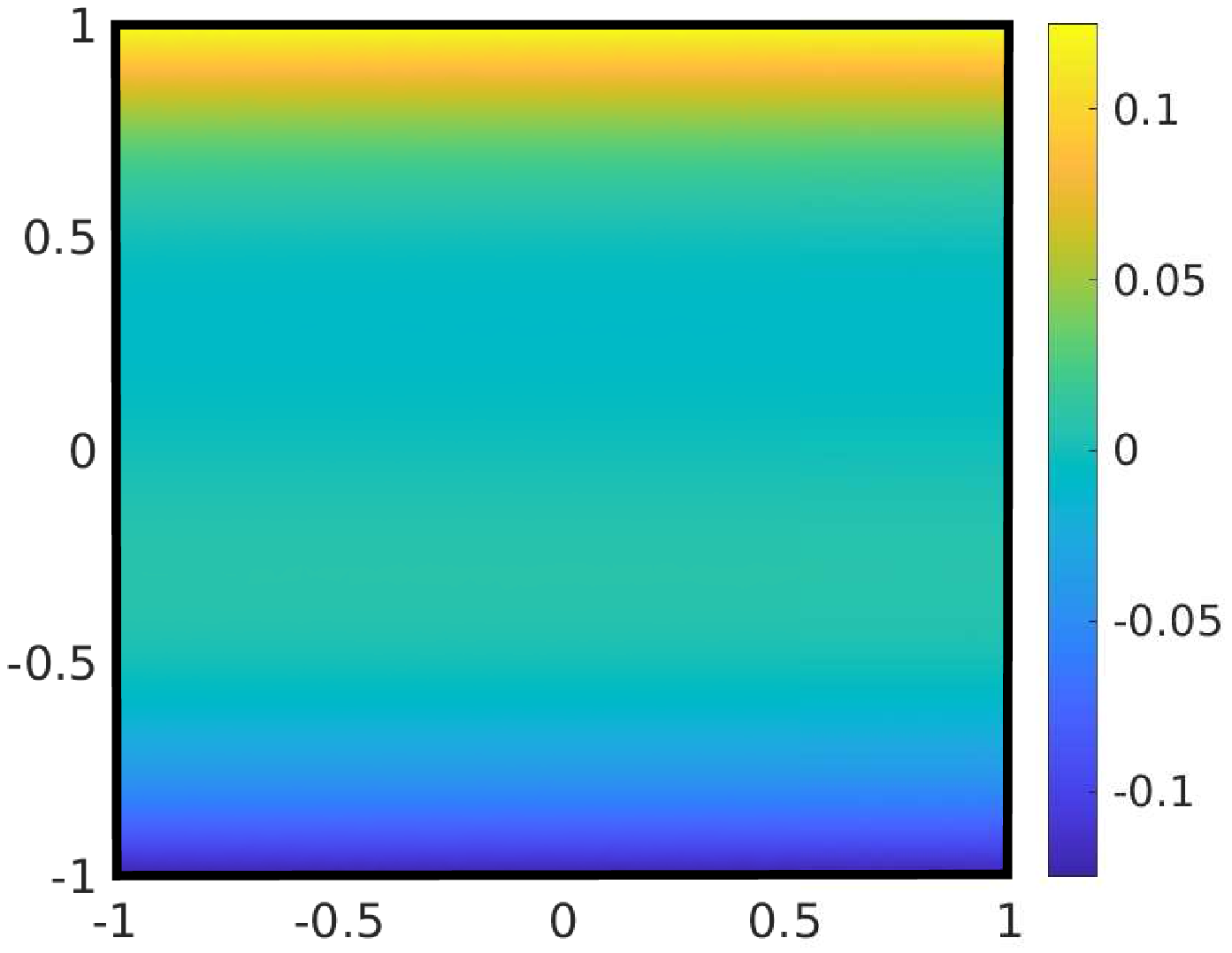}
        \put(35,-2){ { \textbf{x -axis} } }
        \put(-5,35){ 
        \begin{sideways}
        {\textbf{y-axis}} \end{sideways} 
        }
      \end{overpic}
      & \qquad
      \begin{overpic}[width=.45\textwidth]{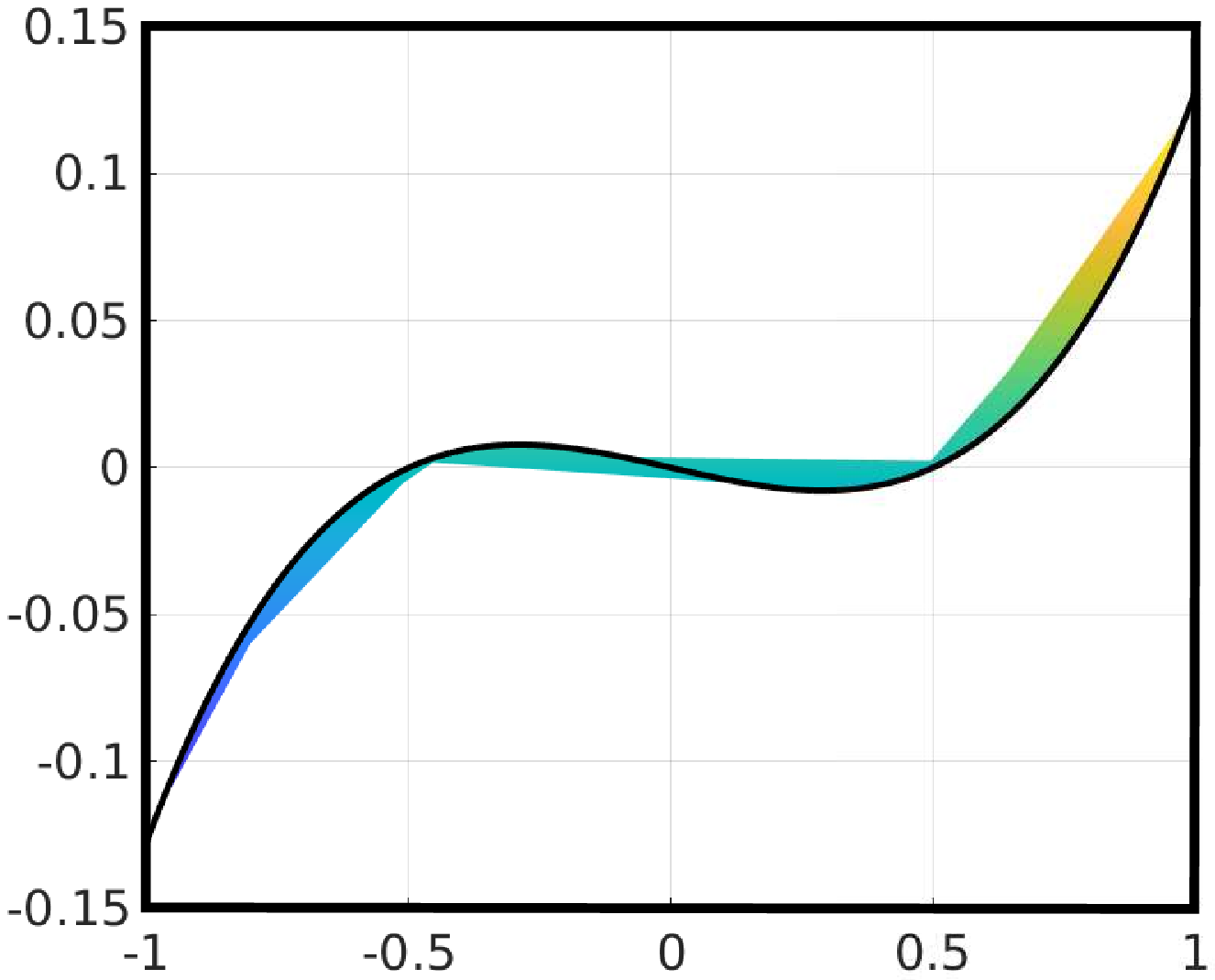}
        \put(45,-2){{\textbf{y-axis}}}
        \put( -5,35){\begin{sideways}{\textbf{z-axis}}\end{sideways}}
      \end{overpic}
     \end{tabular}
  \end{center}
  \caption{Plots of the numerical and analytic solutions for the $x-$component of the magnetic field, computed in a Voronoi tesselation of mesh size $h=0.017$ using the elliptic projector as the alternative to the mass matrix. The plot on the left is of the numerical solution as viewed from above, whereas the plot  on the right shows the numerical solution in a rainbow color bar overlaid with the exact solution in bold black, both are viewed from the side.}
  \label{fig:HartmannSol}
\end{figure}

\begin{figure}[H]
  \begin{center}
    \begin{tabular}{ccc}
      \begin{overpic}[width=.3\textwidth]{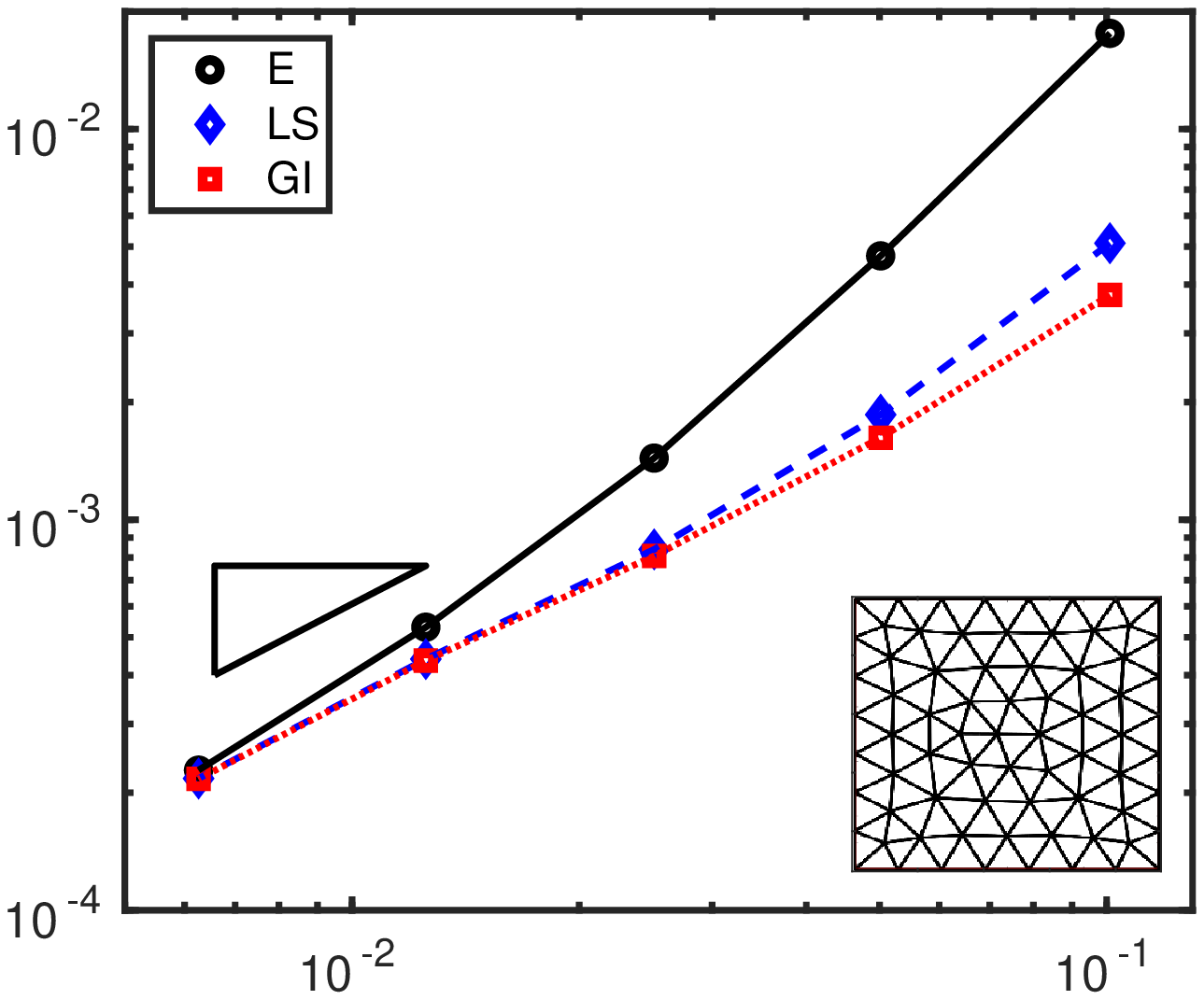}
        \put(30,-2){ { \small\textbf{Mesh Size h} } }
        \put( -5,15){ 
        \begin{sideways}
        {\small\textbf{RelativeError}} \end{sideways} 
        }
        \put(15,28){{\small 1}}
        \put(24,35){{\small 1}}
      \end{overpic}
      & \qquad
      \begin{overpic}[width=.3\textwidth]{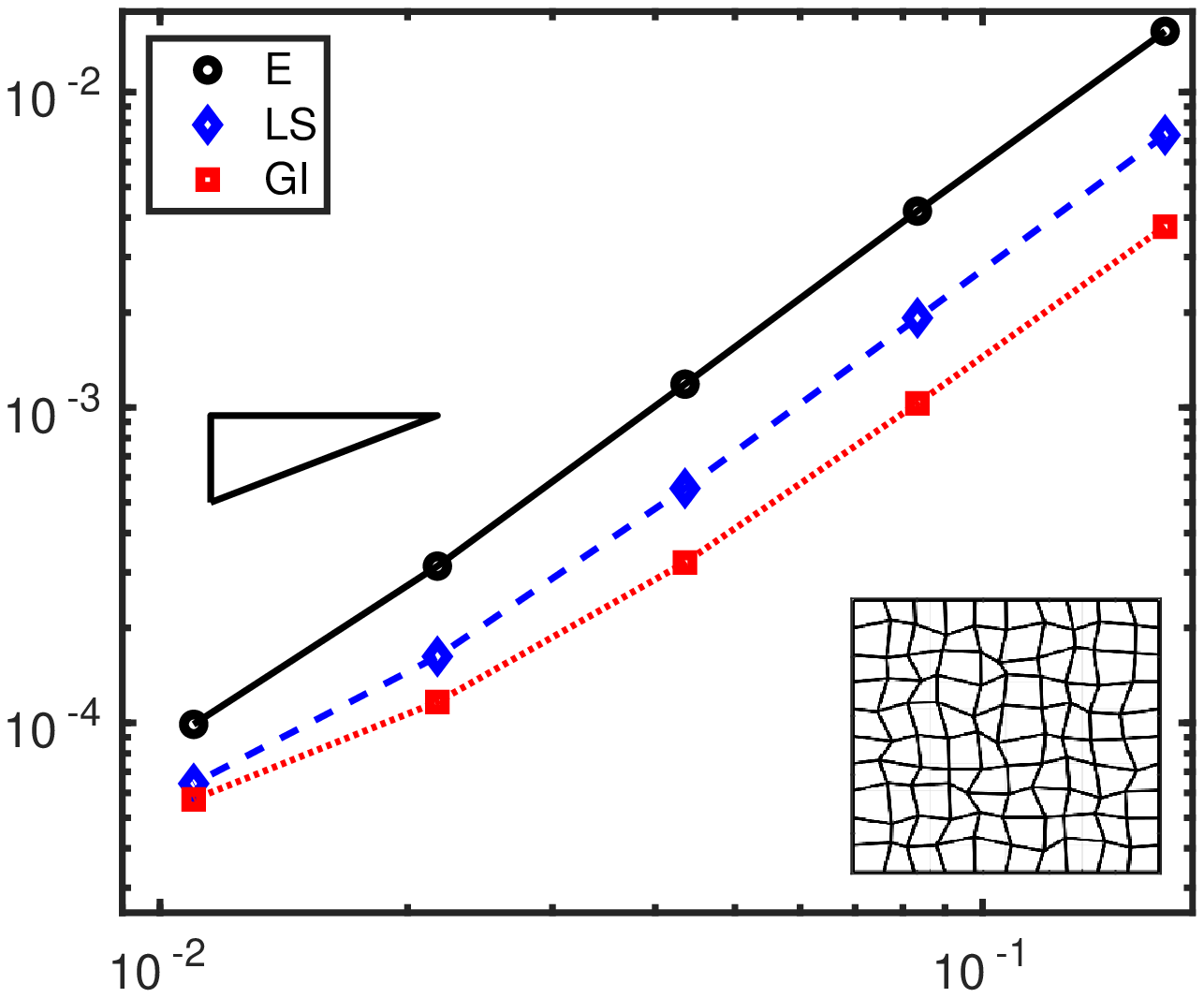}
        \put(30,-2){{\small\textbf{Mesh Size h}}}
        \put( -5,15){\begin{sideways}{\small\textbf{Relative Error}}\end{sideways}}
        \put(15,40){{\small 1}}
        \put(24,46){{\small 1}}
      \end{overpic}
      &\qquad
      \begin{overpic}[width=.3\textwidth]{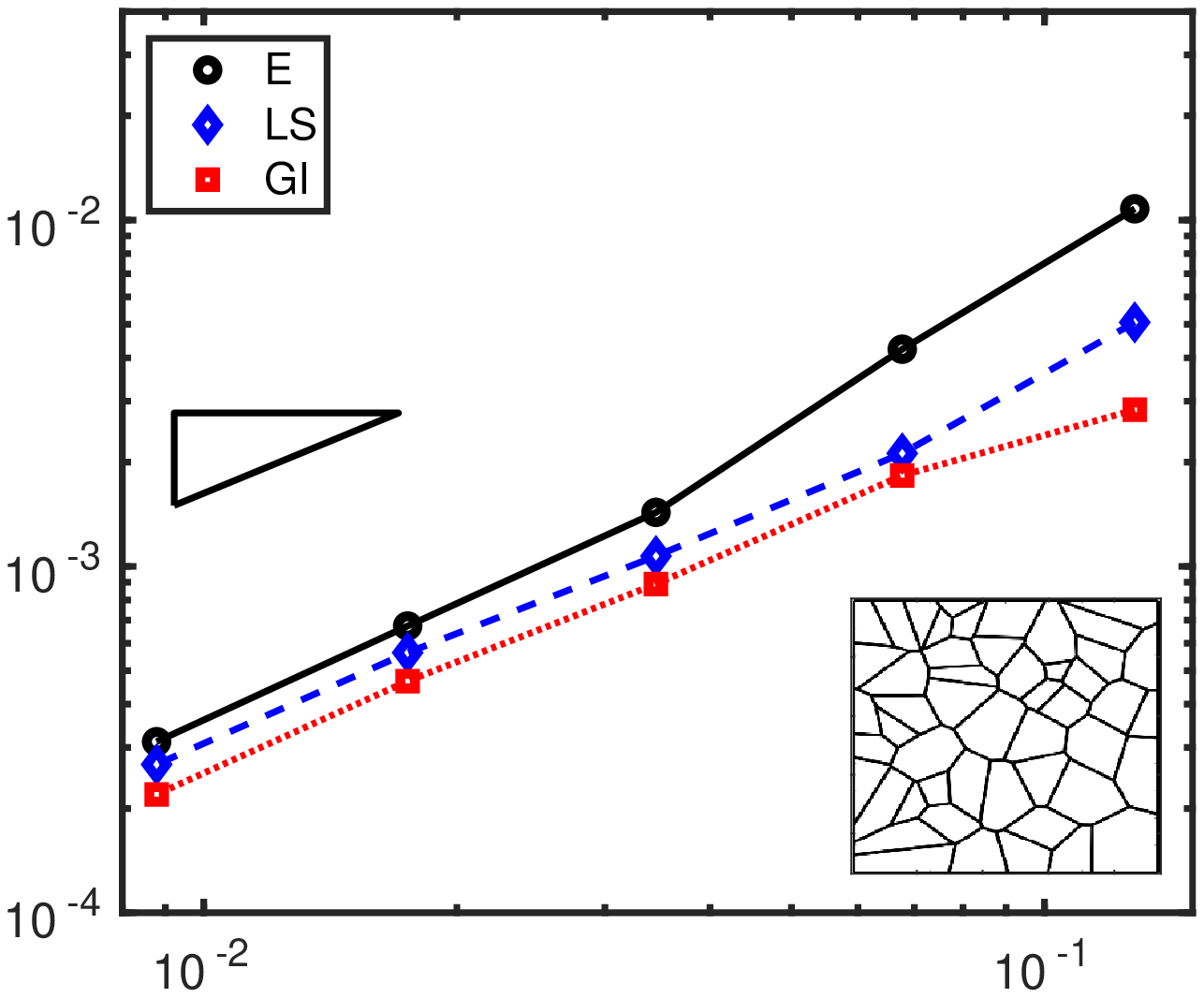}
        \put(30,-2){{\small\textbf{Mesh Size h}}}
        \put(-5,15){\begin{sideways}{\small\textbf{Relative Error}}\end{sideways}}
        \put(13.5,40){{\small 1}}
        \put(22,46){{\small 1}}
      \end{overpic}
     \end{tabular}
  \end{center}
  \caption{Convergence plots for the approximation of the magnetic field on the three different mesh families. The symbols E,LS,GI refers to the three alternatives we have for constructing the nodal mass matrix; the elliptic projector (E), least squares projector (LS) and the Galerkin interpolator (GI) respectively.}
\label{fig:HartamnnConvergencePlots}
\end{figure}

The results, to the naked, eye are satisfactory, Figure~\ref{fig:HartmannSol} gives evidence of this fact.
We further conducted a convergence test that verifies that every alternative to the mass matrix yields a close approximation and provides additional evidence that rate of convergence of the magnetic field is linear, these results are in Figure~\ref{fig:HartamnnConvergencePlots}.
\section{Conclusions}\label{sec:conclusions}
We developed a virtual element method for the Maxwell system of equations
\eqref{eq:IntroSystem} that model the evolution of the electric and
magnetic fields of a magnetized fluid whose flow is prescribed.
It is well documented that, in order to accurately describe the
physics of resistive MHD, it is imperative for the numerical
approximation of the magnetic flux field  to remain divergence free.
This feature is explicitly addressed in this work and
Theorem~\ref{theorem:DivFreeScheme} rigorously proves that the virtual
element scheme in \eqref{eq:VarForEq4} naturally satisfies this
requirement.
The numerical tests in Section~\ref{Sec:NumericalExperiments}
demonstrate that practical implementations of this VEM will satisfy the
divergence free condition on the magnetic flux field.
Moreover, Theorem~\ref{theorem:wellposedness:problem:A} states that
the VEM is wellposed, i.e., that the virtual element approximation
exists and is unique.
We also proved that the VEM is stable through suitable energy
estimates as stated in Theorems
\ref{theorem:ContEnergyEst} and~\ref{theorem:DiscreteEnergyEstimate}.
These estimates were explored numerically in
Section~\ref{Sec:NumericalExperiments}.

\medskip
The performance of the method was investigated experimentally, and a
set of tests using a manufactured solution and the Hartmann flow
problem provide evidence of a quadratic convergence rate for the
approximation of electric field and a linear convergence rate for the
approximation of the magnetic field.

\medskip
Future work will focus on the design and implementation of a VEM for
\eqref{eq:IntroSystem} in three dimensions and combining such
formulation with a flow equation describing the conservation of
momentum and mass.
We will also introduce a non-linear term into Ohm's law to describe
physical effects related to Hall currents.
Such a term is proportional to $\Jv\times\Bv$ with
$\Jv=\nabla\times\Bv$ as dictated by Ampere's law.
We are also planning to develop a higher order accurate VEM by
increasing the order of the local polynomial subspaces, while
preserving the commuting de Rham diagram and the free-divergence
condition for the magnetic flux field.

\medskip
\noindent
\textbf{Acknowledgements.}\\
S. Naranjo Alvarez's work was supported by the National Science
Foundation (NSF) grant \#1545188, ``\emph{NRT-DESE: Risk and
  uncertainty quantification in marine science and policy}'', which provided a one year fellowship and internship support at Los Alamos National Laboratory. In addition, S. Naranjo Alvarez received graduate research funding from V. A. Bokil's  DMS grant \#1720116, an INTERN supplemental award to Professor Bokil's DMS grant \# 1720116 for a second internship at Los Alamos National Laboratory, and teaching support from the Department of Mathematics at Oregon State University. 
  
  Professor V. A. Bokil was partially supported by NSF funding from the DMS grant \# 1720116. 
Dr. V. Gyrya and Dr. G. Manzini were supported by 
the Laboratory Directed Research and Development - Exploratory
Research (LDRD-ER) Program of Los Alamos National Laboratory under
project number 20180428ER.

The authors would like to thank Dr. K. Lipnikov, T-5 Group, Theoretical Division,
Los Alamos National Laboratory, for his advice during the writing of
this article.
Los Alamos National Laboratory is operated by Triad National Security,
LLC, for the National Nuclear Security Administration of
U.S. Department of Energy (Contract No. 89233218CNA000001).

\bibliographystyle{abbrv}
\bibliography{./BIB/pub-mfd,./BIB/pub-varie,./BIB/pub-vem}

\appendix
\section{Proof of \textbf{(V3)} for $\PinP$}
\label{app:sec:proof:V3:PinP}

\newcommand{\xvb}{\overline{\xv}}
\newcommand{\bvb}{\overline{\bv}}
\newcommand{\asb}{\overline{\as}}

We write the elliptic projection of $\vsh\in\Vh(\P)$ as the linear polynomial
$\PinP\vs=\asb+\bvb\cdot(\xv-\xvb)$, where 
\begin{align*}
  \xvb = \frac{1}{\ABS{\partial\P}}\int_{\partial\P}\xv\dV,\quad
  \asb = \frac{1}{\ABS{\partial\P}}\int_{\partial\P}\vsh\dS,\quad
  \bvb = \frac{1}{\mP}\int_{\P}\nabla\vsh\dV,
\end{align*}
and $\ABS{\partial\P}$ is the perimeter of $\P$.
A straightforward calculation yields
\begin{align*}
  \NORM{\PinP\vs}{0,\P}^2 
  &=\int_{\P}\ABS{\asb+\bvb\cdot(\xv-\xvb)}^2\dV
  \leq 2\ABS{\asb}^2\mP + 2\int_{\P}\ABS{\bvb}^2\ABS{\xv-\xvb}^2\dV
  \nonumber\\[0.5em]
  &\leq 2\ABS{\asb}^2\mP + 2\ABS{\bvb}^2\,\int_{\P}\ABS{\xv-\xvb}^2\dV
  \leq 2\ABS{\asb}^2\mP + 2\Cs\ABS{\bvb}^2\mP\hP^2,
\end{align*}
where $\Cs$ is a ``geometric'' constant that may depend on the shape
of $\P$ but does not scale with $\hP$ since
\begin{align*}
  \frac{1}{\mP}\int_{\P}\ABS{\xv-\xvb}^2\dV \simeq \hP^2\mathcal{O}(1).
\end{align*}
Then, first using Jensen's inequality, and, then, Agmon's inequality
yields
\begin{align*}
  \mP\ABS{\asb}^2 
  &= \frac{\mP}{\ABS{\partial\P}^2} \vABS{ \int_{\partial\P}\vsh\dS }^2 
  \leq \frac{\mP}{\ABS{\partial\P}^2}\,\ABS{\partial\P}\int_{\partial\P}\ABS{\vsh}^2\dS 
  \leq \frac{\mP}{\ABS{\partial\P}}\int_{\partial\P}\ABS{\vsh}^2\dS 
  \nonumber\\[0.5em]  
  &\leq \frac{\mP}{\ABS{\partial\P}}\,\Cs^{A} \big( \hP\SNORM{\vsh}{1,\P}^2
  + \hP^{-1}\NORM{\vsh}{0,\P}^2 \big)
  \leq \Cs\Cs^{A}\left(
    \hP^2\SNORM{\vsh}{1,\P}^2 + \NORM{\vsh}{0,\P}^2
  \right),
\end{align*}
where $\Cs^{A}$ is the constant of Agmon's inequality.
In the above inequality we used the fact
$\mP\hP\slash{{\ABS{\partial\P}}}=\mathcal{O}\big(\hP^2\big)$ and
$\mP\hP^{-1}\slash{{\ABS{\partial\P}}}=\mathcal{O}(1)$ since
$\hP\slash{\ABS{\partial\P}}$ and $\mP\slash{\hP^2}$ are uniformly bounded
quantities in view of Assumption~\textbf{(M2)}.
This assumption also implies that the real positive constant $\Cs$ may
only depend on $\rho$, and, from Agmon's inequality, on the number of
polygonal edges, this latter also being uniformly bounded.
Using again Jensen's inequality yields
\begin{align*}
  \mP\hP^2\ABS{\bvb}^2
  =    \frac{\mP\hP^2}{\mP^2}\vABS{ \int_{\P}\nabla\vsh\dV }^2
  =    \frac{\mP\hP^2}{\mP^2}\,\mP\int_{\P}\ABS{\nabla\vsh}^2\dV
  \leq \hP^2\int_{\P}\ABS{\nabla\vsh}^2\dV
  =    \hP^2\SNORM{\vsh}{1,\P}^2.
\end{align*}
Finally, we collect the estimates for $\mP\ABS{\asb}^2$ and
$\mP\hP^2\ABS{\bvb}^2$, and apply the inverse inequality
\begin{align}
\SNORM{\vsh}{1,\P}\leq\Cs^{I}\hP^{-1}\NORM{\vsh}{0,\P},
\end{align}
which follows from a scaling argument,
see~\cite[Chapter~2]{Boffi-Brezzi-Fortin:2013} and the recent work of
Ref.~\cite{Chen-Huang:2017}, and whose constant $\Cs^{I}$ is
independent of $\hP$ to obtain:
\begin{align*}
  \NORM{\PinP\vs}{0,\P}^2
  \leq \Cs\left( \hP^2\SNORM{\vsh}{1,\P}^2 + \NORM{\vsh}{0,\P}^2 \right)
  \leq \Cs\left( \Cs^{I}\hP^2\hP^{-2}\NORM{\vsh}{0,\P}^2 + \NORM{\vsh}{0,\P}^2 \right)
  \leq \Cs\NORM{\vsh}{0,\P}^2.
\end{align*}
Tracing back the constants introduced in the various inequalities, we
find that the final constant $\Cs$ may depend on $\Cs^{A}$, $\Cs^{I}$,
$\rho$, but is independent of $\hP$, and is obviously the same for all
$\vsh\in\Vh(\P)$.
This argument provides the desired upper bound on operator $\PinP$.

\section{Proof of \textbf{(V3)} for $\PiVhPLS$}
\label{app:sec:proof:V3:PiVhPLS}

The Least Squares reconstruction operator applied to $\vsh\in\Vh(\P)$
provides the linear polynomial~\eqref{eq:least-squares-polynomials},
which we conveniently rewrite here:
\begin{align}
  \PiVhPLS\vsh(\xs,\ys) 
  = \as+\bs\frac{\xs-\xs_{\P}}{\hP} + \cs\frac{\ys-\ys_{\P}}{\hP}.
  \label{app:eq:least-squares-polynomials}
\end{align}
The three coefficients $\as$, $\bs$ and $\cs$ are determined by
imposing the conditions in~\eqref{eq:least-squares-conditions}.
Let $\bzeta=(\as,\bs,\cs)^T$ denote the vector collecting the three
unknowns in~\eqref{app:eq:least-squares-polynomials}.
Let $\xv_{\V_i}=(\xs_{\V_i},\ys_{\V_i})^T$ denote the coordinate
vector of the $i$-th vertex $\V_i$ for $i=1,2,\ldots,\NV$, $\NV$ being
the number of vertices of $\P$, and
$\etav=\big(\vsh(\xv_{\V_1}),\vsh(\xv_{\V_2}),\ldots,\vsh(\xv_{\V_{\NV}})\big)^T$
the vector collecting the nodal degrees of freedom of $\vsh$.
Using this notation, we rewrite the linear
system~\eqref{eq:least-squares-conditions} in the more compact form
$\matA\bzeta=\etav$, where the matrix of the system coefficients is
given by:
\begin{align*}
  \matA =
  \left[
    \begin{array}{cc}
      1      & \quad\dfrac{(\xv_{\V_1}   -\xv_{\P})^T}{\hP}  \\[0.5em]
      1      & \quad\dfrac{(\xv_{\V_2}   -\xv_{\P})^T}{\hP}  \\[0.5em]
      \vdots & \vdots                                      \\[0.5em]
      1      & \quad\dfrac{(\xv_{\V_{\NV}}-\xv_{\P})^T}{\hP}
    \end{array}
  \right]
  =
  \left[
    \begin{array}{ccc}
      1      & \quad\dfrac{\xs_{\V_1}   -\xs_{\P}}{\hP}  & \quad\dfrac{\ys_{\V_1}    -\ys_{\P}}{\hP}\\[0.5em]
      1      & \quad\dfrac{\xs_{\V_2}   -\xs_{\P}}{\hP}  & \quad\dfrac{\ys_{\V_2}    -\ys_{\P}}{\hP}\\[0.5em]
      \vdots & \vdots                                  & \vdots                           \\[0.5em]
      1      & \quad\dfrac{\xs_{\V_{\NV}}-\xs_{\P}}{\hP} & \quad\dfrac{\ys_{\V_{\NV}}-\ys_{\P}}{\hP}
    \end{array}
  \right].
\end{align*}
The coefficients of the least squares solution are given by solving
the normal equations, i.e., $\bzeta=(\matA^T\matA)^{-1}\matA^T\etav$.

\medskip
Now, we introduce the discrete norm
\begin{align}
  \TNORM{\vsh}{\P}^2 
  = \mP\sum_{i=1}^{\NV}\ABS{\vsh(\xv_{\V_i})}^2
  = \mP\abs{\etav}^2
  \label{app:eq:discrete:norm}
\end{align}
and we observe that 
\begin{align}
  \TNORM{\PiVhPLS\vsh}{\P}^2
  =    \mP\sum_{i=1}^{\NV}\abs{\PiVhPLS\vsh(\xv_{\V_i})}^2
  =    \mP\sum_{i=1}^{\NV}\Abs{ \as+\bs\frac{\xs_{\V_i}-\xs_{\P}}{\hP} + \cs\frac{\ys_{\V_i}-\ys_{\P}}{\hP} }^2
  =    \mP\abs{\matA\bzeta}^2.
  \label{app:eq:discrete:norm:10}
\end{align}
The norm defined in~\eqref{app:eq:discrete:norm} is spectrally
equivalent to the $\LTWO$ norm, so that there exist two strictly
positive constant $\nu_*$ and $\nu^*$ such that
\begin{align}
  \nu_*\NORM{\vsh}{0,\P}\leq\TNORM{\vsh}{\P}\leq\nu^*\NORM{\vsh}{0,\P}
  \quad\forall\vsh\in\Vh(\P).
  \label{app:eq:norm:equivalence}  
\end{align}
The two norms $\NORM{\vsh}{0,\P}$ and $\TNORM{\vsh}{\P}$, because of
the explicit dependence of the latter on $\mP$, have the same scaling
with respect to $\hP$.
Therefore, the two constants $\nu_*$ and $\nu^*$ may depend on the
geometric shape of $\P$ but must be independent of $\hP$.
  
\medskip
A straightforward calculation starting
from the left inequality of~\eqref{app:eq:norm:equivalence} yields
\begin{align*}
  \begin{array}{rll}
    \NORM{\PiVhPLS\vs}{0,\P}
    &\leq (\nu_*)^{-1}\TNORM{\PiVhPLS\vsh}{\P}                                             & \hspace{1cm}\mbox{\big[use~\eqref{app:eq:discrete:norm:10}\big]}                              \\[0.5em]
    &=    (\nu_*)^{-1}\mP^{\frac12}\abs{\matA\bzeta}                                        & \hspace{1cm}\mbox{\big[substitute~$\bzeta=\big(\matA^T\matA\big)^{-1}\matA^T\etav$\big]}       \\[0.5em]
    &=    (\nu_*)^{-1}\mP^{\frac12}\abs{\matA\big(\matA^T\matA\big)^{-1}\matA^T\etav}        & \hspace{1cm}\mbox{\big[use the continuity of $\matA\big(\matA^T\matA\big)^{-1}\matA^T$\big]}   \\[0.5em]
    &\leq (\nu_*)^{-1}\mP^{\frac12}\abs{\matA\big(\matA^T\matA\big)^{-1}\matA^T}\ABS{\etav}  & \hspace{1cm}\mbox{\big[note that $\matA\big(\matA^T\matA\big)^{-1}\matA^T$ is a projector\big]}\\[0.5em]
    &\leq (\nu_*)^{-1}\mP^{\frac12}\abs{\etav}                                              & \hspace{1cm}\mbox{\big[use~\eqref{app:eq:discrete:norm}\big]}                                 \\[0.5em]
    &=    (\nu_*)^{-1}\TNORM{\vsh}{\P}                                                     & \hspace{1cm}\mbox{\big[use the right inequality of~\eqref{app:eq:norm:equivalence}\big]}      \\[0.5em]
    &\leq \frac{\nu^*}{\nu_*}  \NORM{\vsh}{0,\P}.
  \end{array}
\end{align*}
In the chain of inequalities above, we used the fact that
$\matA\big(\matA^T\matA\big)^{-1}\matA^T$ is the orthogonal projection operator
with respect to the Euclidean inner product onto the span of the
columns of matrix $\matA$.
This projection operator scales like $\mathcal{O}(1)$ with respect to
$\hP$ by definition of matrix $\matA$, and its eigenvalues must be $0$
and $1$.
As a consequence, it is continuous and such that $\ABS{\matA\big(\matA^T\matA)^{-1}\matA^T}=1$.
Finally, we note that in this case the constant $\Cs$ that appears in the assertion of proposition $\textbf{(V3)}$ is equal to $(\nu^*\slash{\nu_*})$, and is independent of $\hh$ as already pointed out above.

\end{document}